# SMARANDACHE SEMIGROUPS

$$p_1 : \begin{cases} x_1 \to x_1 \\ x_2 \to x_3 \\ x_3 \to x_2 \end{cases}$$

$\sigma(1) = 2$      $\mu(1) = 3$
$\sigma(2) = 3$      $\mu(2) = 1$
$\sigma(3) = 1$      $\mu(3) = 2.$

$$p \left| \left( o(G) - \sum_{N(a) \neq G} \frac{o(G)}{o(N(a))} \right) \right. = o(Z(G)).$$

W. B. Vasantha Kandasamy

# Smarandache Semigroups

$$A = \left\{ \begin{pmatrix} 1 & 0 & \ldots & 0 & 0 \\ 0 & 1 & \ldots & 0 & 0 \\ \vdots & \vdots & \ldots & \vdots & \vdots \\ 0 & 0 & \ldots & 0 & 0 \\ 0 & 0 & \ldots & 0 & 1 \end{pmatrix} = I_{n \times n}, \begin{pmatrix} 0 & 0 & \ldots & 0 & 1 \\ 0 & 0 & \ldots & 1 & 0 \\ \vdots & \vdots & \ldots & \vdots & \vdots \\ 0 & 1 & \ldots & 0 & 0 \\ 1 & 0 & \ldots & 0 & 0 \end{pmatrix} \right\}$$

**W.B.Vasantha Kandasamy**
Department of Mathematics
Indian Institute of Technology Madras
Chennai – 600 036 India

e-mail: vasantha@iitm.ac.in
web: http://mat.iitm.ac.in/~wbv

**Definition:**

**Generally, in any human field, a *Smarandache Structure* on a set A means a weak structure W on A such that there exists a proper subset B contained in A which is embedded with a stronger structure S.**

**These types of structures occur in our everyday life, that's why we study them in this book.**

**Thus, as a particular case:**

**A *Smarandache Semigroup* is a semigroup A which has a proper subset B contained in A that is a group (with respect to the same binary operation on A).**

# CONTENTS

# PREFACE

The main motivation and desire for writing this book, is the direct appreciation and attraction towards the Smarandache notions in general and Smarandache algebraic structures in particular. The Smarandache semigroups exhibit properties of both a group and a semigroup simultaneously. This book is a piece of work on Smarandache semigroups and assumes the reader to have a good background on group theory; we give some recollection about groups and some of its properties just for quick reference.

Since most of the properties and theorems given regarding the Smarandache semigroups are new and cannot be found in existing literature the author has taken utmost efforts to see that the concepts are completely understood by illustrating with examples and a great number of problems. Solutions to all the problems need extraordinary effort.

The book is organized in the following way: It has seven chapters. The first chapter on preliminaries gives some important notions and concepts, which are used in this book. Chapters 2 and 3 gives most of the basic concepts on group theory and results in group theory which have been used in this text to study Smarandache notions in groups or Smarandache semigroups. This text does not in any way claim completeness in giving the properties of groups.

Chapter 4 starts with the definition of the Smarandache semigroup and gives some interesting properties of Smarandache semigroups. This chapter is made easy for comprehension by several examples. The problems are a must for the researchers to solve, for they alone will give them the complete conceptual understanding of the Smarandache semigroup.

In chapter five we make use of the newly defined and special types of Smarandache semigroups in proving or disproving the classical theorems or analogs of the classical theorems. This chapter is also substantiated with examples and several problems are given. The sixth chapter is a mixture of both Smarandache notions on groups and the study of properties of Smarandache semigroups.

The final chapter, a special attraction to researchers and algebraists is a list of open research problems. Most of the proposed problems are not very easy to solve, but certainly, this feature will attract not only research students but also their research guides to take up research on Smarandache notions. Smarandache notions are revolutionary because group theory does not make one think of the unthinkable that can naturally occur but Smarandache semigroup explicitly and concretely expresses the possibilities of such occurrences like the validity of Lagrange's theorem, Cauchy's theorem and Sylow's theorem. Since Smarandache semigroups are the overlap of two structures, we are able to see how the mixture of a group and a semigroup behaves.

I deeply acknowledge Dr. Minh Perez for his constant encouragement and support which made this book possible.

**CHAPTER ONE**

# PRELIMINARY NOTIONS

One of the essential and outstanding features of the twenty-first century mathematics has been not only the recognition of the power of abstract approach but also its simplicity by way of illustrative examples and counter examples. Here the Smarandache notions in groups and the concept of Smarandache semigroups which are a class of very innovative and conceptually a creative structure, has been introduced in the context of groups and a complete possible study has been taken in this book. Thus, the main purpose of this book is to make both researcher and an algebraist to know and enjoy the Smarandache analog concept for groups. It is pertinent to mention that Smarandache notions on all algebraic and mathematical structures are interesting to the world of mathematicians and researchers, so at this juncture we felt it would be appropriate to study the Smarandache semigroups and introduce some Smarandache notions in groups.

The introduction of Smarandache concepts in groups and Smarandache semigroups in a way makes a mathematician wonder when some of the classical theorems like Lagrange's become untrue but at the same time enables for a lucid extension of the Cayley's theorem. This book deals only with the algebraic structure of groups in the context of Smarandache structures. To make this book self-contained much effort is taken to see that the chapters two and three give most of the relevant concepts in group theory for which we have introduced the Smarandache notions. In these two chapters we have restrained ourselves by not giving any problems, the problems, which are included, are those results that are essential for our study of Smarandache notions in groups and Smarandache semigroups. Chapters 4, 5 and 6 deals with Smarandache semigroups and Smarandache notions in groups are crowded with examples, counter-examples and problems. The prominent feature of the book is - all abstract concepts are illustrated by examples.

As to the best of the author's knowledge concepts and results about Smarandache semigroups and Smarandache notions in groups are very meager or absent except for the definition. Here all pains are taken to introduce in a best sequential way for the reader to appreciate and contribute to the subject of Smarandache semigroups. The last chapter is completely devoted to research problems some of them are really very difficult for these problems may attract a research student and an algebraist and force them to contribute something to the world of Smarandache notions in groups and Smarandache semigroups. In Chapter one, we introduce some basic notation, Binary relations, mappings and the concept of semigroup and Smarandache semigroup.

## 1.1 Binary Relation

Let A be any non-empty set. We consider the Cartesian product of a set A with itself; $A \times A$. Note that if the set A is a finite set having n elements, then the set $A \times A$ is also a finite set, but has $n^2$ elements. The set of elements (a, a) in $A \times A$ is called the diagonal of A × A. A subset S of A × A is said to define an equivalence relation on A if

    (a, a) ∈ S for all a ∈ A
    (a, b) ∈ S implies (b, a) ∈ S
    (a, b) ∈ S and (b, c) ∈ S implies that (a, c) ∈ S.

Instead of speaking about subsets of A × A we can now speak about a binary relation (one between two elements of A) on A itself, by defining b to be related to a if (a, b) ∈ S. The properties 1, 2, 3 of the subset S immediately translate into the properties 1, 2, 3 of the following definition.

**DEFINITION:** *The binary relation ~ on A is said to be an equivalence relation on A if for all a, b, c in A*

    i.    *a ~ a*
    ii.   *a ~ b implies b ~ a*
    iii.  *a ~ b and b ~ c implies a ~ c*

The first of these properties is called reflexivity, the second, symmetry and the third transitivity. The concept of an equivalence relation is an extremely important one and plays a central role in all mathematics.

**DEFINITION:** *If A is a set and if ~ is an equivalence relation on A, then the equivalence class of a ∈ A is the set {x ∈ A/ a ~ x}. We write this set as cl(a) or [a].*

**THEOREM 1.1.1:** *The distinct equivalence classes of an equivalence relation on A provide us with a decomposition of A as a union of mutually disjoint subsets. Conversely, given a decomposition of A as a union of mutually disjoint, nonempty subsets, we can define an equivalence relation on A for which these subsets are the distinct equivalence classes.*

*Proof:* Let the equivalence relation on A be denoted by ~. We first note that since for any a ∈ A, a ~ a, a must be in cl(a), whence the union of the cl(a)'s is all of A. We now assert that given two equivalence classes they are either equal or disjoint. For, suppose that cl(a) and cl(b) are not disjoint: then there is an element x ∈ cl(a) ∩ cl(b). Since x ∈ cl(a), a ~ x; since x ∈ cl(b), b ~ x, whence by the symmetry of the relation, x ~ b. However, a ~ x and x ~ b by the transitivity of the relation forces a ~ b. Suppose, now that y ∈ cl(b); thus b ~ y.

However, from a ~ b and b ~ y, we deduce that a ~ y, that is, that y ∈ cl(a). Therefore, every element in cl(b) is in cl(a), which proves that cl(b) ⊂ cl(a). The argument is clearly symmetric, whence we conclude that cl(a) ⊂ cl (b). The argument opposite containing relations imply that cl(a) = cl(b). We have thus shown that the distinct cl(a)'s are mutually disjoint and that their union is A. This proves the first half of the theorem. Now for the other half! Suppose that $A = \bigcup A_\alpha$ where the $A_\alpha$ are mutually disjoint, nonempty sets (α is in some index set T). How shall we use them to define an equivalence relation? The way is clear; given an element, a in A it is in exactly one $A_\alpha$. We define for a, b ∈ A, a ~ b if a and b are in the same $A_\alpha$.

We leave it as a problem for the reader to prove that this is an equivalence relation on A and that the distinct equivalence classes are the $A_\alpha$'s.

## 1.2 Mappings

Here we introduce the concept of mapping of one set into another. Without exaggeration that is possibly the single most important and universal notion that runs through all of Mathematics. It is hardly a new thing to any of us, for we have been considering mappings from the very earliest days of our mathematical training.

**DEFINITION:** *If S and T are nonempty sets, then a mapping from S to T is a subset M of $S \times T$ such that for every $s \in S$ there is a unique $t \in T$ such that the ordered pair (s, t) is in M.*

This definition serves to make the concept of a mapping precise for us but we shall almost never use it in this form. Instead we do prefer to think of a mapping as a rule which associates with any element s in S some element t in T, the rule being, associate (or map) $s \in S$ with $t \in T$ if and only if $(s, t) \in M$. We shall say that t is the image of s under the mapping.

Now for some notation for these things. Let $\sigma$ be a mapping from S to T; we often denote this by writing $\sigma: S \to T$ or $S \xrightarrow{\sigma} T$. If t is the image of s under $\sigma$ we shall sometimes write this as $\sigma: s \to t$; more often we shall represent this fact by $t = s\sigma$ or $t = \sigma(s)$. Algebraists often write mappings on the right, other mathematicians write them on the left.

**DEFINITION:** *The mappings $\tau$ of S into T is said to be onto T if given $t \in T$ there exists an element $s \in S$ such that $t = s\tau$.*

**DEFINITION:** *The mapping $\tau$ of S into T is said to be a one to one mapping if whenever $s_1 \neq s_2$; then $s_1\tau \neq s_2\tau$. In terms of inverse images, the mapping $\tau$ is one-to-one if for any $t \in T$ the inverse image of t is either empty or is a set consisting of one element.*

**DEFINITION:** *Two mappings $\sigma$ and $\tau$ of S into T are said to be equal if $s\sigma = s\tau$ for every $s \in S$.*

**DEFINITION:** *If $\sigma: S \to T$ and $\tau: T \to U$ then the composition of $\sigma$ and $\tau$ (also called their product) is the mapping $\sigma \circ \tau: S \to U$ defined by means of $s(\sigma \circ \tau) = (s\sigma)\tau$ for every $s \in S$.*

Note the following example is very important in a sense that we shall be using it in almost all the chapters of this book. The example deals with nothing but mapping of a set of n elements to itself.

*Example 1.2.1:* Let (1, 2, 3) be the set S. Let S(3) denote the set of all mappings of S to itself. The number of elements in S(3) is $27 = 3^3$.

***Example 1.2.2:*** Let $S = (1, 2, 3, ... , n)$. The set of all mappings of S to itself denoted by $S(n)$ has $n^n$ elements in it. Throughout this book by $S(n)$ we mean only the set of all mappings of the set $S = (1, 2, 3, ... , n)$ to itself. We are now interested to study the set $S = (1, 2, 3, ... , n)$ when we take only one to one mappings of S to itself.

***Example 1.2.3:*** Let $S = (1, 2, 3)$ be the set of all one to one mappings of S to itself, denoted by $S_3$. $S_3$ contains only 6 elements i.e. 3! elements. It is very important to note that the set $S_3$ is a proper subset of the set $S(3)$.

***Example 1.2.4***: Let $S = (1, 2, ... ,n)$. By $S_n$ we denote the set of all one to one mappings of S to itself. Clearly, $S_n$ has n! elements and we know the set of all one to one mappings of S is a proper subset of the set of all mappings of S to itself that is $S_n \subset S(n)$.

Now it is very important to note that the composition of mappings of the set to itself in general is not commutative that is if $\sigma$ and $\tau: S \to S$, $\sigma \circ \tau \neq \tau \circ \sigma$. Secondly we see if the mappings $\sigma, \tau, \delta : S \to S$ then we have $(\sigma \circ \tau) \circ \delta = \sigma \circ (\tau \circ \delta)$ that is the composition of mappings in general satisfies the associative law. Further if $\sigma, \tau : S \to S$ is one to one then $\sigma \circ \tau$ and $\tau \circ \sigma$ are also one to one maps from S to S. We say a mapping $\sigma : S \to S$ is the identity map if $\sigma(s) = s$ for all $s \in S$. In this book we denote the identity map just by 1.

**DEFINITION:** *Given an arbitrary set S we call a mapping $\sigma$ from a set $S \times S$ into S a binary operation on S. Given such a mapping $\sigma : S \times S \to S$ we could use it to define a "product" in S by declaring $a \circ b = c$ if $\sigma(a, b) = c$.*

**DEFINITION:** *If $\sigma$ is a one to one mapping from the set S into S we have a one to one mapping $\mu$ from S to S such that $\sigma \circ \mu = \mu \circ \sigma$ is the identity map from S to S. We call $\mu$ the inverse map of $\sigma$ and $\sigma$ is called the converse map of $\mu$.*

***Example 1.2.5:*** Let $\sigma, \mu : S \to S$ where $S = (1, 2, 3)$ be given by

$\sigma(1) = 2$          $\mu(1) = 3$
$\sigma(2) = 3$          $\mu(2) = 1$
$\sigma(3) = 1$          $\mu(3) = 2$.

Now, $(\sigma \circ \mu)(1) = 1$, $(\sigma \circ \mu)(2) = 2$ and $(\sigma \circ \mu)(3) = 3$. Thus we see $\sigma \circ \mu =$ identity map. Similarly $\mu \circ \sigma =$ identity map.

## 1.3 Semigroups and Smarandache Semigroups

In this section, we just recall the definitions of semigroups and Smarandache semigroups. Semigroups are the algebraic structures in which are defined a binary operation which is both closed and associative. Already we have defined binary operation but we have not mentioned explicitly that those maps from $S \times S$ to S where S is an arbitrary set are binary operations. Further, we mention here that semigroups are the most generalized structures with a single binary operation defined on it.

**DEFINITION:** *Let S be a non empty set on which is defined a binary operation 'o', (S, o) is a semigroup if 1) for all a, b ∈ S we have a o b = c ∈ S and 2) a o (b o c) = (a o b) o c for all a, b, c ∈ S.*

***Example 1.3.1:*** Let $Z^+$ = {1, 2, 3, ...} be the set of positive integers define '+', usual addition of integers on $Z^+$ as the binary operation, $(Z^+, +)$ is a semigroup.

***Example 1.3.2:*** Let $Z_n$ = {0, 1, 2, ... , n-1} the set of integers modulo n. $Z_n$ under multiplication is a semigroup denoted by $(Z_n, \times)$

If a semigroup (S, o) has an element e such that s o e = e o s = s for all s ∈ S then we say e is the identity element of S relative to the binary operation o and we call S a monoid or simply a semigroup. A semigroup with the identity element is called a monoid. Clearly the semigroup in example 1.3.1 is not a monoid but the semigroup given in example 1.3.2 is a monoid for 1 acts the multiplicative identity for $Z_n$. $(Z_n, x)$ is a monoid.

***Example 1.3.3:*** Let S = (1, 2). Then S(2) is the semigroup with four elements, the elements of S(2) are the mappings of S to S given by $\sigma_1, \sigma_2, \sigma_3$ and $\sigma_4$; where $\sigma_1(1)$ = 1, $\sigma_2(1)$ = 2, $\sigma_3(1)$ = 1, $\sigma_4(1)$ = 2, $\sigma_1(2)$ = 1, $\sigma_2(2)$ = 2, $\sigma_3(2)$ = 2 and $\sigma_4(2)$ = 1. Thus S(2) = $(\sigma_1, \sigma_2, \sigma_3, \sigma_4)$.

***Example 1.3.4:*** Let S = (1, 2, ... , n). S(n) is the set of all mappings of S to S. S(n) under the composition of mappings is a semigroup with $n^n$ elements in it.

***Example 1.3.5:*** Let $Z^+$ = {1, 2, ... , n,...} be the set of positive integers $(Z^+, \times)$ is a monoid where '$\times$' is the usual multiplication of integers. Now it is important to remark $(Z^+, \times)$ is a monoid where as $(Z^+, +)$ is not a monoid as zero the identity with respect to addition does not belong to $Z^+$. Thus on the same set different binary operations may be defined so as to make the semigroup a monoid or vice versa.

***Example 1.3.6:*** Let $Z_{10}$ = {0, 1,2, ... , 9}. $(Z_{10}, \times)$ is a monoid under the usual multiplication '$\times$' modulo 10. Thus, every monoid is obviously a semigroup with the identity. Throughout this book, even if the semigroup has the identity we choose to call it only a semigroup. Definition of a Smarandache semigroup is taken from the paper of Padilla Raul, for he was the first one to introduce the notion of Smarandache algebraic structures in the year 1998. His first definition on Smarandache algebraic structures is the Smarandache semigroup.

**DEFINITION:** *The Smarandache semigroup is defined to be a semigroup A such that a proper subset of A is a group (with respect to the same binary operation on A).*

Here we just use the concept of the algebraic structure, viz. group and we have not defined it in the first chapter, but as this text is not just for a beginner or a graduate we make this lapse; for the concept of group is elaborately dealt in chapter 2. From the definition of the Smarandache semigroup we see every Smarandache semigroup is obviously a semigroup as semigroup is used to define the concept of Smarandache semigroup. To make this definition explicit we illustrate them by the following examples:

*Example 1.3.7:* Let $Z_8 = \{0, 1, 2, \ldots, 7\}$ be the semigroup under multiplication mod 8. Clearly the set $A = (1, 7) \subset Z_8$ is a group so $Z_8$ is the Smarandache semigroup.

*Example 1.3.8:* Let $Z_4 = \{0, 1, 2, 3\}$ be the semigroup under multiplication modulo 4. Clearly the subset $A = \{1, 3\} \subset Z_4$ under multiplication modulo 4 is a group. So $Z_4$ is a Smarandache semigroup.

*Example 1.3.9:* Let $Z_7 = \{0, 1, 2, \ldots, 6\}$ be the semigroup under multiplication '$\times$' modulo 7. The only proper subsets of $Z_7$ which are groups under '$\times$' are $S = \{1, 2, 3, \ldots, 6\}$ and $B = \{1, 6\}$. Thus, $Z_7$ is a Smarandache semigroup.

*Example 1.3.10:* Let $S_{2\times 2} = \{(a_{ij}) \mid a_{ij} \in Z_2 = \{0, 1\}\}$ be the set of all 2×2 matrices with entries from the prime field $Z_2 = \{0, 1\}$. $S_{2\times 2}$ is a semigroup under matrix multiplication $A = \left\{ \begin{pmatrix} 1 & 0 \\ 0 & 1 \end{pmatrix}, \begin{pmatrix} 0 & 1 \\ 1 & 0 \end{pmatrix} \right\} \subset S_{2\times 2}$ is a subgroup under matrix multiplication.

Thus, $S_{2\times 2}$ is a Smarandache semigroup. We are not expecting the algebraist to understand fully the properties and notions of the Smarandache semigroup by this definition and illustrations, for we are going to deal elaborately about Smarandache semigroups and Smarandache notions in groups in Chapters 4, 5 and 6 with more illustrations.

**PROBLEM 1:** Given an example of a Smarandache semigroup of order 42.

**PROBLEM 2:** Find all subgroups in the Smarandache semigroup $Z_{72} = \{0, 1, 2, \ldots, 71\}$ multiplication modulo 72.

**PROBLEM 3:** Is $Z_{19} = \{0, 1, 2, \ldots, 18\}$ a Smarandache semigroup under multiplication?

**PROBLEM 4:** Give the subgroups $S(5)$, $S(5)$ the Smarandache semigroup of mappings of the set $S = (1, 2, 3, 4, 5)$ to itself, under composition of maps.

**PROBLEM 5:** Find all subgroup of $S_{2\times 2}$ given in example 1.3.10.

With this we end this chapter suggesting the reader the following materials as supplementary reading:

### Supplementary Reading

1. Birkhoff G. and S. Maclane, A Brief Survey of Modern Algebra, 2$^{nd}$ Edition, New York, Macmillan, 1965.

2. Padilla Raul, Smarandache algebraic structures, Bull of Pure and applied Sciences, Delhi, Vol. 17E, No. 1, 119-121, 1998.

## CHAPTER TWO

# ELEMENTARY PROPERTIES OF GROUPS

In this chapter we introduce the notion of groups and recall some of the elementary properties of groups. In later chapters we will be defining a Smarandache analog of these properties whenever possible. Further, we give examples after each definition to make the reader understand easily.

## 2.1 Definition of a Group

It is a well-known fact that groups are the only algebraic structures with a single binary operation that is mathematically so perfect that an introduction of a richer structure within it is impossible. Now we proceed on to define a group.

**DEFINITION:** *A non empty set of elements G is said to form a group if in G there is defined a binary operation, called the product and denoted by '•' such that*

1. *a, b $\in$ G implies that a • b $\in$ G (closed)*
2. *a, b, c $\in$ G implies a • (b • c) = (a • b) • c (associative law)*
3. *There exists an element e $\in$ G such that a • e = e • a = a for all a $\in$ G (the existence of identity element in G).*
4. *For every a $\in$ G there exists an element $a^{-1}$ $\in$ G such that a • $a^{-1}$ = $a^{-1}$ • a = e (the existence of inverse in G).*

**DEFINITION:** *A group G is said to be abelian (or commutative) if for every a, b $\in$ G; a • b = b • a.*

A group, which is not abelian, is called naturally enough, non-abelian. Another natural characteristic of a group G is the number of elements it contains. We call this the order of G and denote it by o(G). The number is most interesting when it is finite. In that case, we say that G is a finite group.

## 2.2 Some Examples of Groups

*Example 2.2.1:* G consists of real numbers 1 and $-1$. G under multiplication is a group of order 2 and it is abelian.

*Example 2.2.2:* Let G be the set of all $2 \times 2$ matrices $\begin{pmatrix} a & b \\ c & d \end{pmatrix}$ where a, b, c, d are real numbers such that $ad - bc \neq 0$. G is a group under matrix multiplication with $\begin{pmatrix} 1 & 0 \\ 0 & 1 \end{pmatrix}$ as its identity. G is a non-commutative group.

*Example 2.2.3:* Let $S_3$ be the set of all 1–1 mappings of the set $\{x_1, x_2, x_3\}$ onto it self, under the product called composition of mappings, $S_3$ is a group of order 6. The map

$$p_1 : \begin{cases} x_1 \to x_1 \\ x_2 \to x_3 \\ x_3 \to x_2 \end{cases}$$

is denoted by

$$p_1 = \begin{pmatrix} x_1 & x_2 & x_3 \\ x_1 & x_3 & x_2 \end{pmatrix}, \; p_2 = \begin{pmatrix} x_1 & x_2 & x_3 \\ x_3 & x_2 & x_1 \end{pmatrix}, \; p_3 = \begin{pmatrix} x_1 & x_2 & x_3 \\ x_2 & x_1 & x_3 \end{pmatrix},$$

$$p_4 = \begin{pmatrix} x_1 & x_2 & x_3 \\ x_2 & x_3 & x_1 \end{pmatrix}, \; p_5 = \begin{pmatrix} x_1 & x_2 & x_3 \\ x_3 & x_1 & x_2 \end{pmatrix} \text{ and } e = \begin{pmatrix} x_1 & x_2 & x_3 \\ x_1 & x_2 & x_3 \end{pmatrix}.$$

$p_i \circ e = e \circ p_i = p_i$ for $i = 1, 2, 3, 4, 5$. $S_3$ is the smallest non-commutative group of order 6.

### 2.3 Some Preliminary Results

Here we just prove some interesting preliminary results about groups, which clearly characterize the algebraic structure of groups.

**THEOREM 2.3.1:** Let G be a group, then the identity element of G is unique.

*Proof:* Given G is a group. To prove the identity element of G is unique, we will show that if two elements e and f in G enjoy the property that for every $a \in G$, $a \cdot e = e \cdot a = a$ and $a \cdot f = f \cdot a = a$ that is $a \cdot e = e \cdot a = a \cdot f = f \cdot a$ implies $e = f$. Since $e \cdot a = a$ for every $a \in G$ in particular we have $e \cdot f = f$. However, on the other hand since $b \cdot f = f$ for every $b \in G$, we must have for $e \in G$, $e \cdot f = e$. Piecing these two bits of information together we get $f = e \cdot f = e$ and so $e = f$.

**THEOREM 2.3.2:** If G is a group, then every $a \in G$ has a unique inverse in G.

*Proof:* Let $a \in G$ suppose we have $x, y \in G$ such that $x \cdot a = a \cdot x = e$ and $y \cdot a = a \cdot y = e$ to prove $x = y$. Suppose that for a in G. $a \cdot x = e$ and $a \cdot y = e$ then obviously $a \cdot x = a \cdot y$. Let us make this our starting point, that is, assume that $a \cdot x = a \cdot y$ for a, x, y in G. There is an element $b \in G$ such that $b \cdot a = e$ (as far as we know yet there may be several such b's). Thus $b \cdot (a \cdot x) = b \cdot (a \cdot y)$ using the associative law this leads to $x = e \cdot x = (b \cdot a) \cdot x = b \cdot (a \cdot x) = b \cdot (a \cdot y) = (b \cdot a) \cdot y = e \cdot y = y$. We have, in fact, proved that $a \cdot x = a \cdot y$ in a group forces $x = y$.

Similarly $x \cdot a = y \cdot a$ implies $x = y$. This says that we can cancel from the same side, in equations in groups. However it is important to note that $a \cdot x = y \cdot a$ does not imply $x = y$.

**THEOREM 2.3.3:** Let G be a group; for every $a \in G$, $(a^{-1})^{-1} = a$.

*Proof:* This simply follows from the fact $a^{-1} \cdot (a^{-1})^{-1} = e = a^{-1} \cdot a$ canceling off the $a^{-1}$ we get $(a^{-1})^{-1} = a$. This is analogous to the very familiar result $-(-5) = 5$.

**THEOREM 2.3.4:** Let G be a group. For a, b ∈ G  $(a \bullet b)^{-1} = b^{-1} \bullet a^{-1}$.

*Proof:* Now $(a \bullet b) \bullet (b^{-1} \bullet a^{-1}) = a \bullet (b \bullet b^{-1}) \bullet a^{-1} = a \bullet e \bullet a^{-1} = a \bullet a^{-1} = e$ so by the very definition of the inverse $(a \bullet b)^{-1} = b^{-1} \bullet a^{-1}$. Consequent of the four theorems proved we give the proof of the following results as problems to the reader.

**PROBLEM 1:** Given G is a group. For a, b ∈ G prove the equation a • x = b and y • a = b have unique solution for x and y in G.

**PROBLEM 2:** Prove in a group G

a • u = a • w implies u = w and
u • a = w • a implies u = w for a, u, w ∈ G.

## 2.4 Subgroups

In general, we shall not be interested in subsets of a group G for they do not reflect the fact that G has an algebraic structure imposed on it. Whatever subsets we do consider, will be those endowed with algebraic properties derived from those of G. Smarandache structures are built in a reverse way. We will see later in this book how a Smarandache semigroup is defined.

**DEFINITION:** A non empty subset H of a group G is said to be a subgroup of G if, under the product in G, H itself forms a group.

The following remark is clear; if H is a subgroup of G and K is a subgroup of H, then K is a subgroup of G.

*Example 2.4.1:* Let G = {1, -1} be the group under multiplication H = {1} is a subgroup of G.

We call this subgroup improper or trivial subgroup of G. Thus for every group G the identity element e of G is a subgroup which we call as trivial or improper subgroup of G. Likewise G the group itself is a subgroup of G called the improper subgroup of G. So H a subset of G is called a proper subgroup of G if H is not the identity subgroup or H is not the whole group G.

*Example 2.4.2:* Let $S_3$ = {e, $p_1$, $p_2$, $p_3$, $p_4$, $p_5$} be the group given in example 2.2.3. Clearly H = {e, $p_1$} and K = {e, $p_4$, $p_5$} are subgroups of $S_3$. Both the subgroups are proper subgroups of $S_3$.

*Example 2.4.3:* Let Z = $\{\ldots -2, -1, 0, 1, 2 \ldots\}$ be the set of integers positive, negative with zero. Z under addition is a group. 2Z is a subgroup of Z which is a proper subgroup of Z, as 2Z = $\{\ldots -4, -2, 0, 2, 4, \ldots\}$ is a proper subset of Z.

**THEOREM 2.4.1:** *A non-empty subset H of the group G is a subgroup of G if and only if*

1. *a, b ∈ H implies that a • b ∈ H.*

2.   $a \in H$ implies that $a^{-1} \in H$.

*Proof:* Clearly if H is a subgroup of G then it is obvious 1 and 2 holds good. Conversely suppose 1 and 2 hold good, to establish H is a subgroup of G, is left for the reader as a problem.

**THEOREM 2.4.2:** *If H is a non-empty finite subset of group G and H is closed under multiplication, then H is a subgroup of G.*

*Proof:* From the theorem 2.4.1 we need but show that whenever $a \in H$, then $a^{-1} \in H$. Suppose that $a \in H$; then $a^2 = a \bullet a \in H$, $a^3 = a^2 \bullet a \in H$, ... , $a^m \in H$, ... since H is closed. Thus the infinite collection of elements a, $a^2$, ... , $a^m$, ... must all fit into H, which is a finite subset of G.

Thus there must be repetition in this collection of elements; that is, for some integers r, s with $r > s > o$ $a^r = a^s$. By the cancellation in G; $a^{r-s} = e$ (whence e is in H) since $r - s - 1 \geq 0$, $a^{r-s-1} \in H$ and $a^{-1} = a^{r-s-1}$ since $a \bullet a^{r-s-1} = a^{r-s} = e$. Thus $a^{-1} \in H$, completing the proof of the theorem.

This theorem has a nice analog in the case of Smarandache semigroups, which will be dealt in the later chapters. As the main aim of this book is to introduce Smarandache notions in groups and use Smarandache semigroups we do not give any problems on group theory further the problems which we give here are some theorems or results on groups which can be easily proved.

**DEFINITION:** Let G be a group. H a subgroup of G: for $a,b \in G$ we say a is congruent to b mod H, written as $a \equiv b \pmod{H}$ if $ab^{-1} \in H$. It is easily verified that the relation $a \equiv b \pmod{H}$ is an equivalence relation.

**DEFINITION:** If H is a subgroup of G, $a \in G$, then Ha = {ha / $h \in H$} Ha is called a right coset of H in G.

**THEOREM 2.4.3:** For all $a \in G$, Ha = {$x \in G$ / $a \equiv x \bmod H$}

*Proof:* Let [a] = {$x \in G$ / $a \equiv x \bmod H$}. We first show that Ha $\subset$ [a]. For, if $h \in H$, then $a (ha)^{-1} = a (a^{-1} h^{-1}) = h^{-1} \in H$ since H is a subgroup of G. By the definition of congruence mod H this implies that ha $\in$ [a] for every $h \in H$, and so Ha $\subset$ [a]. Suppose, now, that $x \in$ [a].

Thus $ax^{-1} \in H$, so $(ax^{-1})^{-1} = xa^{-1}$ is also in H. That is, $xa^{-1} = h$ for some $h \in H$. Multiplying both sides by a from the right we get x = ha and so x $\in$ Ha. Thus [a] $\subset$ Ha. Having proved both the inclusions [a] $\subset$ Ha and Ha $\subset$ [a] we can conclude Ha = [a]. Hence the claim.

**DEFINITION:** *If G is a group and $a \in G$, the order of a is the least positive integer m such that $a^m = e$.*

If no such integer exists we say that a is of infinite order. We use the notation o(a) for the order of a.

**DEFINITION:** A subgroup N of a group G is said to be a normal subgroup of G if for every $g \in G$ and $n \in N$, $gng^{-1} \in N$.

Equivalently by $gNg^{-1}$ we mean the set of all $gng^{-1}$, $n \in N$ then N is a normal subgroup of G if and only if $gNg^{-1} \subset N$ for every $g \in G$.

**THEOREM 2.4.4:** N is a normal subgroup of G if and only if $gNg^{-1} = N$ for every $g \in G$.

*Proof:* If $gNg^{-1} = N$ for every $g \in G$, certainly $gNg^{-1} \subset N$, so N is normal in G. Suppose that N is normal in G. Thus if $g \in G$, $gNg^{-1} \subset N$ and $g^{-1}Ng = g^{-1}N(g^{-1})^{-1} \subset N$. Now, since $g^{-1}Ng \subset N$, $N = g(g^{-1}Ng)g^{-1} \subset gNg^{-1} \subset N$ whence $N = gNg^{-1}$.

**DEFINITION:** *A mapping $\phi$ from a group G to a group $\overline{G}$ is said to be a group homomorphism if for all $a, b \in G$, $\phi(ab) = \phi(a)\phi(b)$.*

*Remarks:* If $\phi(a) = e$ for all $a \in G$. We call $\phi$ a trivial homomorphism. Likewise if $\phi$ is a map from G to $\overline{G}$ such that $\phi(x) = x$ for every $x \in G$ then also we say $\phi$ is a trivial homomorphism or the identity homomorphism from G to $\overline{G}$. The kernel of $\phi$, G $\rightarrow$ $\overline{G}$ is defined by $K_\phi = \{x \in G \, / \, \phi(x) = \overline{e}, \overline{e}$ is the identity element of $\overline{G}\}$.

**DEFINITION:** *A homomorphism G into $\overline{G}$ is said to be an isomorphism if $\phi$ is one to one.*

**DEFINITION:** *By an automorphism of a group G, we shall mean an isomorphism of G onto itself.*

The following is left as a problem for the reader.

**PROBLEM:** If G is a group then prove A(G), the set of all automorphism of G is also a group.

**DEFINITION:** *If $a, b \in G$, then b is said to be a conjugate of a in G if there exists an element $c \in G$ such that $b = c^{-1}ac$. We shall write, for this $a \sim b$ and shall refer to this relation as conjugacy.*

**THEOREM 2.4.5:** *Conjugacy is an equivalence relation on G.*

*Proof:* As usual, in order to establish this, we must prove that

1. $a \sim a$;
2. $a \sim b$ implies $b \sim a$;
3. $a \sim b$, $b \sim c$ implies $a \sim c$ for all $a, b, c$ in G.

We prove each of these in turn

1. Since $a = e^{-1}ae$, $a \sim a$ with $c = e$ serving as the c in the definition of conjugacy.

2. If a ~ b then b = $x^{-1}ax$ for some x ∈ G, hence a = $(x^{-1})^{-1}bx^{-1}$ and since y = $x^{-1}$ ∈ G, a = $y^{-1}by$, b ~ a follows.
3. Suppose that a ~ b and b ~ c where a, b, c ∈ G. Then b = $x^{-1}ax$, c = $y^{-1}by$ for some x, y ∈ G. Substituting for b in the expression for c we obtain c = $y^{-1}(x^{-1}ax)y = (xy)^{-1}a(xy)$; since xy ∈ G, a ~ c is a consequence.

For a ∈ G let C(a) = {x ∈ G / a ~ x} C(a), the equivalence class of a in G under our relation, usually called the conjugate class of a in G; it contains the set of all distinct elements of the form $y^{-1}ay$ as y ranges over G.

**DEFINITION:** *If a ∈ G, then N(a), the normalizer of a in G, is a set N(a) = {x ∈ G / xa = ax}.*

N(a) consists of precisely those elements in G which commute with a. It is left as a problem for the reader to prove N(a) is a subgroup of G.

**THEOREM 2.4.6:** *If G is a finite group then $C_a$ = o(G) / o(N(a)); in other words, the number of elements conjugate to a in G is the index of the normalizer of a in G.*

*Proof:* To begin with the conjugate class of a in G, C(a) consists exactly of all the elements $x^{-1}ax$ as x ranges over G. $C_a$ measures the number of distinct $x^{-1}ax$'s. Our method of proof will be to show that two elements in the same right coset of N(a) in G yield the same conjugate of a whereas two elements in different right cosets N(a) in G give rise to different conjugates of a.

In this way we shall have a one to one correspondence between conjugates of a and right cosets of N(a). Suppose x, y ∈ G are in the same right coset of N(a) in G. Thus y = nx where n ∈ N(a) and so na = an. Therefore, since $y^{-1} = (nx)^{-1} = x^{-1}n^{-1}$, $y^{-1}ay = x^{-1}n^{-1}anx = x^{-1}n^{-1}nax = x^{-1}ax$, whence x and y result in the same conjugate of a.

If, on the other hand x and y are in different right cosets of N(a) in G we claim that $x^{-1}ax \ne y^{-1}ay$. Were this is not the case, from $x^{-1}ax = y^{-1}ay$, we would deduce that $yx^{-1}a = ayx^{-1}$; this in turn would imply that $yx^{-1}$ ∈ N(a). However, this declares x and y to be in the same right coset of N(a) in G, contradicting the fact that they are in different cosets. This proof is now complete.

Since o(G) = $\sum C_a$ using the theorem we have

$$o(G) = \sum \frac{o(G)}{o(N(a))}$$

where the sum runs over one element a in each conjugate class. This is known as the class equation.

**DEFINITION:** *Let G be a group. Z(G) = {x ∈ G | gx = xg for all g ∈ G}. Then Z(G) is called the center of the group G.*

**DEFINITION:** *Let G be a group, A, B be subgroups of G. If x, y ∈ G define x ~ y if y = axb for some a ∈ A and b ∈ B. We call the set AxB = {axb / a ∈ A, b ∈ B} a double coset of A, B in G.*

It is left as a problem for the reader to prove the relation defined above is an equivalence relation on G. The equivalence class of x ∈ G is the set AxB = {axb / a ∈ A, b ∈ B}.

If A, B are finite subgroups of G, how many elements are there in the double coset AxB? It is again left for the reader to verify.

$$o(AxB) = \frac{o(A)o(B)}{o(A \cap xBx^{-1})}$$

**DEFINITION**: Let G be a group. A and B subgroups of G, we say A and B are conjugate with each other if for some $g \in G$, $A = gBg^{-1}$.

Clearly if A and B are conjugate subgroups of G then o(A) = o(B).

**DEFINITION:** *Let $G_1, \ldots, G_n$ be any n groups. Let $G = G_1 \times \ldots \times G_n = \{(g_1, g_2, \ldots, g_n) / g_i \in G_i\}$ be the set of all ordered n-tuples, that is, the cartesian product of $G_1, G_2, \ldots, G_n$.*

*We define a product in G via $(g_1, g_2, \ldots, g_n)(g'_1, g'_2, \ldots, g'_n) = (g_1 g'_1, g_2 g'_2, \ldots, g_n g'_n)$ that is, component-wise multiplication. The product in the $i^{th}$ component is carried in the group G. Then G is a group in which $(e_1, e_2, \ldots, e_n)$ is the unit element, where each $e_i$ is the unit element of $G_i$, and where $(g_1, g_2, \ldots, g_n)^{-1} = (g_1^{-1}, g_2^{-1}, \ldots, g_n^{-1})$. We call this group G the external direct product of $G_1, \ldots, G_n$.*

In $G = G_1 \times \ldots \times G_n$ let $\overline{G}_i = \{(e_1, e_2, \ldots, e_{i-1}, g_i, e_{i+1}, \ldots, e_n) / g_i \in G_i\}$. Then $\overline{G}_i$ is a normal subgroup of G and is isomorphic to $G_i$. Moreover $G = \overline{G}_1 \overline{G}_2 \ldots \overline{G}_n$ and every $g \in G$ has a unique decomposition; $g = \overline{g}_1 \overline{g}_2 \ldots \overline{g}_n$ where $\overline{g}_1 \in \overline{G}_1$, $\overline{g}_2 \in \overline{G}_2, \ldots \overline{g}_n \in \overline{G}_n$. We leave the verification of these facts to the reader.

**DEFINITION**: *Let G be a group and $N_1, N_2, \ldots, N_n$ normal subgroups of G such that $G = N_1 N_2 \ldots N_n$. Given $g \in G$ then $g = m_1 m_2 \ldots m_n$, $m_i \in N_i$, g written in this way is unique. We then say that G is the internal direct product of $N_1, N_2, \ldots, N_n$.*

It is left for the reader to verify the following facts.

Let G be the internal direct product of $N_1, \ldots, N_n$. Then for $i \neq j$, $N_i \cap N_j = (e)$ and if $a \in N_i$, $b \in N_j$ then ab = ba.

Let G be a group and suppose that G is the internal direct product of $N_1, \ldots, N_n$ and $T = N_1 \times N_2 \times \ldots \times N_n$, the external direct product of $N_1, \ldots, N_n$. Then prove G and T are isomorphic.

**Supplementary Reading:**

# Chapter Three

# SOME CLASSICAL THEOREMS IN GROUP THEORY

In this chapter, we just recall some theorems in group theory with proofs. The main purpose for giving the proof is that when we try to adopt them for Smarandache semigroups, we would prove either that the classical result or theorem is true or the theorem is not true. Thus, the Lagrange's theorem, Cayley's theorem, Cauchy's theorem and Sylow's theorems are the main theorems which we are interested in proving or disproving in case of Smarandache semigroup or proving the validity in case of Smarandache semigroups.

Throughout this chapter by a symmetric group of degree n denoted by $S_n$ we mean the set of all 1-1 mappings of the set (1,2, ... ,n) to itself and the group operation being the composition of maps. By the dihedral group we mean the group $D_{2n} = \{a, b / a^2 = b^n = 1, bab = a\}$ where this group contains exactly 2n elements.

## 3.1 Lagrange's Theorem

The famous theorem by Lagrange mainly uses the concept of cosets to prove that if G is a finite group and H is a subgroup of G, then o(H) is a divisor of o(G). It might be difficult at this point, for the student to see the extreme importance of this result. As the subject is penetrated more deeply, one will become increasingly aware of its basic character. Here we give the proof of Lagrange's theorem.

**THEOREM 3.1.1: (LAGRANGE).** If G is a finite group and H is a subgroup of G then o(H) is a divisor of o(G).

*Proof:* Suppose G is a finite group and H is a subgroup of G. Let $h_1, h_2, \ldots, h_r$ be a complete list of the elements of H, r = o(H). If H = G, there is nothing to prove. Suppose, that H ≠ G, thus there is an a ∈ G, a ∉ H. List all the elements so far in two rows as

$$h_1, h_2, \ldots, h_r$$
$$h_1a, h_2a, \ldots, h_ra.$$

We claim that all the entries in the second line are different from each other and are different from the entries in the first line. If any two in the second line were equal, then $h_ia = h_ja$ with i ≠ j, but by the cancellation law this would lead to $h_i = h_j$ a contradiction. If an entry in the second line were equal to one in the first line, then $h_ia = h_j$ resulting in a = $h_i^{-1} h_j$ ∈ H since H is a subgroup of G, this violates a ∉ H.

Thus we have, so far, listed 2o(H) elements; if these elements account for all the elements of G, we are done. If not there is an element b ∈ G, b ∉ Ha and b ∉ H that did not occur in these two lines. Consider the new list

$$h_1, h_2, \ldots, h_r$$
$$h_1a, h_2a, \ldots, h_ra$$
$$h_1b, h_2b, \ldots, h_rb$$

As before (we are now waving our hands) we could show that, no two entries in the third line are equal to each other, and that no entry in the line occurs in the first or second line. Thus, we have listed 3o(H) elements. Continuing in this way, every new element introduced, in fact, produced o(H) new elements.

Since G is a finite group, we must eventually exhaust all the elements of G. But if we ended up using k lines to list all elements of the group, we would have written ko(H) distinct elements and so ko(H) = o(G). Hence the claim.

It is essential to point out that for every divisor of the order of a finite group G we need not in general have a subgroup. That the converse to Lagrange's theorem is false - a group G need not have a subgroup of order m if m is a divisor of o(G).

Consider the group $S_4$, the symmetric group of degree 4 which has $A_4$ the alternating subgroup of order 12. Clearly 6/12 but $A_4$ has no subgroup of order 6. Thus, we see the converse of Lagrange's theorem in general is not true. Hence, there are very few results, which assert the existence of subgroups of prescribed order in arbitrary finite groups. The Lagrange's theorem has some very important Corollaries.

**COROLLARY 3.1.2:** *If G is a finite group and a $\in$ G, then o(a) | o(G).*

*Proof:* With Lagrange's theorem already in hand, it seems most natural to prove the corollary by exhibiting a subgroup of G whose order is o(a). The element a itself furnishes us with this subgroup by considering the cyclic subgroup generated by a that is, (a) of G; (a) consists of e, a, $a^2$, .... How many elements are there in (a)?

We assert that this number is the order of a. Clearly, since $a^{o(a)}$ = e, this subgroup has at most o(a) elements. If it should actually have fewer than this number of elements, then $a^i = a^j$ for some integers $0 \le i < j < o(a)$. Then $a^{j-i}$ = e, yet $0 < j - i <$ o(a) which would contradict the very meaning of o(a). Thus the cyclic subgroup generated by a has o(a) elements, whence, by Lagrange's theorem, o(a) | o(G).

**COROLLARY 3.1.3:** *If G is a finite group and a $\in$ G, then $a^{o(G)}$ = e.*

*Proof:* By Corollary 3.1.2, o(a) | o(G); thus o(G) = mo(a). Therefore, $a^{o(G)} = a^{mo(a)} = (a^{o(a)})^m = e^m = e$.

### 3.2 Cauchy's Theorem

In this section we give the two Cauchy's theorems one for abelian groups and the other for non-abelian groups. The main result on finite groups is that if the order of the group is n (n < $\infty$) if p is a prime dividing n by Cauchy's theorem we will always be able to pick up an element a $\in$ G such that $a^p$ = e. In fact we can say

Sylow's theorem is a partial extension of Cauchy's theorem for he says this finite group G has a subgroup of order $p^\alpha (\alpha \geq 1$, p, a prime).

**THEOREM 3.2.1:** (**CAUCHY'S THEOREM FOR ABELIAN GROUPS**). *Suppose G is a finite abelian group and p / o(G), where p is a prime number. Then there is an element a ≠ e ∈ G such that $a^p = e$.*

*Proof:* We proceed by induction over o(G). In other words, we assume that the theorem is true for all abelian groups having fewer elements than G. From this we wish to prove that the result holds for G. To start the induction we note that the theorem is vacuously true for groups having a single element.

If G has no subgroups H ≠ (e), G, must be cyclic of prime order. This prime must be p, and G certainly has p-1 elements a ≠ e satisfying $a^p = a^{o(G)} = e$. So suppose G has a subgroup N ≠ (e), G. If p/o(N), by our induction hypothesis, since o(N)< o(G) and N is abelian, there is an element b ∈ N, b ≠ e, satisfying $b^p = e$; since b ∈ N ⊂ G we would have exhibited an element of the type required. Therefore, we may assume that p ∤ o(N). Since G is abelian, N is a normal subgroup of G, so G/N is a group. Moreover, o(G/N) = o(G)/o(N), since p ∤ o(N),

$$p \left| \frac{o(G)}{o(N)} \right. < o(G).$$

Also, since G is abelian, G/N is abelian. Thus by our induction hypothesis there is an element X ∈ G/N satisfying $X^p = e_1$, ; the unit element of G/N, X ≠ $e_1$. By the very form of elements of G/N, X = Nb, b ∈ G, so that $X^p = (Nb)^p = Nb^p$. Since $e_1$ = Ne, $X^p = e_1$, X ≠ $e_1$ translates into $Nb^p = N$, Nb ≠ N. Thus $b^p$ ∈ N, b ∉ N.

Using one of the corollaries to Lagrange's theorem, $(b^p)^{o(N)} = e$. That is, $(b^p)^{o(N)}$ = e. Let c = $b^{o(N)}$. Certainly $c^p = e$. In order to show that c is an element that satisfies the conclusions of the theorem we must finally show that c ≠ e. However, if c = e, $b^{o(N)}$ = e, and so $(Nb)^{o(N)}$ = N. Combining this with $(Nb)^p = N$, p ∤ o(N), p a prime number, we find that Nb = N, so b ∈ N, a contradiction. Thus c ≠ e, $c^p$ = e, and we have completed the induction. This proves the result.

**THEOREM 3.2.2:** (**CAUCHY**) *If p is a prime number and p | o(G), then G has an element of order p.*

*Proof:* We seek an element a ≠ e ∈ G satisfying $a^p$ = e. To prove its existence we proceed by induction on o(G); that is, we assume the theorem to be true for all groups T such that o(T) < o(G). We need not worry about starting the induction for the result is vacuously true for groups of order 1.

If for any subgroup W of G, W ≠ G, were it to happen that p | o(W), then by our induction hypothesis there would exist an element of order p in W, and thus there would be such an element in G. Thus we may assume that p is not a divisor of the order of any proper subgroup of G. In particular, if a ∉ Z(G), since N(a) ≠ G, p ∤ o(N(a)). Let us write down the class equation:

$$o(G) = o(Z(G)) + \sum_{N(a) \neq G} \frac{o(G)}{o(N(a))}.$$

Since $p \mid o(G)$, $p \nmid o(N(a))$ we have that $p \left| \frac{o(G)}{o(N(a))} \right.$, and so $p \left| \sum_{N(a) \neq G} \frac{o(G)}{o(N(a))} \right.$; Since we also have that $p \mid o(G)$, we conclude that

$$p \left| \left( o(G) - \sum_{N(a) \neq G} \frac{o(G)}{o(N(a))} \right) = o(Z(G)). \right.$$

Z(G) is thus a subgroup of G whose order is divisible by p. But, after all, we have assumed that p is not a divisor of the order of any proper subgroup of G, so that Z(G) cannot be a proper subgroup of G. We are forced to accept the only possibility left us, namely, that Z(G) = G. But then G is abelian; now we invoke the result already established for abelian groups to complete the induction. This proves the theorem.

### 3.3 Cayley's Theorem

Though one may marvel at the number of groups of varying types carrying many different properties, except for Cayley's we would not have seen them to be imbedded in the class of groups this was done by Cayley's in his famous theorem. Smarandache semigroups also has a beautiful analog for Cayley's theorem which will be given in Chapter 5.

By A(S) we mean the set of all one to one maps of the set S into itself. Clearly A(S) is a group having n! elements if $o(S) = n < \infty$, if S is an infinite set, A(S) has infinitely many elements.

**THEOREM 3.3.1: (CAYLEY)** *Every group is isomorphic to a subgroup of A(S) for some appropriate S.*

*Proof:* Let G be a group. For the set S we will use the elements of G; that is, put S = G. If $g \in G$, define $\tau_g : S(= G) \to S(= G)$ by $x\tau_g = xg$ for every $x \in G$. If $y \in G$, then $y = (yg^{-1})g = (yg^{-1})\tau_g$, so that $\tau_g$ maps S onto itself. Moreover, $\tau_g$ is one to one, for if x, y $\in$ S and $x\tau_g = y\tau_g$, then xg = yg, which, by the cancellation property of groups, implies that x = y. We have proved that for every $g \in G$, $\tau_g \in A(S)$.

If g, h $\in$ G, consider $\tau_{gh}$. For any $x \in S = G$, $x\tau_{gh} = x(gh) = (xg)h = (x\tau_g)\tau_h$. Note that we used the associative law in a very essential way here. From $x\tau_{gh} = x\tau_g\tau_h$ we deduce that $\tau_{gh} = \tau_g\tau_h$. Therefore, if $\psi: G \to A(S)$ is defined by $\psi(g) = \tau_g$, the relation $\tau_{gh} = \tau_g\tau_h$ tells us that $\psi$ is a homomorphism. What is the kernel K of $\psi$? If $g_0 \in K$, then $\psi(g_0) = \tau_{g_0}$ is the identity map on S, so that for $x \in G$, and, in particular, for $e \in G$, $e\tau_{g_0} = e$. But $e\tau_{g_0} = eg_0 = g_0$. Thus comparing these two expressions for $e\tau_{g_0}$ we conclude that $g_0 = e$, whence $K = (e)$. We know a homomorphism $\psi$ of G into A(S) with kernel K is an isomorphism of G into A(S) if and only if $K = (e)$, proving the theorem.

## 3.4 Sylow's Theorems

The Norwegian mathematician Peter Ludvig Mejdell Sylow was the contributor of Sylow's theorems. Sylow's theorems serve double purpose. One hand they form partial answers to the converse of Lagrange's theorem and on the other hand they are the complete extension of Cauchy's Theorem. Thus Sylow's work interlinks the works of two great mathematicians Lagrange and Cauchy. The following theorem is one, which makes use of Cauchy's theorem. It gives a nice partial converse to Lagrange's theorem and is easily understood.

**THEOREM 3.4.1: (SYLOW'S THEOREM FOR ABELIAN GROUPS)** *If G is an abelian group of order o(G), and if p is a prime number, such that $p^{\alpha} \mid o(G)$, $p^{\alpha+1} \nmid o(G)$, then G has a subgroup of order $p^{\alpha}$.*

*Proof:* Given G is an abelian group of order o(G) and p is a prime number such that $p^{\alpha}/o(G)$ and $p^{\alpha+1} \nmid o(G)$. Suppose $\alpha = 0$, then the subgroup (e) satisfies the conclusion of the result. So suppose $\alpha \neq 0$. Then $p \mid o(G)$. By Cauchy's theorem for abelian groups, there is an element $a \neq e \in G$ satisfying $a^p = e$.

Let $S = \{x \in G \mid x^{p^n} = e \text{ for some integer n}\}$. Since $a \in S$, $a \neq e$, it follows that $S \neq (e)$. We now assert that S is a subgroup of G. Since G is finite we must only verify that S is closed. If $x, y \in S$, $x^{p^n} = e, y^{p^m} = e$, so that $(xy)^{p^{n+m}} = x^{p^{n+m}} y^{p^{n+m}} = e$ (we have used that G is abelian), proving that $xy \in S$. We next claim that $o(S) = p^{\beta}$ with $\beta$ an integer $0 < \beta \leq \alpha$. For, if some prime $q \mid o(S)$, $q \neq p$, by the result of Cauchy's theorem for abelian groups there is an element $c \in S$, $c \neq e$, satisfying $c^q = e$.

However, $c^{p^n} = e$ for some n since $c \in S$. Since $p^n$, q are relatively prime, we can find integers $\lambda, \mu$ such that $\lambda q + \mu p^n = 1$, so that $c = c^1 = c^{\lambda q + \mu p^n} = (c^q)^{\lambda}(c^{p^n})^{\mu} = e$, contradicting $c \neq e$. By Lagrange's theorem $o(S) \mid o(G)$, so that $\beta \leq \alpha$. Suppose that $\beta < \alpha$; consider the abelian group G/S. Since $\beta < \alpha$ and $o(G/S) = o(G)/o(S)$, $p \mid o(G/S)$, there is an element $Sx$, $(x \in G)$ in G/S satisfying $Sx \neq S, (Sx)^{p^n} = S$ for some integer $n > 0$. But $S = (Sx)^{p^n} = Sx^{p^n}$, and so $x^{p^n} \in S$; consequently $e = (x^{p^n})^{o(S)} = (x^{p^n})^{p^{\beta}} = x^{p^{n+\beta}}$. Therefore, x satisfies the exact requirements needed to put it in S; in other words, $x \in S$. Consequently $Sx = S$ contradicting $Sx \neq S$. Thus $\beta < \alpha$ is impossible and we are left with the only alternative, namely, that $\beta = \alpha$. S is required subgroup of order $p^{\alpha}$.

**COROLLARY 3.4.2:** *If G is abelian of order o(G) and $p^{\alpha} \mid o(G)$, $p^{\alpha+1} \nmid o(G)$, there is a unique subgroup of G of order $p^{\alpha}$.*

*Proof:* Suppose T is another subgroup of G of order $p^\alpha$, $T \neq S$. Since G is abelian ST = TS, so that ST is a subgroup of G. We know if S and T are finite subgroups of G of order o(S) and o(T) respectively then.

$$o(ST) = \frac{o(S)o(T)}{o(S \cap T)} = \frac{p^\alpha p^\alpha}{o(S \cap T)}$$

and since $S \neq T$, $o(S \cap T) < p^\alpha$, leaving us with $o(ST) = p^\gamma$, $\gamma > \alpha$. Since ST is a subgroup of G, $o(ST) \mid o(G)$; thus $p^\gamma \mid o(G)$ violating the fact that $\alpha$ is the largest power of p which divides o(G). Thus no such subgroup T exists, and S is the unique subgroup of order $p^\alpha$.

**DEFINITION:** *Let G be a finite group. A subgroup G of order $p^\alpha$, where $p^\alpha \mid o(G)$ but $p^\alpha \nmid o(G)$, is called a p-Sylow subgroup of G. Thus we see that for any finite group G if p is any prime which divides o(G); then G has a p-Sylow subgroup.*

Thus, the classical three parts of theorems due to Sylow with proofs will be given in this chapter. It is interesting to note that out of three proofs were given to the first Sylow's theorem, which clearly enables us to understand that Sylow's theorem is that important that it merits this multi front approach. However, in this text we give only the proof, which uses induction and the class equation.

**THEOREM 3.4.3: (FIRST PART OF SYLOW'S THEOREM).** *If p is a prime number and $p^\alpha \mid o(G)$ and $p^{\alpha+1} \nmid o(G)$, G is a finite group, then G has a subgroup of order $p^\alpha$.*

*Proof:* We give the proof using induction on the order of the group G, that for every prime p dividing the order of G, G has a p-Sylow subgroup. If the order of the group G is 2, the only relevant prime is 2 and the group certainly has a subgroup of order 2, namely itself. So we suppose the result to be correct for all groups of order less than o(G). From this we want to show that the result is valid for G. Suppose, then, that $p^\alpha \mid o(G)$, $p^{\alpha+1} \nmid o(G)$ where p is a prime, $\alpha \geq 1$. If $p^\alpha \mid o(H)$ for any subgroup H of G, where $H \neq G$, then by the induction hypothesis, H would have a subgroup T of order $p^\alpha$.

However, since T is a subgroup of H and H is a subgroup of G, T too is a subgroup of G. But then T would be the sought after subgroup of order $p^\alpha$. We therefore may assume that $p^\alpha \nmid o(H)$ for any subgroup H of G, where $H \neq G$. We restrict our attention to a limited set of such subgroups. Recall that if $a \in G$, then $N(a) = \{x \in G \mid xa = ax\}$ is a subgroup of G; moreover, if $a \notin Z(G)$, the center of G, then $N(a) \neq G$. Recall, too, that the class equation of G states that $o(G) = \sum \frac{o(G)}{o(N(a))}$ where the sum runs over one element a from each conjugate class. We separate this sum into two pieces, those a which lie in Z(G), and those which don't. This gives, $o(G) = z + \sum_{a \notin Z} \frac{o(G)}{o(N(a))}$ where $z = o(Z(G))$. Now invoke the reduction we have made namely, that $p^\alpha \nmid o(H)$ for any subgroup H of G, where $H \neq G$, to those subgroups N(a) for $a \notin Z(G)$. Since in this case, $p^\alpha \mid o(G)$ and $p^\alpha \nmid o(N(a))$, we must have that $p \left| \dfrac{o(G)}{o(N(a))} \right.$. Restating this result $p \left| \dfrac{o(G)}{o(N(a))} \right.$, for every $a \in G$ where $a \notin Z(G)$. Look at the class equation with this information. Since $p^\alpha / o(G)$ we have that $p/o(G)$, also $p \left| \displaystyle\sum_{a \notin Z} \dfrac{o(G)}{o(N(a))} \right.$.

Thus the class equation gives us that $p/z$. Since $p/z = o(Z(G))$ by Cauchy's Theorem $Z(G)$ has an element $b \neq e$ of order $p$. Let $B = (b)$, the subgroup of $G$ generated by $b$. $B$ is of order $p$; moreover, since $b \in Z(G)$, $B$ must be normal in $G$. Hence we can form the quotient group $\overline{G} = G/B$. We look at $\overline{G}$, first of all its order if $o(G)/o(B) = o(G)/p$, hence is certainly less than $o(G)$. Secondly we have $p^{\alpha-1}/o(\overline{G})$ but $p^\alpha \nmid o(\overline{G})$. Thus by the induction hypothesis $\overline{G}$ has subgroup $\overline{P}$ of order $p^{\alpha-1}$. Let $P = \{x \in G / xB \in \overline{P}\}$. It is left for the reader to prove; $P$ is a subgroup of $G$.

Moreover, $\overline{P} \approx P/B$ (Prove!); thus $p^{\alpha-1} = o(\overline{P}) = \dfrac{o(P)}{o(B)} = \dfrac{o(P)}{p}$. This results in $o(P) = p^\alpha$. Therefore $P$ is the required p-Sylow subgroup of $G$. This completes the induction and so proves the theorem.

**THEOREM 3.4.4: (SECOND PART OF SYLOW'S THEOREM)** If $G$ is a finite group, $p$ a prime and $p^n \mid o(G)$ but $p^{n+1} \nmid o(G)$, then any two subgroup of $G$ of order $p^n$ are conjugate.

*Proof:* Let $A, B$ be subgroups of $G$, each of order $p^n$. We want to show that $A = gBg^{-1}$ for some $g \in G$. Decompose $G$ into double cosets of $A$ and $B$; $G = \bigcup AxB$. We know $c(AxB) = \dfrac{o(A)o(B)}{o(A \cap xBx^{-1})}$. If $A \neq xBx^{-1}$ for every $x \in G$ then $o(A \cap xBx^{-1}) = p^m$ where $m < n$. Thus $o(AxB) = \dfrac{o(A)o(B)}{p^m} = \dfrac{p^{2n}}{p^m} = p^{2n-m}$ and $2n - m \geq n + 1$. Since $p^{n+1} \mid o(AxB)$ for every $x$ and since $o(G) = \sum o(AxB)$, we would get the contradiction $p^{n+1} \mid o(G)$. Thus $A = gBg^{-1}$ for some $g \in G$. This is the assertion of the theorem.

**THEOREM 3.4.5: (THIRD PART OF SYLOW'S THEOREM)** *The number of p-Sylow subgroups in G, for a given prime, is of the form $1 + kp$.*

Proof: Let $P$ be a p- Sylow subgroup of $G$. We decompose $G$ into double cosets of $P$ and $P$. Thus $G = \bigcup PxP$. We know that

$$o(PxP) = \dfrac{o(P)^2}{o(P \cap xPx^{-1})}.$$

Thus, if $P \cap xPx^{-1} \neq P$ then $p^{n+1} \mid o(PxP)$, where $p^n = o(P)$. Paraphrasing this: if $x \notin N(P)$ then $P^{n+1} \mid o(PxP)$. Also, if $x \in N(P)$, then $PxP = P(Px) = P^2x = Px$, so $o(PxP) = p^n$ in this case. Now

$$o(G) = \sum_{x \in N(P)} o(PxP) + \sum_{x \notin N(P)} o(PxP),$$

where each sum runs over one element from each double coset. However, if $x \in N(P)$, since $PxP = Px$, the first sum is merely $\Sigma_{x \in N(p)}\, o(Px)$ over the distinct cosets of $P$ in $N(P)$. Thus this first sum is just $o(N(P))$. What about the second sum? We saw that each of its constituent terms is divisible by $p^{n+1}$, hence

$$p^{n+1} \;\Big|\; \sum_{x \notin N(P)} o(PxP).$$

We can thus write this second sum as

$$\sum_{x \notin N(P)} o(PxP) = p^{n+1} u.$$

Therefore $o(G) = o(N(P)) + p^{n+1}u$, so

$$\frac{o(G)}{o(N(P))} = 1 + \frac{p^{n+1} u}{o(N(P))}.$$

Now $o(N(P)) \mid o(G)$ since $N(P)$ is a subgroup of $G$, hence $p^{n+1}u \mid o(N(P))$ is an integer. Also, since $p^{n+1} \nmid o(G)$, $p^{\alpha+1}$ cannot divide $o(N(P))$. But then $p^{n+1}u \mid o(N(P))$ must be divisible by $p$, so we can write as $kp$, where $k$ is an integer. Feeding this information back into our equation above, we have

$$\frac{o(G)}{o(N(P))} = 1 + kp.$$

Recalling that $o(G) \mid o(N(P))$ is the number of p-Sylow subgroups in $G$, we have the theorem.

**Supplementary Reading**

1. Burnside W., Theory of groups of finite order, Cambridge Univ. Press 1911. New York.

2. Hall, Marshall, Theory of groups, New York, The Macmillan Company, 1961.

3. Mc. Kay James H, Another Proof of Cauchy's group Theorem, American Math. Monthly, Vol 66; 119 (1959).

4. Herstein, I.N., Topics in Algebra, New York, Blaisdell (1964).

5. John B. Fraleigh, A First Course in Abstract Algebra, Addison Wesley, 1967.

# CHAPTER FOUR
# SMARANDACHE SEMIGROUPS

Padilla Raul introduced the notion of Smarandache semigroups in the year 1998 in the paper entitled Smarandache Algebraic Structures. Since groups are the perfect structures under a single closed associative binary operation, it has become infeasible to define Smarandache groups. Smarandache semigroups are the analog in the Smarandache ideologies of the groups.

Now in this chapter we define new classes of Smarandache semigroups like Smarandache Lagrange semigroups, Smarandache p-Sylow subgroups, Smarandache subsemigroups, Smarandache hyper subsemigroups Smarandache simple semigroups and Smarandache Cauchy semigroup. Finally, the concept of Smarandache cosets was introduced in 2001 in the paper Smarandache cosets which has appeared in the online "Smarandache Notions Journal and is accessible at:
http://www.gallup.unm.edu/~smarandache/Cosets.pdf

## 4.1 Definition of Smarandache Semigroup

Here we first recall the definition of Smarandache semigroups as given by Raul (1998) and introduce in this section concepts like Smarandache commutative semigroup, Smarandache weakly commutative semigroup, Smarandache cyclic and weakly cyclic semigroups.

**DEFINITION:** *The Smarandache semigroup (S-semigroup) is defined to be a semigroup A such that a proper subset of A is a group (with respect to the same induced operation).*

**DEFINITION:** *Let S be a S-semigroup. If every proper subset of A in S, which is a group is commutative then we say the S-semigroup S to be a Smarandache commutative semigroup.*

*Remark:* It is important to note that if S is a commutative semigroup and if S is a S-semigroup then S is a Smarandache commutative semigroup. Here we are interested in finding whether there exists proper subsets of S-semigroups which are subgroups of which some of them are commutative and some non-commutative. This leads us to define:

**DEFINITION:** *Let S be S-semigroup, if S contains at least a proper subset A that is a commutative subgroup under the operations of S then we say S is a Smarandache weakly commutative semigroup.*

**DEFINITION:** *Let S be S-semigroup if every proper subset A of S which is a subgroup is cyclic then we say S is a Smarandache cyclic semigroup.*

**DEFINITION:** *Let S be a S-semigroup if there exists at least a proper subset A of S, which is a cyclic subgroup under the operations of S then we say S is a Smarandache weakly cyclic semigroup.*

**DEFINITION:** *Let S be a S-semigroup. If the number of distinct elements in S is finite, we say S is a finite S-semigroup otherwise we say S is a infinite S-semigroup.*

We are more interested in this book only about S-semigroups of finite order. We use the term subgroup or group in a S-semigroup in a synonymous way

### 4.2 Examples of S-semigroups

The lucidity and understanding of an algebraic concept is made easy only when it is illustrated by many examples. So this book tries to give many examples of S-semigroups.

**Example 4.2.1** Let $Z_{12}$ = {0, 1, 2, ... , 9, 10, 11} be the semigroup under multiplication modulo 12. Clearly, $Z_{12}$ is a S-semigroup. The subsets, which form the subgroups under multiplication mod 12, are given by the following tables:

| × | 1 | 5 |
|---|---|---|
| 1 | 1 | 5 |
| 5 | 5 | 1 |

*1 is the multiplicative identity*

| × | 9 | 3 |
|---|---|---|
| 9 | 9 | 3 |
| 3 | 3 | 9 |

*9 is the multiplicative identity*

| × | 4 | 8 |
|---|---|---|
| 4 | 4 | 8 |
| 8 | 8 | 4 |

*4 is the multiplicative identity*

| × | 1 | 7 |
|---|---|---|
| 1 | 1 | 7 |
| 7 | 7 | 1 |

*1 is the multiplicative identity*

| × | 1 | 11 |
|---|---|----|
| 1 | 1 | 11 |
| 11| 11| 1  |

*1 is the multiplicative identity*

| × | 1  | 5  | 7  | 11 |
|---|----|----|----|----|
| 1 | 1  | 5  | 7  | 11 |
| 5 | 5  | 1  | 11 | 7  |
| 7 | 7  | 11 | 1  | 5  |
| 11| 11 | 7  | 5  | 1  |

*1 is the multiplicative identity*

Thus we see $Z_{12}$ is only a Smarandache weakly cyclic semigroup as $Z_{12}$ has 6 proper subsets which are subgroups under multiplication modulo 12. Here of the six subgroups 5 are cyclic subgroups of order 2 and one is a non-cyclic subgroup of order 4. Further it is very important and fascinating to note that every subgroup of $Z_{12}$ does not have the same element 1 as its multiplicative identity.

***Example 4.2.2:*** Let $Z_6 = \{0, 1, 2, 3, 4, 5\}$ is the semigroup under multiplication mod 6. Clearly, $Z_6$ is a S-semigroup having only the cyclic group of order 2 viz.

| × | 4 | 2 |
|---|---|---|
| 4 | 4 | 2 |
| 2 | 2 | 4 |

| × | 1 | 5 |
|---|---|---|
| 1 | 1 | 5 |
| 5 | 5 | 1 |

$Z_6$ is a Smarandache cyclic semigroup of order 6.

***Example 4.2.3:*** Let S(3) be the set of all maps from the three element set (1, 2, 3) to itself. Clearly, S(3) under the operations of composition of maps 'o' is a semigroup. Further S(3) is a S-semigroup. The subsets, which are subgroups of S(3), are given by the following tables using the following notation,

$$1 = \begin{pmatrix} 1 & 2 & 3 \\ 1 & 2 & 3 \end{pmatrix}, \; p_1 = \begin{pmatrix} 1 & 2 & 3 \\ 1 & 3 & 2 \end{pmatrix}, \; p_2 = \begin{pmatrix} 1 & 2 & 3 \\ 3 & 2 & 1 \end{pmatrix},$$

$$p_3 = \begin{pmatrix} 1 & 2 & 3 \\ 2 & 1 & 3 \end{pmatrix}, \; p_4 = \begin{pmatrix} 1 & 2 & 3 \\ 2 & 3 & 1 \end{pmatrix}, \; p_5 = \begin{pmatrix} 1 & 2 & 3 \\ 3 & 1 & 2 \end{pmatrix}.$$

| o   | 1   | $p_1$ |
|-----|-----|-------|
| 1   | 1   | $p_1$ |
| $p_1$ | $p_1$ | 1   |

| o   | 1   | $p_2$ |
|-----|-----|-------|
| 1   | 1   | $p_2$ |
| $p_2$ | $p_2$ | 1   |

| o   | 1   | $p_3$ |
|-----|-----|-------|
| 1   | 1   | $p_3$ |
| $p_3$ | $p_3$ | 1   |

| o | 1 | $p_4$ | $p_5$ |
|---|---|---|---|
| 1 | 1 | $p_4$ | $p_5$ |
| $p_4$ | $p_4$ | $p_5$ | 1 |
| $p_5$ | $p_5$ | 1 | $p_4$ |

and

| o | 1 | $p_1$ | $p_2$ | $p_3$ | $p_4$ | $p_5$ |
|---|---|---|---|---|---|---|
| 1 | 1 | $p_1$ | $p_2$ | $p_3$ | $p_4$ | $p_5$ |
| $p_1$ | $p_1$ | 1 | $p_5$ | $p_4$ | $p_3$ | $p_2$ |
| $p_2$ | $p_2$ | $p_4$ | 1 | $p_5$ | $p_1$ | $p_3$ |
| $p_3$ | $p_3$ | $p_5$ | $p_4$ | 1 | $p_2$ | $p_1$ |
| $p_4$ | $p_4$ | $p_2$ | $p_3$ | $p_1$ | $p_5$ | 1 |
| $p_5$ | $p_5$ | $p_3$ | $p_1$ | $p_2$ | 1 | $p_4$ |

Thus, we see S(3) has subgroup; all them are not cyclic, so S(3) is a Smarandache weakly cyclic semigroup. It is absorbing to note unlike in the example 4.2.1, here for every subgroup the identity element is the same viz. $1 = \begin{pmatrix} 1 & 2 & 3 \\ 1 & 2 & 3 \end{pmatrix}$.

***Example 4.2.4*** Let $Z_8 = \{0, 1, 2, 3, \ldots, 7\}$ be the semigroup under multiplication modulo 8. The semigroup $Z_8$ is a S-semigroup. It has the following subsets, which are subgroups given by the following tables:

| × | 1 | 3 |
|---|---|---|
| 1 | 1 | 3 |
| 3 | 3 | 1 |

| × | 1 | 5 |
|---|---|---|
| 1 | 1 | 5 |
| 5 | 5 | 1 |

| × | 1 | 7 |
|---|---|---|
| 1 | 1 | 7 |
| 7 | 7 | 1 |

| × | 1 | 3 | 5 | 7 |
|---|---|---|---|---|
| 1 | 1 | 3 | 5 | 7 |
| 3 | 3 | 1 | 7 | 5 |
| 5 | 5 | 7 | 1 | 3 |
| 7 | 7 | 5 | 3 | 1 |

Now $Z_8$ has 4 subsets, which are subgroups of which 3 are cyclic and one is not cyclic but abelian. Thus, $Z_8$ is a Smarandache abelian semigroup, which is not a Smarandache cyclic semigroup. Here it is pertinent to note 1 which is the unit of $Z_8$ acts as the unit for all subgroups of $Z_8$.

***Example 4.2.5*** Consider the semigroup $Z_9 = \{0, 1, 2, 3, ... , 8\}$ under multiplication modulo 9. Now $Z_9$ is a S-semigroup. The following subsets of $Z_9$, which are subgroups are given by the following tables:

| × | 1 | 8 |
|---|---|---|
| 1 | 1 | 8 |
| 8 | 8 | 1 |

| × | 1 | 2 | 4 | 5 | 7 | 8 |
|---|---|---|---|---|---|---|
| 1 | 1 | 2 | 4 | 5 | 7 | 8 |
| 2 | 2 | 4 | 8 | 1 | 5 | 7 |
| 4 | 4 | 8 | 7 | 2 | 1 | 5 |
| 5 | 5 | 1 | 2 | 7 | 8 | 4 |
| 7 | 7 | 5 | 1 | 8 | 4 | 2 |
| 8 | 8 | 7 | 5 | 4 | 2 | 1 |

This example is unique in its own way for it has only 2 proper subsets which are subgroups under multiplication; they are of order 2 and 6. Clearly the order of the S-semigroup is 9 and 2 does not divide 9 and 6 also does not divide 9. Here we cannot even say the order of the subgroup and that of the S-semigroup are relatively prime. For $(2, 9) = 1$ and $(6, 9) = 3$.

***Example 4.2.6:*** Let $Z_{25} = \{0, 1, 2, 3, ... , 23, 24\}$ be the semigroup under multiplication modulo 25. The subset of $Z_{25}$ which are subgroups of $Z_{25}$ are given by A = \{1, 24\} and B = \{1, 2, 3, 4, 6, 7, 8, 9, 11, 12, 13, 14, 16, 17, 18, 19, 21, 22, 23, 24\} that is $Z_{25} \setminus \{0, 5, 10, 15, 20\}$ = B is also a subgroup under multiplication modulo 25. Thus $o(A) = 2$ and $o(B) = 20$. Clearly, B is not a cyclic subgroup of $Z_{25}$.
Here $Z_{25}$ is not a Smarandache cyclic semigroup.

## 4.3 Some Preliminary Theorems

Here we give some basic theorems about some S-semigroups. This will not only make the S-semigroup properties easy and enjoyable but also make one see how these structures satisfy various new properties.

**LEMMA 4.3.1:** *Let S be a Smarandache cyclic semigroup. Then S is a Smarandache commutative semigroup.*

*Proof:* By the very definition, we know all cyclic groups are abelian. So a Smarandache cyclic semigroup is a Smarandache commutative semigroup. Hence the claim.

**THEOREM 4.3.2:** *Let G be a Smarandache commutative semigroup. G in general need not be a Smarandache cyclic semigroup.*

*Proof:* We prove this by a counter example. Consider the S-semigroup given by $Z_{25} = \{0, 1, 2, ... ,23, 24\}$, the semigroup under multiplication modulo 25 given in example 4.2.6. Clearly, $Z_{25}$ is a Smarandache abelian semigroup but $Z_{25}$ is not a Smarandache cyclic semigroup as $Z_{25} \setminus \{0, 5, 10, 15, 20\} = B$ is a proper subset of $Z_{25}$ which forms an abelian group under multiplication but is not a cyclic group. Hence the claim.

The examples introduced in these sections have enabled us to prove the following results, about the semigroup of integers $Z_n$ under multiplication modulo n and S(n) the semigroup of mappings of the set of n elements to itself.

**THEOREM 4.3.3:** *$Z_n$ be the semigroup under multiplication modulo n; $n \geq 3$; $Z_n$ is a S-semigroup.*

*Proof:* $Z_n = \{0, 1, 2, \ldots, n-1\}$. $Z_n$ is a semigroup under multiplication modulo n. Clearly, we have the set $A = \{1, n-1\}$ is proper subset of $Z_n$, which is a subgroup under multiplication given by the following table:

| × | 1 | n-1 |
|---|---|-----|
| 1 | 1 | n-1 |
| n-1 | n-1 | 1 |

Hence, $Z_n$ is a S-semigroup of order n.

From the above theorem, we have a nice property about the S-semigroup $Z_n$.

**THEOREM 4.3.4:** *The S-semigroup $Z_n$ is a Smarandache weakly cyclic semigroup*

*Proof:* By the above theorem, $Z_n$ always has a cyclic subgroup of order 2 given by $A = \{1, n-1\}$. Hence $Z_n$ is a Smarandache weakly cyclic semigroup.

**THEOREM 4.3.5:** *S(n) is the S-semigroup.*

*Proof:* Clearly S(n) is the semigroup of order $n^n$. $S_n$ is a S-semigroup for it contains the symmetric group of degree n, i.e $S_n$ is a proper subset which is the group of permutations on $(1, 2, 3, \ldots, n)$. Hence S(n) is a S-semigroup.

Clearly S(n) is a not Smarandache commutative semigroup. Further we have the following engrossing results about S(n).

**THEOREM 4.3.6:** *S(n) the S-semigroup is not a Smarandache commutative (abelian) semigroup ($n \geq 3$).*

*Proof:* Now the semigroup S(n) contains $S_n$ the symmetric group of degree n as a proper subset which is a group; further $S_n$ is non abelian so S(n) is not a Smarandache commutative semigroup.

*Remark:* The condition $n \geq 3$ is important for S(2) is abelian.

**COROLLARY:** *S(n) is a Smarandache weakly cyclic semigroup.*

*Proof:* S(n) is a S-semigroup for $S_n$ the proper subset having n! elements is a subgroup of S(n). Now consider any element

$$\begin{pmatrix} 1 & 2 & 3 & 4 & 5 & \ldots & n-1 & n \\ 2 & 3 & 4 & 5 & 6 & \ldots & n & 1 \end{pmatrix} = p.$$

Thus, p generates a cyclic group of order n, which is cyclic. Hence the claim.

**PROBLEM 1:** Find the subgroups in the semigroup $Z_{21}$ under multiplication modulo 21. Is $Z_{21}$ a Smarandache cyclic semigroup?

**PROBLEM 2:** Find S-semigroup $Z_n$ ($Z_n$ semigroup under multiplication modulo n) which is a Smarandache cyclic semigroup for n > 30.

**PROBLEM 3:** Find all subsets which are subgroups of $Z_{120}$, the S-semigroup under multiplication modulo 120.

**PROBLEM 4:** How many cyclic subgroups are there in $Z_{120}$? ($Z_{120}$ given in problem 3). Is $Z_{120}$ a Smarandache cyclic semigroup?

**PROBLEM 5:** Find all cyclic subgroups in the semigroup $Z_{20} \times Z_{15}$ where $Z_{20} \times Z_{15}$ is the Cartesian product of the semigroups $Z_{20}$ and $Z_{15}$ under multiplication. Prove $Z_{20} \times Z_{15}$ is a S-semigroup.

**PROBLEM 6:** Find all abelian subgroups in S(15); where S(15) is the S-semigroup of mappings of the set S =(1, 2, …, 15) to itself.

**PROBLEM 7:** Find all abelian groups which are not cyclic in S(15) (given in problem 6).

**PROBLEM 8:** Prove S(15) is only

1. Smarandache weakly cyclic semigroup.
2. Smarandache weakly abelian semigroup.

**PROBLEM 9:** Find the largest subgroup in $S(20) \times S(6)$ where S(20) and S(6) are S-semigroups of order $20^{20}$ and $6^6$ respectively.

**PROBLEM 10:** Prove the largest subgroup in the S-semigroup $Z_{49} = \{0, 1, 2, \ldots, 48\}$ under multiplication modulo 49 is of order 42.

**PROBLEM 11:** Let $S_{2\times 2} = \{(a_{ij})/ a_{ij} \in Z_5 = \{0, 1, 2, 3, 4\}\}$. Prove $S_{2\times 2}$ is not a Smarandache commutative semigroup under matrix multiplication. Find the order of $S_{2\times 2}$.

**PROBLEM 12:** Is $R_{2\times 2} = \{(a_{ij})/a_{ij} \in Z_4 = \{0, 1, 2, 3\}\}$, the semigroup under matrix multiplication a Smarandache weakly cyclic semigroup? Prove your claim.

### 4.4 Smarandache Subsemigroup

In this section, we introduce the concept of Smarandache subsemigroup and obtain some interesting results about this study and this analysis of Smarandache subsemigroup has lead us to the definition of Smarandache hyper subsemigroup which is defined in the following section 4.5.

**DEFINITION:** *Let S be a S-semigroup. A proper subset A of S is said to be a Smarandache subsemigroup of S if A itself is a S-semigroup, that is A is a semigroup of S containing a proper subset B such that B is the group under the operations of S. Note we do not accept A to be a group. A must only be a semigroup.*

*Example 4.4.1:* Let $S(5)$ be the set of all mappings of the set $S = (1, 2, 3, 4, 5)$ to itself. Clearly $S(5)$ is a semigroup under composition of mappings. $S(5)$ is a S-semigroup as the proper subset $S_5 \subset S(5)$ is a group, that is the symmetric group of degree 5. Now take A = the semigroup generated by the elements, viz.,

$$\begin{pmatrix} 1 & 2 & 3 & 4 & 5 \\ 1 & 1 & 1 & 1 & 1 \end{pmatrix}, \begin{pmatrix} 1 & 2 & 3 & 4 & 5 \\ 2 & 3 & 4 & 5 & 1 \end{pmatrix}$$

Clearly A is a proper subset of $S(5)$ which is a semigroup. A is a S-semigroup as A contains the subgroup B generated by $\begin{pmatrix} 1 & 2 & 3 & 4 & 5 \\ 2 & 3 & 4 & 5 & 1 \end{pmatrix}$.

Thus $S(5)$ has proper Smarandache subsemigroups.

*Example 4.4.2:* Let $Z_{10} = \{0, 1, 2, \ldots, 9\}$ be the semigroup under multiplication modulo 10. $Z_{10}$ is a S-semigroup. $Z_{10}$ also has Smarandache subsemigroup. For $A = \{0, 1, 9\} \subset Z_{10}$ is a Smarandache subsemigroup of $Z_{10}$.

We are now interested to study the following; that is the reverse structure of the Smarandache subsemigroup suppose we have say $A_1, \ldots, A_n$ to be n proper subsets of a S-semigroup S which are the only subgroups of S then our question now is; does for each $A_i$, $i = 1, 2, \ldots, n$ we have a subsemigroup $P_i \subset S$ such that $A_i$ is properly contained in each $P_i$. The answer to our question is this is not always true; which is substantiated by the following example.

*Example 4.4.3:* Let $Z_7 = \{0, 1, 2, \ldots, 6\}$ be the S-semigroup under multiplication mod 7. The only subsets which are subgroups of $Z_7$ are $A = \{1, 6\}$ and $A_2 = \{1, 2, 3, 4, 5, 6\}$. Clearly take $P_1 = \{0, 1, 6\}$ then $P_1$ is a subsemigroup containing the subset $A_1$, which is a subgroup. Thus $P_1$ is a Smarandache subsemigroup of $Z_7$. But for $A_2$ we cannot find a proper subsemigroup in $Z_7$ such that $A_2$ is a proper subset of it. So all subgroups of a S-semigroup may not in general be contained in a Smarandache subsemigroup.

Thus we like to see for what S-semigroups we have all semigroups to be contained in proper subsemigroups.

**THEOREM 4.4.1:** *Let $Z_n = \{0, 1, 2, \ldots, n-1\}$ be the S-semigroup under multiplication modulo n-1, n is a composite number. Then every proper subset of $Z_n$ which are subgroup of $Z_n$ is properly contained in a proper Smarandache subsemigroup.*

*Proof:* Given $Z_n = \{0, 1, 2, \ldots, n-1\}$ where n is a composite number and $Z_n$ is a S-semigroup. If $A_1, \ldots, A_m$ are proper subsets of $Z_n$ which are subgroups of $Z_n$. Clearly none of the semigroups $A_1, A_2, \ldots, A_m$ have n-1 elements in them. All subgroups have elements strictly less than n-1. Hence $A_i \cup \{0\} \subset Z_n$ for i = 1, 2, ... , m is a Smarandache subsemigroup of $Z_n$. Hence the claim.

**THEOREM 4.4.2:** *Let $Z_p$ = {0, 1, 2, 3, ... , p-1}, p a prime be the S-semigroup of order p. Clearly $Z_p$ has a proper subset which cannot be properly contained in a proper subsemigroup of $Z_p$.*

*Proof:* Given $Z_p$ = {0, 1, 2,.., p-1} is a S-semigroup of order p, p a prime. The subgroups of $Z_p$ are $A_1$ = {1, p-1} and $A_2$ = {1, 2, ... , p-1}. Clearly $A_1 \cup \{0\}$ is a subsemigroup of $Z_p$. Hence $A_1 \cup \{0\}$ is a Smarandache subsemigroup. But $A_2$ cannot be strictly contained in any proper subset of $Z_p$. So $A_2$ cannot be contained in a subsemigroup in $Z_p$. Hence the claim.

**DEFINITION:** *Let S be a S-semigroup. If $A \subset S$ is a proper subset of S and A is a subgroup which cannot be contained in any proper subsemigroup of S we say A is the largest subgroup of S.*

This terminology will be used throughout this book.

**THEOREM 4.4.3:** *Let S(n) be the symmetric semigroup of order $n^n$, where S = (1, 2, 3, ... , n) and S(n) is the set of all maps from S to S. Now $S_n$ is the largest subgroup in S(n).*

*Proof:* We know $S_n$ is the largest subgroup of S(n). For no element in S(n) \ $S_n$ has inverse for it to become a group. Hence the claim.

**DEFINITION:** *We call the semigroup S(n) as the Smarandache symmetric semigroup of order $n^n$.*

**COROLLARY 4.4.4:** *The Smarandache symmetric semigroup S(n) has its largest group $S_n$ to be contained in the proper subset A where A is a subsemigroup properly contained in S(n).*

*Proof:* We know from the above theorem $S_n$ is the largest group. Take

$$A = S_n \cup \left\{ \begin{pmatrix} 1 & 2 & 3 & \ldots & n \\ 1 & 1 & 1 & \ldots & 1 \end{pmatrix}, \begin{pmatrix} 1 & 2 & 3 & \ldots & n \\ 2 & 2 & 2 & \ldots & 2 \end{pmatrix}, \ldots, \begin{pmatrix} 1 & 2 & 3 & \ldots & n \\ n & n & n & \ldots & n \end{pmatrix} \right\}$$

Clearly A is a proper subset and a subsemigroup of S(n). Hence the claim.

## 4.5 Smarandache Hyper Subsemigroups

In this section we introduce a new concept called Smarandache hyper subsemigroups and define what are called Smarandache simple semigroup when they have no proper Smarandache hyper subsemigroups. There are S-semigroups, which do not have Smarandache hyper subsemigroups. We obtain nice results about these Smarandache hyper subsemigroups.

**DEFINITION:** *Let S be a S-semigroup. If A be a proper subset of S which is subsemigroup of S and A contains the largest group of S then we say A to be the Smarandache hyper subsemigroup of S.*

*Example 4.5.1:* Let $S(8)$ be the Smarandache symmetric semigroup of all mappings of the set $S = (1, 2, 3, \ldots, 8)$. Now $S(8)$ has a Smarandache hyper subsemigroup for take

$$S_8 \cup \left\{ \begin{pmatrix} 1 & 2 & 3 & \ldots & 8 \\ 1 & 1 & 1 & \ldots & 1 \end{pmatrix}, \right.$$

$$\left. \begin{pmatrix} 1 & 2 & 3 & \ldots & 8 \\ 2 & 2 & 2 & \ldots & 2 \end{pmatrix}, \begin{pmatrix} 1 & 2 & 3 & \ldots & 8 \\ 3 & 3 & 3 & \ldots & 3 \end{pmatrix}, \ldots, \begin{pmatrix} 1 & 2 & 3 & \ldots & 8 \\ 8 & 8 & 8 & \ldots & 8 \end{pmatrix} \right\}.$$

Clearly A is a subsemigroup of $S(8)$ and has the largest group $S_8$ in it.

It is interesting to know the relation between Smarandache hyper subsemigroup and the Smarandache subsemigroup.

**THEOREM 4.5.1:** *Let S be a S-semigroup. Every Smarandache hyper subsemigroup is a Smarandache subsemigroup but every Smarandache subsemigroup is not a Smarandache hyper subsemigroup.*

*Proof:* Given S is a S-semigroup, and $A \subset S$ is a Smarandache hyper subsemigroup of S. Clearly A contains the largest subgroup in S so A is a Smarandache subsemigroup of S.

Conversely, if S contains B to be Smarandache subsemigroup to show B in general is not a Smarandache hyper subsemigroup. We prove this only by a counter example. Let $Z_{16} = \{0, 1, 2, 3, \ldots, 15\}$ be the S-semigroup under multiplication modulo 16. Now $A = \{0, 1, 15\}$ is Smarandache subsemigroup which is clearly not a Smarandache hyper subsemigroup. Hence the claim.

**THEOREM 4.5.2:** *$S(n)$ the Smarandache symmetric semigroup of mappings of the set $S = \{1, 2, \ldots, n\}$ to itself has Smarandache hyper subsemigroup.*

*Proof:* The only largest subgroup of $S(n)$ is $S_n$ take $A =$

$$S_n \cup \left\{ \begin{pmatrix} 1 & 2 & 3 & \ldots & n \\ 1 & 1 & 1 & \ldots & 1 \end{pmatrix} \begin{pmatrix} 1 & 2 & 3 & \ldots & 8 \\ 2 & 2 & 2 & \ldots & 2 \end{pmatrix}, \begin{pmatrix} 1 & 2 & 3 & \ldots & 8 \\ 3 & 3 & 3 & \ldots & 3 \end{pmatrix} \right.$$

$$\left. , \ldots, \begin{pmatrix} 1 & 2 & 3 & \ldots & 8 \\ 8 & 8 & 8 & \ldots & 8 \end{pmatrix} \right\}.$$

Clearly A is a Smarandache hyper subsemigroup of $S(n)$ as it contains the largest subgroup $S_n$.

*Example 4.5.2:* Let $Z_{11} = \{0, 1, 2, \ldots, 10\}$ be the S-semigroup under multiplication of order 11. Now the largest subgroup of $Z_{11}$ is $A = \{1, 2, 3, \ldots, 10\}$. Clearly $A = Z_{11} \setminus$

{0} which cannot be contained in any proper subset of $Z_{11}$ containing A. Hence $Z_{11}$ has no Smarandache hyper subsemigroup.

This example leads us to define Smarandache simple semigroups.

**DEFINITION:** *Let S be a S-semigroup. We say S is a Smarandache simple semigroup if S has no proper subsemigroup A, which contains the largest subgroup of S or equivalently S has no Smarandache hyper subsemigroup.*

**THEOREM 4.5.3:** *$Z_p$ = {0, 1, 2, ... , p-1} where p is a prime is a S-semigroup under multiplication modulo p. But $Z_p$ is a Smarandache simple semigroup.*

*Proof:* $Z_p$ = {0, 1, 2, ... , p-1} is a S-semigroup. To show $Z_p$ is Smarandache simple semigroup we have to show the largest subgroup of $Z_p$ cannot be contained in a proper subsemigroup of $Z_p$.

Now the largest subgroup of $Z_p$ is A = {1, 2, 3, ... , p-1}. Clearly A = $Z_p$ \ {0} which cannot be strictly contained in a proper subset of $Z_p$ which is subsemigroup of $Z_p$ other than itself. Hence, $Z_p$ is simple.

## 4.6 Smarandache Lagrange Semigroup

In this section, we define the concept of Smarandache Lagrange semigroup and Smarandache weakly Lagrange semigroup. Further, as the classical Lagrange's theorem for groups is not true in case of S-semigroup, we have defined Smarandache Lagrange semigroup and Smarandache weakly Lagrange semigroup and prove Lagrange's theorem to be true for Smarandache Lagrange semigroups.

Also the converse of Lagrange's theorem is not true in case of S-semigroups and its subgroups.

**DEFINITION:** *Let S be a finite S-semigroup. If the order of every subgroup of S divides the order of the S-semigroup S then we say S is a Smarandache Lagrange semigroup.*

***Example 4.6.1:*** Let $Z_4$ = {0, 1, 2, 3} be the semigroup under multiplication modulo 4. Clearly, $Z_4$ is a S-semigroup. Further the only subgroup in $Z_4$ is A = {1,3} and o(A) / 4 so $Z_4$ is a Smarandache Lagrange semigroup.

However, we see in general most of the S-semigroups are not Smarandache Lagrange semigroups. For example, we consider $Z_9$.

***Example 4.6.2:*** Let $Z_9$ = {0, 1, 2, 3, 4, 5, 6, 7, 8} be the semigroup under multiplication mod 9. Now A = {1, 8} is a subgroup of $Z_9$. Also B = {1, 2, 4, 5, 7, 8} is also a subgroup of $Z_9$ the order of both of them do not divide 9. Hence, $Z_9$ is not a Smarandache Lagrange semigroup.

Therefore, we are interested in defining the concept of Smarandache weakly Lagrange semigroup.

**DEFINITION:** *Let S be a finite S-semigroup. If there exists at least one subgroup A that is a proper subset (A ⊂ S) having the same operations of S whose order divides the order of S then we say that S is a Smarandache weakly Lagrange semigroup.*

***Example 4.6.3:*** $Z_{10}$ = {0, 1,2, ... , 9} is a S-semigroup under multiplication modulo 10. Clearly A = {1, 9} is a subset of $Z_{10}$ which is a subgroup such that o(A) / 10. Therefore, $Z_{10}$ is a Smarandache weakly Lagrange semigroup.

It is still engaging and important to note that all S-semigroups are not Smarandache weakly Lagrange semigroup. For the S-semigroup given in example 4.6.2 is not even a Smarandache weakly Lagrange semigroup. This leads to formulate the following theorem.

**THEOREM 4.6.1:** *Every Smarandache Lagrange semigroup is a Smarandache weakly Lagrange semigroup and not conversely.*

*Proof:* By the very definition of Smarandache Lagrange semigroup and Smarandache weakly Lagrange semigroup we see every Smarandache Lagrange semigroup is a Smarandache weakly Lagrange semigroup. To prove the converse we consider the following example.

***Example 4.6.4:*** Let S(3) be the semigroup under the composition of mappings of the 3 element set S = (1, 2, 3). Clearly S(3) is a S-semigroup for it contains $S_3$ the symmetric group of degree 3 which is of order 6. Clearly 6 ∤ o(S(3)) as o(S(3)) = 27. So S(3) is not a Smarandache Lagrange semigroup, but S(3) is a Smarandache weakly Lagrange semigroup for consider the set

$$A = \left\{ \begin{pmatrix} 1 & 2 & 3 \\ 1 & 2 & 3 \end{pmatrix}, \begin{pmatrix} 1 & 2 & 3 \\ 2 & 3 & 1 \end{pmatrix}, \begin{pmatrix} 1 & 2 & 3 \\ 3 & 1 & 2 \end{pmatrix} \right\}$$

Clearly A is a proper subset of S(3) which is also a group under composition of maps that is a cyclic group of order 3. o(A) / S(3). Hence S(3) is a Smarandache weakly Lagrange semigroup.

Thus by the above example we see the converse of the theorem 4.6.1 is not true in general.

**THEOREM 4.6.2:** *Every Smarandache symmetric semigroup S(n) is a Smarandache weakly Lagrange semigroup and not a Smarandache Lagrange semigroup for n ≥ 3.*

*Proof:* Now to prove S(n) is a Smarandache weakly Lagrange semigroup and not a Smarandache Lagrange semigroup we have to prove

1. S(n) has a subgroup which divides the order of S(n) and
2. S(n) has a subgroup which does not divide the order of S(n).

To prove this theorem we consider two cases one when n is odd and other when n is even.

**CASE 1:**

Let n be odd. Now S(n) is a semigroup of order $n^n$ and S(n) is a S-semigroup, for S(n) contains $S_n$ the symmetric group of degree n. Further the order of the group $S_n$ is n!. Now to show S(n) is not a Smarandache Lagrange semigroup we see $o(S_n) \nmid o(S(n))$ for $o(S(n)) = n \times ... \times n$ where n is odd but $o(S_n) = n!$ so n-1 is even hence $o(S_n) \nmid S(n)$ as no even integer can divide $n^n$ when n is odd. Hence S(n) is not a Smarandache Lagrange semigroup.

To show S(n) is a Smarandache weakly Lagrange semigroup. For consider A = {the permutation of 1 2 3 4 ... n} = the subgroup generated by g where

$$g = \begin{pmatrix} 1 & 2 & 3 & 4 & 5 & ... & n-1 & n \\ 2 & 3 & 4 & 5 & 6 & ... & n & 1 \end{pmatrix}.$$

Clearly, the number of elements in the group generated by g is n, and $n/n^n$. So S(n) is a Smarandache weakly Lagrange semigroup when n is odd.

**CASE 2:**

Let n be even. Now S(n) is a S-semigroup of order $n^n$. Let $S_n$ be the symmetric group of degree n, so $o(S_n)$ is n!. $o(S_n)$ does not divide $o(S(n))$ for n is even so n-1 is odd; $n - 1 \nmid n^n$ when n is even; so S(n) is not a Smarandache Lagrange semigroup.

But S(n) is a Smarandache weakly Lagrange semigroup, for A = the group generated by g

$$g = \begin{pmatrix} 1 & 2 & 3 & 4 & 5 & ... & n-1 & n \\ 2 & 3 & 4 & 5 & 6 & ... & n & 1 \end{pmatrix}$$

is a cyclic group of order n. $n/n^n$. Hence the claim.

*Remark:* $n \geq 3$ is essential in the theorem for if n = 2 we have S(2) contains only 4 elements and $S_2$ the subgroup has 2 elements so $o(S_2) / o(S(2))$. Hence S(2) is a Smarandache Lagrange semigroup.

### 4.7 Smarandache p-Sylow Subgroups

In this section, we define the concept of Smarandache p-Sylow subgroups of a S-semigroups and get some appealing results about them.

**DEFINITION:** *Let S be a finite S-semigroup. p be a prime such that p divides the order of S. If there exists a subgroup A of S of order p or $p^t$ (t > 1) we say S has a Smarandache p-Sylow subgroup.*

*Note:* It is very important to note that p/o(S) but $p^t \nmid o(S)$ still we may have Smarandache p-Sylow subgroups having $p^t$ elements; this is illustrated by the following example.

**Example 4.7.1:** Let $Z_{10} = \{0, 1, 2, \ldots, 9\}$ be the semigroup under multiplication modulo 10. Clearly, $Z_{10}$ is a S-semigroup. For $A = \{1, 9\} \subset Z_{10}$ is a subgroup of $Z_{10}$. Now 2/10 but $2^2 \nmid 10$ but $Z_{10}$ has Smarandache 2-Sylow subgroups of order 4. Now take $B = \{6, 2, 4, 8\} \subset Z_{10}$. B is a subgroup under multiplication modulo 10 given by the following table

| × | 6 | 2 | 4 | 8 |
|---|---|---|---|---|
| 6 | 6 | 2 | 4 | 8 |
| 2 | 2 | 4 | 8 | 6 |
| 4 | 4 | 8 | 6 | 2 |
| 8 | 8 | 6 | 2 | 4 |

B is a subgroup with 6 as the multiplicative identity. Further $C = \{1, 3, 7, 9\}$ is also a Smarandache 2- Sylow subgroup of $Z_{10}$ given by the following table

| × | 1 | 3 | 7 | 9 |
|---|---|---|---|---|
| 1 | 1 | 3 | 7 | 9 |
| 3 | 3 | 9 | 1 | 7 |
| 7 | 7 | 1 | 9 | 3 |
| 9 | 9 | 7 | 3 | 1 |

*Thus unlike as in the classical definition of p-Sylow subgroups we see the definition of Smarandache p-Sylow subgroups is different for if p is a prime dividing the order of the S-semigroup S, then S can have Smarandache p-Sylow subgroups of order $p^\alpha$ where $p^\alpha \nmid o(S)$. From above example it is evident for 2/10 and $4 \nmid 10$.*

**Example 4.7.2:** Let $Z_{16} = \{0, 1, 2, \ldots, 15\}$ be a S-semigroup of order 16. Clearly, this S-semigroup has 2-Sylow semigroups of order 2, 4 and 8 given by the following tables

| ×  | 1  | 15 |
|----|----|----|
| 1  | 1  | 15 |
| 15 | 15 | 1  |

*1 is the identity element of this subgroup A = {1, 15}.*

| ×  | 1  | 3  | 9  | 11 |
|----|----|----|----|----|
| 1  | 1  | 3  | 9  | 11 |
| 3  | 3  | 9  | 11 | 1  |
| 9  | 9  | 11 | 1  | 3  |
| 11 | 11 | 1  | 3  | 9  |

*1 is the identity element of this subgroup A = {1, 3, 9, 11}.*

| × | 1 | 5 | 9 | 13 |
|---|---|---|---|----|
| 1 | 1 | 5 | 9 | 13 |
| 5 | 5 | 9 | 13 | 1 |
| 9 | 9 | 13 | 1 | 5 |
| 13 | 13 | 1 | 5 | 9 |

*1 is the identity element of this subgroup B = {1, 5, 9, 13}.*

$A_1 = \{1, 3, 5, 7, 9, 11, 13, 15\}$ can be verified to be a 2-Sylow subgroup of order 8 in $Z_{10}$. Thus unlike our classical group the S-semigroup of order 16 has 2 Sylow subgroups of order 2, 4 and 8 respectively.

***Example 4.7.3:*** Let $Z_{12} = \{0, 1, 2, 3, \ldots, 11\}$ be the S-semigroup. The operation on $Z_{12}$ is multiplication modulo 12. The p-Sylow subgroups of $Z_{12}$ are for p = 2, given by $A_1, \ldots, A_6$ where $A_1 = \{1, 7\}$, $A_2 = \{1, 5\}$, $A_3 = \{1, 11\}$, $A_4 = \{4, 8\}$, $A_5 = \{3, 9\}$ and $A_6 = \{1, 5, 7, 11\}$ are groups of order 2. Finally $A_6 = \{1, 5, 7, 11\}$ is a subgroup of order 4 given by the following table

| × | 1 | 5 | 7 | 11 |
|---|---|---|---|----|
| 1 | 1 | 5 | 7 | 11 |
| 5 | 5 | 1 | 11 | 7 |
| 7 | 7 | 11 | 1 | 5 |
| 11 | 11 | 7 | 5 | 1 |

In fact, this S-semigroup is also a Smarandache Lagrange semigroup.

## 4.8 Smarandache Cauchy Element of a S-semigroup

In this section, we define the concept called Smarandache Cauchy element in a S-semigroup. The main motivation for this is to try to prove or disprove the classical Cauchy's theorem in case of S-semigroup. This is carried out in chapter 5.

**DEFINITION:** *Let S be a finite S-semigroup. An element $a \in A$, $A \subset S$ (A a proper subset of S and A is the subgroup under the operation of S) is said to be a Smarandache Cauchy element of S if $a^r = 1$ (r > 1) and 1 is the unit element of A and r divides the order of S otherwise a is not a Smarandache Cauchy element of S.*

***Example 4.8.1:*** Let S(3) be the semigroup got from the mappings of the set S = (1, 2, 3) to itself. Every element of S(3) is not a Smarandache Cauchy element of S(3).

For consider the element

$$a = \begin{pmatrix} 1 & 2 & 3 \\ 1 & 3 & 2 \end{pmatrix} \in A = \left\{ \begin{pmatrix} 1 & 2 & 3 \\ 1 & 2 & 3 \end{pmatrix}, \begin{pmatrix} 1 & 2 & 3 \\ 1 & 3 & 2 \end{pmatrix} \right\} \subset S(3).$$

Clearly, A is a subgroup. Now, $a^2 = \begin{pmatrix} 1 & 2 & 3 \\ 1 & 2 & 3 \end{pmatrix}$ = identity element of S(3). But $o(S(3)) = 27$ and $a^2 = 1$ Clearly $2 \nmid 27$. Hence $a = \begin{pmatrix} 1 & 2 & 3 \\ 1 & 3 & 2 \end{pmatrix}$ is not a Smarandache Cauchy element of S(3).

Thus, we see all elements of the form $a^m = 1$ (m>1) in a S-semigroup S need not in general be Smarandache Cauchy elements of S.

**DEFINITION:** *Let S be a finite S-semigroup if every element in every subgroup of S is a Smarandache Cauchy element of S then we say S is a Smarandache Cauchy semigroup.*

But it is interesting to note that there exists S-semigroup in which no element is a Smarandache Cauchy element.

**THEOREM 4.8.1:** *Let $Z_P = \{0, 1, 2, \ldots, p-1\}$, p is a prime be the S-semigroup under multiplication. No element in $Z_p$ is a Smarandache Cauchy element of $Z_p$.*

*Proof:* Consider $Z_p = \{0, 1, 2, \ldots, p-1\}$, $Z_p$ is a S-semigroup under multiplication modulo p. $A = \{1, 2, \ldots, p-1\}$ is a proper subset of $Z_p$ which is a group under multiplication. Every element in A is such that $a^r = 1$ where r lies between 2 and p-1 as p is a prime with every element less than itself cannot divide p, hence no element in $Z_p$ is a Smarandache Cauchy element of $Z_p$ as p is a prime. Hence the claim.

## 4.9 Smarandache Coset

This section is devoted to introduction of Smarandache right coset (left coset) in a S-semigroup A. We prove by examples as in the case of groups the number of elements in each coset for a given subgroup is not equal in general.

**DEFINITION:** *Let A be a S-semigroup. H is a proper subset of A ($H \subset A$) be a group under the operations of A. For any $a \in A$ define the Smarandache right coset $Ha = \{ha \;/\; h \in H\}$, Ha is called the Smarandache right coset of H in A.*

Similarly Smarandache left coset of H in A can be defined. If the S-semigroup is a commutative semigroup then we will see the concept of the Smarandache left coset and the Smarandache right coset will coincide and will get the same set.

**DEFINITION:** *Let S be a S-semigroup. $H \subset S$ be a subgroup of S. We say aH is the Smarandache coset of H in S for $a \in S$ if Ha = aH, that is $\{ha \;/\; h \in H\} = \{ah \;/\; h \in H\}$.*

**Example 4.9.1:** Let $Z_{12} = \{0, 1, 2, \ldots, 11\}$ be the S-semigroup under multiplication modulo 12. Clearly, $Z_{12}$ is a Smarandache commutative semigroup. Let $A = \{3, 9\}$, is given by the following table:

| × | 9 | 3 |
|---|---|---|
| 9 | 9 | 3 |
| 3 | 3 | 9 |

Clearly A is subgroup with 9 as the identity element as $9^2 \equiv 9 \pmod{12}$. For $4 \in Z_{12}$ the right (left) coset of A is $4A = \{0\}$. Thus we see the number of elements in $4A$ is not 2 but one viz $\{0\}$. Further for $5 \in Z_{12}$ we have $5A = \{3, 9\} = A$. This property too is unique and distinctly different from cosets in a group. For in case of cosets of H in a group G we see that $aH = H$ if and only if $a \in H$ but in case of Smarandache cosets we can have $aH = H$ even if $a \notin H$ but a is in the semigroup S.

*Example 4.9.2:* Let $Z_p = \{0, 1, 2, \ldots, p-1\}$ where p is a prime be the S-semigroup under multiplication modulo 12. Now the subgroups of $Z_p$ are $A = \{1, 2, \ldots, p-1\}$ and $\{1, p-1\}$. So the cosets of A are $\{0\}$ and A.

Thus we have interesting results about Smarandache cosets which will be proved in the next chapter.

**PROBLEM 1:** Is $Z_{18} = \{0, 1, 2, \ldots, 17\}$ a Smarandache Lagrange semigroup?

**PROBLEM 2:** Give an example of a Smarandache weakly commutative semigroup other than S(n).

**PROBLEM 3:** Is $R_{3\times 3} = \{(a_{ij}) \mid a_{ij} \in Z_2 = \{0, 1\}\}$ the semigroup under matrix multiplication a Smarandache weakly Lagrange semigroup?

**PROBLEM 4:** Find the number of elements in $R_{3\times 3}$ given in the problem 3.

**PROBLEM 5:** Let $A_{3\times 3} = \{(a_{ij}) \mid a_{ij} \in Z_2 = \{0, 1\}\}$ denote only upper triangular matrices with the possibility of the diagonal elements having the value to be zero. $A_{3\times 3}$ under matrix multiplication is a S-semigroup.

a. Is $A_{3\times 3}$ a Smarandache Lagrange semigroup?
b. Find the number of elements in $A_{3\times 3}$.

**PROBLEM 6:** Let $Z_{35} = \{0, 1, 2, \ldots, 34\}$ be the S-semigroup under multiplication modulo 35. Does $Z_{35}$ have Smarandache 5-Sylow semigroup or a Smarandache 7-Sylow semigroup?

**PROBLEM 7:** Let $Z_{20} = \{0, 1, 2, \ldots, 19\}$ be the S-semigroup under multiplication modulo 20. Find all subsets in $Z_{20}$ which are subgroups under multiplication.

**PROBLEM 8:** Find all Smarandache p-Sylow subgroups in $Z_{125} = \{0, 1, 2, \ldots, 124\}$ where $Z_{125}$ is a S-semigroup under multiplication mod 125.

**PROBLEM 9:** Find all the Smarandache p-Sylow subgroups of $R_{3\times 3}$ given in problem 3.

**PROBLEM 10:** Find all the Smarandache p-Sylow subgroups of S(4).

**PROBLEM 11:** Find all the Smarandache p-Sylow subgroups of S(25).

**PROBLEM 12:** When will the number of Smarandache p-Sylow subgroups be more in S(n) when n is prime or when n is a composite number?

**PROBLEM 13:** Find two subgroups of order two in $Z_{33} = \{0, 1, 2, \ldots, 32\}$ the S-semigroup under multiplication mod 33.

**PROBLEM 14:** Find all subsets in $Z_{33}$, the S-semigroup of order 33, that are subgroups under multiplication modulo 33.

**PROBLEM 15:** Find all the Smarandache Cauchy elements in $Z_{20}$.

**PROBLEM 16:** Give an example of S-semigroup in which every invertible element is a Smarandache Cauchy element.

**PROBLEM 17:** Give an example of a S-semigroup of order 11 in which no element is a Smarandache Cauchy element.

**PROBLEM 18:** Does there exist a S-semigroup of order 24 in which no element is a Smarandache Cauchy element?

**PROBLEM 19:** Does there exist an example of a S-semigroup of order 30 in which every element is a Smarandache Cauchy element?

**PROBLEM 20:** Find all the Smarandache Cauchy elements in $S_{3\times 3} = \{(a_{ij}) \mid a_{ij} \in Z_4 = \{0, 1, 2, 3\}\}$ the S-semigroup under matrix multiplication.

**PROBLEM 21:** Find whether $S_{3\times 3}$ (given in Problem 20) is Smarandache weakly Lagrange semigroup.

**PROBLEM 22:** Find at least three Smarandache Cauchy elements in $S_{3\times 3}$. ($S_{3\times 3}$ given in Problem 20.)

**PROBLEM 23:** Does $S_{3\times 3}$ (given in Problem 20) contain a subset of order 24, which is a subgroup of $S_{3\times 3}$ under multiplication?

**PROBLEM 24:** Let S(5) be the Smarandache symmetric semigroup. Does S(5) contain a subgroup of order 12?

**PROBLEM 25:** Let S(7) be the Smarandache symmetric semigroup. Find a subgroup in S(7) of order 120.

**PROBLEM 26:** Find a subgroup of order 60 in S(5). (S(5) given in Problem 24).

**PROBLEM 27:** Find all the subgroups in S(15).

**PROBLEM 28:** Prove at least $S(3) \times S(4)$ is a S-semigroup and a Smarandache weakly Lagrange semigroup.

**PROBLEM 29:** Find all the Smarandache p-Sylow subgroups in $S(3) \times S(4)$, the direct product of the semigroups which is a S-semigroup.

**PROBLEM 30:** Find all the Smarandache p-Sylow subgroups in $Z_7 \times Z_9$ and prove $Z_7 \times Z_9$ is a S-semigroup.

**PROBLEM 31:** Find all the Smarandache Cauchy elements in $Z_5 \times Z_{12}$ and prove $Z_5 \times Z_{12}$ is a Smarandache weakly Lagrange semigroup.

**PROBLEM 32:** Find all the Smarandache Cauchy elements in $S(5) \times S(6)$.

**PROBLEM 33:** Find all Smarandache p-Sylow subgroups of $Z_7 \times Z_8$. Does $Z_7 \times Z_8$ contain Smarandache Cauchy elements? If so find them.

**PROBLEM 34:** Let $S(4)$ be the S-semigroup and $A_4$ is the subgroup of $S(4)$. Are the right cosets of $A_4$ the same as the left cosets of $A_4$? Justify your answer.

**PROBLEM 35:** Let $Z_{120} = \{0, 1, 2, \ldots, 119\}$ be the S-semigroup under multiplication modulo 120. Does there exist one element x in the subgroup H in $Z_{120}$ such that $o(Hx) = o(H)$ for $x \in Z_{120} \setminus \{0\}$?

**PROBLEM 36:** Does there exist a subgroup H in $Z_{42} = \{0, 1, 2, \ldots, 41\}$ such that $o(aH) = o(H)$ for all $a \in Z_{42}$ and $a \neq 0$?

**PROBLEM 37:** Let $Z_p = \{0, 1, \ldots, p-1\}$ be the S-semigroup of order p, p is a prime prove the subgroup $A = \{1, p-1\}$ partitions $Z_p$ into equivalence classes.

**PROBLEM 38:** Let $Z_n = \{0, 1, \ldots, n-1\}$ be the S-semigroup of order n, n a composite number. Prove when does the subgroup $A = \{1, p-1\}$ partition $Z_n$ into equivalence classes. Find how many distinct equivalence classes exist?

**PROBLEM 39:** Find the Smarandache right and left cosets for the subgroup $A \times B \subset Z_9 \times Z_{21}$ where $A = \{1, 8\}$ and $B = \{1, 20\}$ when $z = (3, 7)$ and $y = (6, 3)$. Is $z(A \times B) = (A \times B)z$?

**PROBLEM 40:** Find the Smarandache cosets for the subgroup $A_1 \times B_1 \subset Z_9 \times Z_{21}$ when $A_1 = Z_9$ and $B_1 = \{1, 20\}$ for $x = (3, 7)$ and $y = (6, 3)$ Make a comparison about Smarandache cosets in example 39 and 40.

**PROBLEM 41:** Prove for $Z_{19} = \{0, 1, 2, \ldots, 18\}$ the S-semigroup under multiplication, the Smarandache cosets got by $A = \{1, 18\}$ is such that $Z_{19} = \bigcup_{x_i \in Z_{19}} x_i A$ ; $x_i A \cap x_j A = \phi$ if $i \neq j$.

**PROBLEM 42:** Can you extend the concept in problem 41 and prove $Z_n = \bigcup_{x_i \in Z_n} x_i A$ where $A = \{1, n-1\}$ and $x_i A \cap x_j A = \phi$ if $i \neq j$ for any n?

**PROBLEM 43:** Let $S(5)$ be the Smarandache symmetric semigroup for the group $A_5$ and $S_5$, find the Smarandache coset decomposition of $S(5)$ relative to $A_5$ and $S_5$.

**PROBLEM 44:** Find for $x = \begin{pmatrix} 1 & 2 & 3 & 4 & 5 \\ 1 & 2 & 1 & 3 & 5 \end{pmatrix}$ $xA_5$ and $A_5x$. $S(5) \supset A_5$ and $x \in S(5)$ given in Problem 43.

**PROBLEM 45:** For the same problem 43 does there exist a $x \in S(5) \setminus S_5$ such that, $xA_5 = A_5x$?

**PROBLEM 46:** Let $S(6)$ be the S-semigroup. $S_6$ the subgroup in $S(6)$. Find an element $x$ in $S(6)$ and $x \notin S_6$ such that $xS_6 = S_6x$.

**PROBLEM 47:** Does the Smarandache coset decomposition of $S(6)$ by $S_6$ divide $S(6)$ into equivalence classes of same length? Prove your answer.

**PROBLEM 48:** Does there exist any subgroup G in $S(6)$ which will decompose $S(6)$ into equivalence classes?

**PROBLEM 49:** Let $S(9)$ be the S-semigroup G be the group generated by $g = \begin{pmatrix} 1 & 2 & 3 & \ldots & 8 & 9 \\ 2 & 3 & 4 & \ldots & 9 & 1 \end{pmatrix}$. Clearly, G is a subgroup of order 9. Can $xG = Gx$ for any $x \in S(9) \setminus G$?

**PROBLEM 50:** In problem 49 will G divide - $S(9)$ into equivalence classes?

### Supplementary Reading

1. Padilla Raul, Smarandache algebraic structures, Bull. of Pure and Applied Sciences, Delhi, Vol. 17E, No. 1, 119-121, 1998.

2. W.B.Vasantha Kandasamy, Smarandache cosets, Smarandache Notions Journal, American Research Press, 2001. Internet address: http://www.gallup.unm.edu/~smarandache/Loops.pdf

# CHAPTER FIVE

# THEOREMS FOR S-SEMIGROUPS

In this chapter, we introduce the classical theorems for S-semigroups. Already in Chapter 4 all essential definitions were made about S-semigroup, their properties studied and those concepts illustrated with examples. So now, we proceed to introduce those classical theorems of groups for S-semigroups.

## 5.1 Lagrange's Theorem for S-semigroups

Lagrange's Theorem states that if G is a finite group having a subgroup H then o(H)/o(G). However for S-semigroups we see Lagrange's Theorem does not hold in general. So we introduced in Chapter 4 Smarandache Lagrange semigroup and we will prove in this section the Lagrange's Theorem holds good only for this class of Smarandache Lagrange semigroups.

**THEOREM 5.1.1:** *Let S be a finite Smarandache Lagrange semigroup. If H is a subgroup of S then o(H) / o(S).*

*Proof:* By the very definition of the Smarandache Lagrange semigroup, we know every subgroup of S divides the order of S. Hence the result.

Thus, it is very important to note that only the class of Smarandache Lagrange semigroups satisfies the Smarandache Lagrange's Theorem all other semigroups which are not Smarandache Lagrange semigroups will not even satisfy the classical Lagrange's theorem for groups which is stated in chapter 3 of this book.

Hence, the class of Smarandache weakly Lagrange semigroups do not satisfy the classical theorem. So in general when we go for the S-semigroup structure we see that no more the order of the subgroups of the S-semigroup divide the order of the S-semigroup.

Now the converse of the Lagrange's theorem for S-semigroups is doublefold. One, order of the subgroup of a finite S-semigroup does not divide the order of the S-semigroup, for which the class of Smarandache weakly Lagrange semigroup is an evident example. Secondly, to disprove the converse of the Smarandache Lagrange's theorem for Smarandache Lagrange semigroup we have to find a finite Smarandache Lagrange semigroup for which we should find a divisor and show that there does not exist a subgroup of that order.

**THEOREM 5.1.2:** *Let S be a finite Smarandache Lagrange semigroup. If m/o(S), S need not in general contain a subgroup of order m.*

*Proof:* Consider the S-semigroup given by $Z_{12}$ = {0, 1, 2, ... , 11}. Clearly, $Z_{12}$ is a Smarandache Lagrange semigroup of order 12. The only subgroups of $Z_{12}$ are $A_1$ = {1, 11}, $A_2$ = {1, 5}, $A_3$ = {1, 7}, $A_4$ = {4, 8 / 4 is the multiplicative identity} $A_5$ = {3, 9 / 9 is the multiplicative identity} and $A_6$ = {1, 5, 7, 11}. Now every subgroup is of order 2 or 4 and 2/12 and 4/12. However, 3/12 and $Z_{12}$ does not contain a subgroup of order 3. Hence the claim.

Therefore, we now give the weak from of the Lagrange's theorem for S-semigroups which is as follows:

**THEOREM 5.1.3:** *S be a finite Smarandache weakly Lagrange's semigroup then there exists at least one subgroup A contained in S such that o(A)/o(S).*

*Proof:* The theorem is true only for finite Smarandache weakly Lagrange semigroups. By the very definition of the Smarandache weakly Lagrange semigroup we have the claim in the theorem to be true.

Now we have to consider only a S-semigroup, which is not a Smarandache weakly Lagrange semigroup to prove the following theorem.

**THEOREM 5.1.4:** *Let S be a finite S-semigroup. In general, order of a subgroup need not divide the o(S).*

*Proof:* By a counter example. Clearly, this finite semigroup S must not be a Smarandache Lagrange semigroup or a Smarandache weakly Lagrange semigroup. Consider $Z_9$ = {0, 1, 2, ... , 8} be the semigroup of order 9. Clearly, $Z_9$ is a S-semigroup under multiplication modulo 9. The subgroups in $Z_9$ are {1, 8} = A and {1, 2, 4, 5, 7, 8} = B. Clearly o(A) ∤ o($Z_9$) and o(B) ∤ o($Z_9$). Hence the above theorem.

Thus we see in case of finite S-semigroup S we have subgroups in S none of whose order divides the order of S, these S-semigroups are really very innovative because they make the analog of the classical Lagrange's theorem false. This forces us or instigates us to define a new class of S-semigroups which we choose to call the Smarandache non Lagrange semigroups.

**DEFINITION:** *Let S be a finite S-semigroup we say S is a Smarandache non-Lagrange semigroup if none of the subgroups of S divide the order of the S-semigroup S.*

The class of Smarandache non Lagrange semigroups is non empty, by the following theorem.

**THEOREM 5.1.5:** *Let $Z_p$ = {0, 1, 2, ... , p-1}, p an odd prime. $Z_p$ is a S-semigroup under multiplication modulo p. $Z_p$ for every prime p, is a Smarandache non Lagrange semigroup.*

*Proof:* $Z_P$ = {0, 1, 2, ... , p-1}, p is a prime and $Z_p$ is a S-semigroup under multiplication modulo p. The only subgroups of $Z_p$ are A = {1, p-1} and B = {1, 2, 3, ... , p-1} = $Z_p$ \ {0} . Clearly, 2 ∤ p and p–1 ∤ p. Since in $Z_p$ every element other than 1 and 0 generates $Z_p$ under multiplication. Hence the claim.

Since the number of primes is infinite, we have infinitely many Smarandache non-Lagrange semigroups. Interestingly we have divided the class of all S-semigroups as Smarandache Lagrange semigroups, Smarandache weakly Lagrange semigroups and Smarandache non Lagrange semigroups. Class of Smarandache Lagrange semigroups is properly contained in the class of Smarandache weakly Lagrange semigroups and Smarandache non Lagrange semigroups are disjoint with the class of Smarandache weakly Lagrange semigroups.

## 5.2 Cayley's Theorem for S-semigroups

To prove the classical Cayley's theorem of group for S-semigroups we need the concept of S-semigroup homomorphism using which we will prove our result.

**DEFINITION:** *Let S and S' be any two S-semigroups. A map $\phi$ from S to S' is said to be a S-semigroup homomorphism if $\phi$ restricted to a subgroup $A \subset S \to A' \subset S'$ is a group homomorphism. The S-semigroup homomorphism is an isomorphism if $\phi: A \to A'$ is one to one and onto. Similarly, one can define S-semigroup automorphism on S.*

It is surprising to note that two S-semigroups S and S' can be isomorphic even if $o(S) \neq o(S')$.

***Example 5.2.1:*** Let $Z_{12} = \{0, 1, 2, \ldots, 11\}$ be a S-semigroup under multiplication modulo 12, $Z_7 = \{0, 1, 2, \ldots, 6\}$ be the S-semigroup under multiplication. Define $\phi: Z_{12} \to Z_7$ by

$$\phi(1) = 1$$
$$\phi(11) = 6 \text{ and}$$
$$\phi(0) = 0$$
$$\phi(x) = 0 \text{ for all } x \neq 1 \text{ and } 11.$$

Clearly $\phi: Z_{12} \to Z_7$ is a S-semigroup homomorphism for $\phi$ restricted to A where $A = \{1, 11\} \subset Z_{12} \to A' = \{1, 6\} \subset Z_7$ that is $\phi: A \to A'$ is an isomorphism of the subgroups. Using this definition of S-semigroup homomorphism we have now the analog of the Cayley's theorem, which is as follows.

**THEOREM 5.2.1: (CAYLEY'S THEOREM FOR S-SEMIGROUP)** *Every S-semigroup is isomorphic to a S-semigroup S(N); of mappings of a set N to itself, for some appropriate set N.*

*Proof:* Let S be a S-semigroup. That is A the proper subset of S which is a group under the operations of S. That is $\phi \neq A \subset S$. Now let N be any set, S(N) denote the set of all mappings from N to N. Clearly S(N) is a S-semigroup. We have in fact proved in the chapter 5 S(N) is a Smarandache weakly Lagrange semigroup of order $N^N$.

Now we use the classical theorem of Cayley for groups. By the classical Cayley's theorem for groups we can always find an isomorphism from the group A to a subgroup $S_N \subset S(N)$ for an appropriate N. Thus S is Smarandache homomorphic with S(N) for an appropriate N, that is A is isomorphic to a subgroup in $S_N \subset S(N)$. Hence the theorem.

## 5.3 Cauchy's Theorem for S-semigroups

Here we prove that Cauchy's theorem in general is not true for S-semigroups. In order to make possible an analog for Cauchy's theorem we have introduced in chapter 4 the concept of Cauchy element in a S-semigroup "An element x in A ⊂ S where S is a finite S-semigroup and A a subgroup in S is said to be a Smarandache Cauchy element of S if $x^r = 1$ (r > 1) and 1 is the unit element of A and r/o(S)". The Smarandache Cauchy semigroup is defined as a S-semigroup in which every element in every subgroup is a Smarandache Cauchy element of S.

**THEOREM 5.3.1:** *Let S be a finite Smarandache Cauchy semigroup. If $a \in S$ and $a^m = 1$ then m/o(S).*

*Proof:* The above theorem is true for all finite Smarandache Cauchy semigroup as we have by the very definition of Smarandache Cauchy semigroup the order of every invertible element $a \in S$ divides the order of S.

**THEOREM 5.3.2:** *Let S be any finite S-semigroup. If $a \in S$ and $a^m = 1$ then in general $m \nmid o(S)$.*

*Proof:* We prove this by giving an example. Let $S = Z_{11} = \{0, 1, 2, 3, \ldots, 10\}$ be the S-semigroup under multiplication modulo 11. Now we have baring 0 and 1. $2^{10} = 1$ where $2 \in A = \{1, 2, 3, 4, 5, 6, 7, 8, 9, 10\} \subset Z_{11}$ but $10 \nmid 11$. Hence the theorem.

## 5.4 p-Sylow Theorem for S-semigroups

In this chapter, we see how far the analogs of the three classical Sylow's theorems for groups can be given equivalent formulation for the S-semigroups. In chapter 4 we have already defined the concept of Smarandache p-Sylow semigroup.

**THEOREM 5.4.1 (FIRST SYLOW'S THEOREM FOR SMARANDACHE p-SYLOW SUBGROUP)** *Let S be a finite S-semigroup if p is a prime such that p/o(S), it does not imply S has a Smarandache p-Sylow subgroup.*

*Proof:* The proof of the theorem is by giving an example. Consider $Z_{12} = \{0, 1, 2, 3, \ldots, 11\}$; clearly $Z_{12}$ is a S-semigroup, 2/12 and 3/12. It is verified that $Z_{12}$ has no subgroup of order 3. Hence the claim.

This property leads us to define a new concept called Smarandache p-Sylow semigroup.

**DEFINITION:** *Let S be a finite S-semigroup. If for every prime p dividing o(S) we have a Smarandache p-Sylow subgroup then we say S is a Smarandache p-Sylow semigroup.*

**THEOREM 5.4.2:** *Let S be a finite Smarandache p-Sylow semigroup if p/o(S) (p any prime) then there exist a Smarandache p-Sylow subgroup of order p or $p^\alpha$.*

*Proof:* The result follows from the very definition of the Smarandache p-Sylow semigroup S.

*Example 5.4.1:* Let $S = (1, 2, 3, 4, 5)$ be the set with 5 elements $S(5)$ is the S-semigroup of order $5^5$. Clearly 5 is the only prime which divides $o(S(5))$. We have

$$A = \left\{ \begin{pmatrix} 1 & 2 & 3 & 4 & 5 \\ 2 & 3 & 4 & 5 & 1 \end{pmatrix} \right\}$$

i.e. the cyclic group generated by

$$\begin{pmatrix} 1 & 2 & 3 & 4 & 5 \\ 2 & 3 & 4 & 5 & 1 \end{pmatrix}$$

is a Smarandache p-Sylow subgroup of $S(5)$. Thus we have $S(5)$ to be a Smarandache p-Sylow semigroup. Using this result we prove $S(n)$ for the set of n elements $(1, 2, \ldots, n)$ is a Smarandache p-Sylow semigroup.

**THEOREM 5.4.3:** *The S-semigroup $S(n)$ is a Smarandache p-Sylow semigroup.*

*Proof:* Now let $p_1, p_2, \ldots, p_r$ be the primes which divide $o(S(n)) = n^n$; that is in short $p_1, p_2, \ldots, p_r$ are the distinct primes which divides n.(if n is a prime n = p) Now $S(n)$ is nothing but mappings of the set $S = (1, 2, \ldots, n)$ to itself. Further $S_n \subset S(n)$ where $S_n$ is the symmetric group of degree n.

Now for any prime $p_i$ ($p_i < n$, $p_i / n$) we have a permutation g which is such that

$$g = \begin{pmatrix} 1 & 2 & \ldots & p_{i-1} & p_1 & p_{i+1} & \ldots & n \\ 2 & 3 & \ldots & p_i & 1 & p_{i+1} & \ldots & n \end{pmatrix}$$

that is fixes $p_{i+1}, p_{i+2}, \ldots, n$ and translates each r to $r + 1$, $r = 1, 2, \ldots, p_{i-1}$ and $p_i = 1$. Now 'g' generates a cyclic group of order $p_i$ this is true for $i = 1, 2, \ldots, r$. Hence the claim.

**THEOREM 5.4.4: (SECOND PART OF SYLOW'S THEOREM).** *Let S be a S-semigroup. Two Smarandache p-Sylow subgroups in S need not be conjugate.*

*Proof:* The proof is given by the following example. Consider $Z_8 = \{0, 1, 2, \ldots, 7\}$ the S-semigroup of order 8 under multiplication modulo 8. The Smarandache 2-Sylow subgroups of $Z_8$ are $A = \{1, 7\}$, $B = \{1, 5\}$, $C = \{1, 3\}$ and $D = \{1, 3, 5, 7\}$. Clearly A is conjugate to B and C but D is not conjugate to A or B or C. Hence the claim.

To overcome this problem we leave it for the reader to introduce some more new concepts.

**THEOREM 5.4.5 (THIRD PART OF SYLOW'S THEOREM):** *Let S be a finite S-semigroup. If $p/o(S)$ and suppose S has Smarandache p-Sylow subgroup, then in general*

$$\frac{o(S)}{o(N(P))} \neq 1 + kp.$$

*Proof:* For this result, also we prove by an example. Now consider the S-semigroup $Z_8 = \{0, 1, 2, \ldots, 7\}$ this has Smarandache 2-Sylow subgroups of order 2 and 4 so we cannot find k which is such that $8/2 = 1 + 2.k$ that is $4 = 1 + 2k$ so k has no integer value. This is the case when the Smarandache 2-Sylow subgroup $A = \{1, 7\}$. Now $B = \{1, 3, 7, 5\}$ $8/4 = 1 + 2k = 2$, $k = ½$ which is not an integer. Thus the third part of p-Sylow theorem is not true in case of Smarandache p-Sylow subgroups.

Still we are interested in studying the situation when S is a finite S-semigroup and p a prime number with $p < o(S)$, $p \nmid o(S)$, we may still find a subgroup of order p or $p^\alpha$ ($\alpha > 1$), to characterize or make a note of such happening which is never possible in case of groups we give the following definition.

**DEFINITION:** *Let S be a finite S-semigroup. If for a prime p, $p < o(S)$, $p \nmid o(S)$ there exist a subgroup in S of order p then we call that subgroup a Smarandache non-p-Sylow subgroup of the S-semigroup S.*

We have in plenty such Smarandache non-p-Sylow subgroup for example.

*Example 5.4.2:* Let $Z_{23} = \{0, 1, 2, \ldots, 22\}$ is the S-semigroup of order 23. This has subgroup of order 2 (where 2 is an even prime) given by $A = \{1, 22\}$. Thus $Z_{23}$ has Smarandache non 2 Sylow subgroup.

**THEOREM 5.4.6:** *Let $Z_m = \{0, 1, \ldots, m-1\}$ be the S-semigroup of order m where m is an odd number. $Z_m$ has Smarandache non 2-Sylow subgroup.*

*Proof:* Clearly, we have the set $A = \{1, m-1\}$ to be a subgroup of order 2 in $Z_m$ and $2 \nmid m$ as m is odd. Hence the claim.

*Example 5.4.3:* Let $S(7)$ be the S-semigroup of the symmetric semigroup. Clearly $o(S(7)) = 7^7$ but $S(7)$ has Smarandache non p-Sylow subgroups for $p = 2, 3$ and $5$. For

$$A = \left\{ \begin{pmatrix} 1 & 2 & 3 & 4 & 5 & 6 & 7 \\ 1 & 2 & 3 & 4 & 5 & 6 & 7 \end{pmatrix}, \begin{pmatrix} 1 & 2 & 3 & 4 & 5 & 6 & 7 \\ 2 & 1 & 3 & 4 & 5 & 6 & 7 \end{pmatrix} \right\}$$

is a group of order 2 or a Smarandache non 2-Sylow subgroup of $S(7)$. Also

$$B = \left\{ \begin{pmatrix} 1 & 2 & 3 & 4 & 5 & 6 & 7 \\ 1 & 2 & 3 & 4 & 5 & 6 & 7 \end{pmatrix}, \begin{pmatrix} 1 & 2 & 3 & 4 & 5 & 6 & 7 \\ 2 & 3 & 4 & 1 & 5 & 6 & 7 \end{pmatrix}, \begin{pmatrix} 1 & 2 & 3 & 4 & 5 & 6 & 7 \\ 3 & 4 & 1 & 2 & 5 & 6 & 7 \end{pmatrix}, \begin{pmatrix} 1 & 2 & 3 & 4 & 5 & 6 & 7 \\ 4 & 1 & 2 & 3 & 5 & 6 & 7 \end{pmatrix} \right\}$$

is a Smarandache non 2-Sylow subgroup of order 4 in $S(7)$. Now let

$$C = \left\{ \begin{pmatrix} 1 & 2 & 3 & 4 & 5 & 6 & 7 \\ 1 & 2 & 3 & 4 & 5 & 6 & 7 \end{pmatrix}, \begin{pmatrix} 1 & 2 & 3 & 4 & 5 & 6 & 7 \\ 2 & 3 & 1 & 4 & 5 & 6 & 7 \end{pmatrix}, \\ \begin{pmatrix} 1 & 2 & 3 & 4 & 5 & 6 & 7 \\ 3 & 1 & 2 & 4 & 5 & 6 & 7 \end{pmatrix} \right\}$$

is a cyclic subgroup of order 3. So C is a Smarandache non 3-Sylow subgroup of S(7). Finally, S(7) has a cyclic group of order 5 given by

$$D = \left\{ \begin{pmatrix} 1 & 2 & 3 & 4 & 5 & 6 & 7 \\ 1 & 2 & 3 & 4 & 5 & 6 & 7 \end{pmatrix}, \begin{pmatrix} 1 & 2 & 3 & 4 & 5 & 6 & 7 \\ 2 & 3 & 4 & 5 & 1 & 6 & 7 \end{pmatrix}, \\ \begin{pmatrix} 1 & 2 & 3 & 4 & 5 & 6 & 7 \\ 3 & 4 & 5 & 1 & 2 & 6 & 7 \end{pmatrix}, \begin{pmatrix} 1 & 2 & 3 & 4 & 5 & 6 & 7 \\ 4 & 5 & 1 & 2 & 3 & 6 & 7 \end{pmatrix}, \\ \begin{pmatrix} 1 & 2 & 3 & 4 & 5 & 6 & 7 \\ 5 & 1 & 2 & 3 & 4 & 6 & 7 \end{pmatrix} \right\}$$

hence this is a Smarandache non 5-Sylow subgroup of the S-semigroup S(7).

This leads us to find a nice theorem on Smarandache symmetric semigroup S(n) n a positive integer.

**THEOREM 5.4.4:** *Let S(n) be the Smarandache symmetric semigroup, n a prime. Then S(n) has Smarandache non p-Sylow subgroups for all primes p; p < n.*

*Proof:* We know the primes less then n are 2, 3, 5, 7, 11,... p < n. Our claim is for every p(p < n) we have correspondingly a cyclic group of order p. This is got by permuting p elements cyclically in the set (1, 2, 3,…,p-1, p, p+1, ... , n) as follows.

The subgroup A generated by

$$g = \begin{pmatrix} 1 & 2 & \ldots & p-1 & p & p+1 & \ldots & n \\ 2 & 3 & \ldots & p & 1 & p+2 & \ldots & n \end{pmatrix}.$$

Clearly, A is a cyclic group of order p generated by g. Since the choice of p is arbitrary from the set of primes p < n our claim is true for all primes p, p < n, but $p \nmid o(S(n)) = n^n$ as n is a prime. Hence we have for every prime p, p < n (n a prime) there exist a Smarandache non p-Sylow subgroup in S(n).

### 5.5 Smarandache Cosets

This section is completely devoted in proving that in case of Smarandache cosets we do not have one to one correspondence between any two Smarandache right cosets of A in a S-semigroup. Further we prove that in general the Smarandache right cosets of any subgroup $A \subset S$ does not partition S.

**THEOREM 5.5.1:** *Let S be a S-semigroup. A $\subset$ S be a proper subset which is a group under the operations of S. There does not exist in general a one to one correspondence between any two Smarandache right cosets of A in the S-semigroup S.*

*Proof:* We prove the theorem by an example. Let $S = Z_{10} = \{0, 1, 2, \ldots, 9\}$ be the S-semigroup of order 10 under multiplication modulo 10.

Take $A = \{1, 9\} \subset Z_{10}$. Clearly, A is a subgroup of $Z_{10}$. We see $0A = \{0\}$, $3A = \{3, 7\}$, $5A = \{5\}$. So we cannot imagine of any one to one correspondence between right cosets of A in the S-semigroup $Z_{10}$. Similarly we take $B = \{2, 4, 6, 8\}$ which is a subgroup with 6 as the identity we get $5A_2 = \{0\}$ $3A_2 = A_2$. So there does not exist a one to one correspondence between the right cosets of B in $Z_{10}$. Hence the claim.

**THEOREM 5.5.2:** *The Smarandache right cosets of A in a S-semigroup does not in general partition S into either equivalence classes of same order or does not partition S at all.*

*Proof:* Consider the S-semigroup $Z_{10}$ given in theorem 5.5.1: where $A = \{1, 9\}$ and $B = \{6, 2, 4, 8\}$ are the subgroups of $S = Z_{10}$. The equivalence classes corresponding to $A = \{1, 9\}$ are $\{0\}$, $\{5\}$, $\{1, 9\}$, $\{2, 8\}$, $\{3, 7\}$ and $\{4, 6\}$. So A partitions $Z_{10}$ but not into equivalence classes of same length.

Now $B = \{6, 2, 4, 8\}$ is a subgroup of $Z_{10}$. Now the Smarandache coset division by B gives just $\{0\}$ and $\{6, 2, 4, 8\}$ only. Therefore, subsets do not even account for every element in $Z_{10}$. Hence the claim.

Finally we conclude this chapter with an interesting example study and test the validity of all the results proved in this chapter.

***Example 5.5.1:*** Let $S_{2 \times 2} = \left\{ \begin{pmatrix} 0 & 0 \\ 0 & 0 \end{pmatrix}, \begin{pmatrix} 0 & 1 \\ 0 & 0 \end{pmatrix}, \begin{pmatrix} 1 & 0 \\ 0 & 0 \end{pmatrix}, \begin{pmatrix} 0 & 0 \\ 1 & 0 \end{pmatrix}, \begin{pmatrix} 0 & 0 \\ 0 & 1 \end{pmatrix}, \begin{pmatrix} 1 & 1 \\ 0 & 0 \end{pmatrix}, \begin{pmatrix} 0 & 0 \\ 1 & 1 \end{pmatrix}, \begin{pmatrix} 1 & 0 \\ 1 & 0 \end{pmatrix}, \begin{pmatrix} 0 & 1 \\ 0 & 1 \end{pmatrix}, \begin{pmatrix} 1 & 0 \\ 0 & 1 \end{pmatrix}, \begin{pmatrix} 0 & 1 \\ 1 & 0 \end{pmatrix}, \begin{pmatrix} 1 & 1 \\ 1 & 0 \end{pmatrix}, \begin{pmatrix} 1 & 1 \\ 0 & 1 \end{pmatrix}, \begin{pmatrix} 0 & 1 \\ 1 & 1 \end{pmatrix}, \begin{pmatrix} 1 & 0 \\ 1 & 1 \end{pmatrix}, \begin{pmatrix} 1 & 1 \\ 1 & 1 \end{pmatrix} \right\}$

be the collection of all 2×2 matrices with entries from the prime field of characteristic 2 viz $Z_2 = \{0, 1\}$.

Now clearly the number of elements in $S_{2 \times 2}$ is $16 = 2^4$. Now first we show $S_{2 \times 2}$ is not a Smarandache Lagrange semigroup. To prove this we have to find a subgroup A of $S_{2 \times 2}$ of order n where n $\nmid$ 16.

Now $A = \left\{ \begin{pmatrix} 1 & 0 \\ 0 & 1 \end{pmatrix}, \begin{pmatrix} 0 & 1 \\ 1 & 1 \end{pmatrix}, \begin{pmatrix} 1 & 1 \\ 1 & 0 \end{pmatrix} \right\}$

is subgroup of $S_{2 \times 2}$ of order 3. Hence the claim. However, $S_{2 \times 2}$ is a Smarandache weakly Lagrange semigroup. To prove this we have to show $S_{2 \times 2}$ contains a subgroup B of order m where m/16.

To this end, we have
$$B = \left\{ \begin{pmatrix} 1 & 0 \\ 0 & 1 \end{pmatrix}, \begin{pmatrix} 0 & 1 \\ 1 & 0 \end{pmatrix} \right\}$$

which is a subgroup of order 2. Hence, $S_{2\times 2}$ is a Smarandache weakly Lagrange semigroup. To show $S_{2\times 2}$ is not a Smarandache Cauchy semigroup we have to get an element g in $S_{2\times 2}$ such that

$$g^n = \begin{pmatrix} 1 & 0 \\ 0 & 1 \end{pmatrix} \text{ but } n \nmid 16.$$

Consider
$$\begin{pmatrix} 1 & 1 \\ 1 & 0 \end{pmatrix} \in S_{2\times 2}$$

Clearly
$$\begin{pmatrix} 1 & 1 \\ 1 & 0 \end{pmatrix}^3 = \begin{pmatrix} 1 & 0 \\ 0 & 1 \end{pmatrix}$$

and $3 \nmid 16$. But $S_{2\times 2}$ has Smarandache Cauchy elements also as $\begin{pmatrix} 1 & 1 \\ 0 & 1 \end{pmatrix} \in S_{2\times 2}$ is such that $\begin{pmatrix} 1 & 1 \\ 0 & 1 \end{pmatrix}^2 = \begin{pmatrix} 1 & 0 \\ 0 & 1 \end{pmatrix}$. Hence $\begin{pmatrix} 1 & 1 \\ 0 & 1 \end{pmatrix}$ is a Cauchy element of $S_{2\times 2}$ as $2 \mid 16$.

Now to show $S_{2\times 2}$ contains both Smarandache p-Sylow subgroups and Smarandache non p-Sylow subgroups.

Clearly,
$$A = \left\{ \begin{pmatrix} 1 & 1 \\ 0 & 1 \end{pmatrix}, \begin{pmatrix} 1 & 0 \\ 0 & 1 \end{pmatrix} \right\}$$

is a subgroup of order 2 so A is Smarandache 2-Sylow subgroup of $S_{2\times 2}$. Consider the subgroup

$$C = \left\{ \begin{pmatrix} 1 & 0 \\ 0 & 1 \end{pmatrix}, \begin{pmatrix} 0 & 1 \\ 1 & 1 \end{pmatrix}, \begin{pmatrix} 1 & 1 \\ 1 & 0 \end{pmatrix} \right\}$$

C is a subgroup of order 3 in $S_{2\times 2}$ $3 \nmid 16$. So C is a Smarandache non 3-Sylow subgroup of $S_{2\times 2}$.

Finally, to show that the right cosets partition $S_{2\times 2}$ for the group C that is the following equivalence classes.

To prove this, consider $C = \left\{ \begin{pmatrix} 1 & 0 \\ 0 & 1 \end{pmatrix}, \begin{pmatrix} 0 & 1 \\ 1 & 1 \end{pmatrix}, \begin{pmatrix} 1 & 1 \\ 1 & 0 \end{pmatrix} \right\}$.

Now $x = \begin{pmatrix} 1 & 0 \\ 0 & 0 \end{pmatrix} \in S_{2\times 2}$, $xC = \left\{ \begin{pmatrix} 1 & 0 \\ 0 & 0 \end{pmatrix}, \begin{pmatrix} 0 & 1 \\ 0 & 0 \end{pmatrix}, \begin{pmatrix} 1 & 1 \\ 0 & 0 \end{pmatrix} \right\}$.

Now take $y = \begin{pmatrix} 0 & 1 \\ 0 & 0 \end{pmatrix}$ we get $yC = \left\{ \begin{pmatrix} 0 & 1 \\ 0 & 0 \end{pmatrix}, \begin{pmatrix} 1 & 1 \\ 0 & 0 \end{pmatrix}, \begin{pmatrix} 1 & 0 \\ 0 & 0 \end{pmatrix} \right\}$.

Next,

For $z = \begin{pmatrix} 1 & 1 \\ 0 & 0 \end{pmatrix}$ we have $zC = \left\{ \begin{pmatrix} 0 & 1 \\ 0 & 0 \end{pmatrix}, \begin{pmatrix} 1 & 1 \\ 0 & 0 \end{pmatrix}, \begin{pmatrix} 1 & 0 \\ 0 & 0 \end{pmatrix} \right\}$.

For $u = \begin{pmatrix} 0 & 1 \\ 1 & 0 \end{pmatrix}$ we have $uC = \left\{ \begin{pmatrix} 0 & 1 \\ 1 & 0 \end{pmatrix}, \begin{pmatrix} 1 & 1 \\ 0 & 1 \end{pmatrix}, \begin{pmatrix} 1 & 0 \\ 1 & 1 \end{pmatrix} \right\}$.

For $v = \begin{pmatrix} 0 & 0 \\ 1 & 0 \end{pmatrix}$ we have $vC = \left\{ \begin{pmatrix} 0 & 0 \\ 1 & 0 \end{pmatrix}, \begin{pmatrix} 0 & 0 \\ 1 & 1 \end{pmatrix}, \begin{pmatrix} 0 & 0 \\ 0 & 1 \end{pmatrix} \right\}$.

For $w = \begin{pmatrix} 1 & 0 \\ 1 & 0 \end{pmatrix}$ we have $wC = \left\{ \begin{pmatrix} 1 & 0 \\ 1 & 0 \end{pmatrix}, \begin{pmatrix} 0 & 1 \\ 0 & 1 \end{pmatrix}, \begin{pmatrix} 1 & 1 \\ 1 & 1 \end{pmatrix} \right\}$.

Thus $\begin{pmatrix} 0 & 0 \\ 0 & 0 \end{pmatrix} C = \begin{pmatrix} 0 & 0 \\ 0 & 0 \end{pmatrix}$ leads to the partition of $S_{2\times 2}$ baring $\begin{pmatrix} 0 & 0 \\ 0 & 0 \end{pmatrix}$ each of same length 3 and all the sets are disjoint.

Now the subgroup $D = \left\{ \begin{pmatrix} 1 & 0 \\ 0 & 1 \end{pmatrix}, \begin{pmatrix} 0 & 1 \\ 1 & 0 \end{pmatrix} \right\}$ divides $S_{2\times 2}$ into equivalence classes given by

$\left\{ \begin{pmatrix} 0 & 0 \\ 0 & 0 \end{pmatrix} \right\}, \left\{ \begin{pmatrix} 1 & 1 \\ 1 & 1 \end{pmatrix} \right\}, \left\{ \begin{pmatrix} 1 & 1 \\ 0 & 0 \end{pmatrix} \right\}, \left\{ \begin{pmatrix} 0 & 0 \\ 1 & 1 \end{pmatrix} \right\},$

$\left\{ \begin{pmatrix} 0 & 1 \\ 0 & 1 \end{pmatrix}, \begin{pmatrix} 1 & 0 \\ 1 & 0 \end{pmatrix} \right\}, \left\{ \begin{pmatrix} 1 & 0 \\ 0 & 0 \end{pmatrix}, \begin{pmatrix} 0 & 1 \\ 0 & 0 \end{pmatrix} \right\}, \left\{ \begin{pmatrix} 0 & 0 \\ 1 & 0 \end{pmatrix}, \begin{pmatrix} 0 & 0 \\ 0 & 1 \end{pmatrix} \right\},$

$\left\{ \begin{pmatrix} 1 & 1 \\ 0 & 1 \end{pmatrix}, \begin{pmatrix} 1 & 1 \\ 1 & 0 \end{pmatrix} \right\}, \left\{ \begin{pmatrix} 1 & 0 \\ 1 & 1 \end{pmatrix}, \begin{pmatrix} 0 & 1 \\ 1 & 1 \end{pmatrix} \right\}, \left\{ \begin{pmatrix} 1 & 0 \\ 0 & 1 \end{pmatrix}, \begin{pmatrix} 0 & 1 \\ 1 & 0 \end{pmatrix} \right\}.$

This is a partition unlike in group each class is of different lengths.

**PROBLEM 1:** Is $Z_{42} = \{0, 1, 2, \ldots, 41\}$ a Smarandache weakly semigroup or a Smarandache Lagrange semigroup?

**PROBLEM 2:** Is $Z_{30} = \{0, 1, 2, \ldots, 29\}$ a Smarandache Lagrange semigroup?

**PROBLEM 3:** Does there exist an example of a Smarandache Lagrange semigroup of order 14?

**PROBLEM 4:** Is $Z_{122} = \{0, 1, 2, \ldots, 121\}$ a Smarandache weakly Lagrange semigroup?

**PROBLEM 5:** Prove $Z_{2n} = \{0, 1, 2, \ldots, 2n-1\}$ is a Smarandache weakly Lagrange semigroup.

**PROBLEM 6:** Can $Z_{15}$ be a Smarandache Lagrange semigroup? Justify your answer.

**PROBLEM 7:** Prove $P_{2\times 2} = \{(a_{ij}) \mid a_{ij} \in Z_3 = \{0, 1, 2\}\}$ the semigroup under matrix multiplication is a Smarandache weakly Lagrange semigroup.

**PROBLEM 8:** Can $2\times 2$ matrices with entries from $Z_p$, p a prime, under matrix multiplication be Smarandache Lagrange semigroup?

**PROBLEM 9:** Does there exist a Smarandache p-Sylow semigroup in a S-semigroup of order 30?

**PROBLEM 10:** For the primes 2, 3, 5, 7 find a Smarandache p-Sylow semigroup of order m where 2/m, 3/m, 5/m and 7/m.

**PROBLEM 11:** Give an example of a Smarandache Cauchy semigroup of order 24.

**PROBLEM 12:** Does there exist a Smarandache Cauchy semigroup of order 210?

**PROBLEM 13:** Give an example of S-semigroup for which there is no Smarandache p-Sylow subgroups.

**PROBLEM 14:** Find all the Smarandache p-Sylow subgroups of S(20).

**PROBLEM 15:** How many Smarandache p-Sylow subgroups are there in S(12)?

**PROBLEM 16:** Find all Smarandache 3-Sylow subgroups of S(18)?

**PROBLEM 17:** Give an example of Smarandache Cauchy semigroup of order 20.

**PROBLEM 18:** Can there exist a Smarandache Cauchy semigroup of order 127? Justify your answer.

**PROBLEM 19:** Does there exist a Smarandache p-Sylow subgroup of order 37? Justify your answer.

**PROBLEM 20:** Find a Smarandache Cauchy semigroup of order 81.

**PROBLEM 21:** Give an example of a Smarandache non p-Sylow subgroup of order 18.

**PROBLEM 22:** Give an example of a Smarandache non p-Sylow subgroup order 72.

**PROBLEM 23:** Verify the classical Smarandache Cauchy theorem for $Z_{105}$.

**PROBLEM 24:** Verify the classical Smarandache Cauchy theorem for the group S(27).

**PROBLEM 25:** Verify Smarandache Sylow theorems for the S-semigroup $S(3) \times S(8)$.

**PROBLEM 26:** Divide into Smarandache cosets the S-semigroup $Z_{121}$ by any subgroup in $Z_{121}$.

**PROBLEM 27:** Divide $S(3) \times S(7)$ into Smarandache cosets by the subgroups $S_3 \times A_7$, $A_3 \times S_7$.

**PROBLEM 28:** Let $S_{2\times 2}$ the set of all $2\times 2$ matrices over the ring $Z_4 = \{0, 1, 2, 3\}$ under multiplication.

1. Prove $S_{2\times 2}$ is not Smarandache Lagrange semigroup.
2. Prove $S_{2\times 2}$ is Smarandache weakly Lagrange semigroup.
3. Prove $S_{2\times 2}$ has Smarandache p-Sylow subgroup.
4. Find in $S_{2\times 2}$ a Smarandache non p-Sylow subgroup.
5. Find a Smarandache Cauchy element in $S_{2\times 2}$.
6. An element which is not Smarandache Cauchy element in $S_{2\times 2}$.
7. Find a subgroup A of order 3 and find the Smarandache coset decomposition of $S_{2\times 2}$.

**PROBLEM 29:** Find all Smarandache non p-Sylow subgroups and Smarandache p-Sylow subgroups of $Z_{125}$.

**PROBLEM 30:** Find all Smarandache non p-Sylow subgroup of $S(6)$.

**PROBLEM 31:** Find all Smarandache p-Sylow subgroup of $Z_{120}$.

## Supplementary Reading

1. Herstein, I.N., Topics in Algebra, New York, Blaisdell (1964).

2. Padilla Raul, Smarandache algebraic structures, Bull of Pure and applied Sciences, Delhi, Vol. 17E, No. 1, 119-121, 1998.

3. W.B.Vasantha Kandasamy, Smarandache cosets, Smarandache Notions Journal, American Research Press, 2001. Internet address: http://www.gallup.unm.edu/~smarandache/Cosets.pdf

4. W.B.Vasantha Kandasamy, Smarandache loops, Smarandache Notions Journal, American Research Press, 2001. Internet address: http://www.gallup.unm.edu/~smarandache/Loops.pdf

# CHAPTER SIX

# SMARANDACHE NOTIONS IN GROUPS

In the earlier chapters we studied the S-semigroup and obtained some attractive results about them. In that study the algebraic structure under consideration was a S-semigroup, that is a semigroup, which has a proper subset in it which is a subgroup under the operations of the semigroup. Now in this chapter our algebraic object will be a group or a S-semigroup and in the group we will introduce Smarandache notions like Smarandache inverse of an element in a group and Smarandache conjugate elements in a group. Smarandache double coset, Smarandache normal subgroup, Smarandache quotient group and Smarandache direct products and study them.

### 6.1 Smarandache Inverse in Groups

This section is completely devoted to the introduction of Smarandache inverse in a group and we obtain some interesting properties about them. Generally in a group every element has a unique inverse but in case of groups Smarandache inverse may or may not exists. Certain elements of a group may have Smarandache inverse and some may not have Smarandache inverse.

**DEFINITION:** *Let G be a group. An element $x \in G \setminus \{1\}$ is said to have a Smarandache inverse y in G if xy = 1 and for a, b $\in G \setminus \{1, x, y\}$ we have xa = y (or ax = y), yb = x (or by = x) with ab = 1.*

The pair (x, y) is called the Smarandache inverse pair and (a, b) is called the related or relating pair of the Smarandache inverse pair (x, y). Further if x is the Smarandache inverse of y we can equivalently say y is the Smarandache inverse of x and (y, x) is the Smarandache inverse pair and the related pair remains the same viz (a, b). Thus x is the Smarandache inverse of y or y is the Smarandache inverse of x when (x, y) is the Smarandache inverse pair no distinction between these two statements exist as in the case of inverse elements in a group G; for x is the inverse of y is the same as y is the inverse of x.

**DEFINITION:** *Let G be a group. x be the Smarandache inverse of y and (x, y) the Smarandache inverse pair with the related pair (a, b). If the pair (a, b) happens to be a Smarandache inverse pair not necessarily with (x, y) as the related pair then we say (a, b) is the Smarandache co inverse pair.*

**Example 6.1.1:** Let $Z'_5 = Z_5 \setminus \{0\}$ (where $Z_5 = \{0, 1, 2, 3, 4\}$) be the group under multiplication modulo 5. Clearly $2 \in Z'_5$ has the Smarandache inverse pair 3 and $2.3 \equiv 1 \pmod 5$ and $4 \in Z'_5$ is such that $2.4 \equiv 3 \pmod 5$ and $3.4 \equiv 2 \pmod 5$ with $4^2 \equiv 1 \pmod 5$. Clearly, $4 \in Z'_5$ has no Smarandache inverse in $Z'_5$ for their exist no $x \in Z'_5 \setminus \{1, 4\}$ such that $x.4 \equiv 4 \pmod 5$. Thus, $4 \in Z'_5$ has no Smarandache inverse. (4, 4) is called the related pair for the Smarandache inverse pair (2, 3).

*Example 6.1.2:* Let $G = \langle g \mid g^6 = 1 \rangle$ be the cyclic group of order 6. Now g in G has $g^5 \in G$ such that $g \cdot g^5 = 1$ further $g^2, g^4 \in G \setminus \{1, g, g^5\}$ is such that $g^5 \cdot g^2 = g$ and $g \cdot g^4 = g^5$ with $g^2 \cdot g^4 = 1$. Clearly $g^3 \cdot g^3 = 1$ but $g^3 \in G$ has no Smarandache inverse for we cannot find $g^i \in G$ with $g^3 \cdot g^i = g^3$. So for $(g, g^5)$ the Smarandache inverse pair the related pair is $(g^2, g^4)$ and $g^3 \in G$ has no Smarandache inverse.

*Remark:* Clearly the pair $(g^2, g^4)$ which is the related pair for the Smarandache inverse pair $(g, g^5)$ is never a Smarandache inverse for the pair $(g^2, g^4)$ has no Smarandache inverse as $g^2 \cdot g^i = g^4$ and $g^j \cdot g^4 = g^2$, $g^i \cdot g^j \in G \setminus \{g^2, g^4, 1\}$ is never possible as i = 2 and j = 4 is the only solution for $g^2 \cdot g^i = g^4$ and $g^j \cdot g^4 = g^2$.

**THEOREM 6.1.1:** *Let G be a group. Every Smarandache inverse has an inverse in the group G but every inverse in G need not have a Smarandache inverse.*

*Proof:* By the very definition of the Smarandache inverse, it should have an inverse. Hence the first part of the theorem. On the other hand, we have every element in G has inverse, but they are not in general Smarandache inverses. For $g^2 \in G$ in example 6.1.2. $g^4$ is the inverse in G as $g^2 \cdot g^4 = 1$ but $g^2$ has no Smarandache inverse. Hence the claim.

*Example 6.1.3:* Let $G = \langle g/g^7 = 1 \rangle$ be the cyclic group of prime order 7. Clearly $g \cdot g^6 = 1$, $g^2 \cdot g^5 = 1$ and $g^3 \cdot g^4 = 1$. For $g \in G$ we have $g^6 \in G$. $g \cdot g^6 = 1$, now $g^2, g^5 \in G \setminus \{1, g, g^6\}$ is such that $g \cdot g^5 = g^6$ and $g^6 \cdot g^2 = g$ with $g^2 \cdot g^5 = 1$. Similarly for $g^2 \in G$, $g^5$ is such that $g^2 \cdot g^5 = 1$. We have $g^3, g^4 \in G \setminus \{1, g^2, g^5\}$ is such that $g^2 \cdot g^3 = g^5$ and $g^4 \cdot g^5 = g^2$ with $g^3 \cdot g^4 = 1$. Now $g^3 \in G$ has a Smarandache inverse $g^4 \in G$ is such that $g^3 \cdot g^4 = 1$. Further g, $g^6 \in G \setminus \{1, g^3, g^4\}$, $g^3 \cdot g = g^4$ and $g^4 \cdot g^6 = g^3$ with $g \cdot g^6 = 1$. Thus it is nice to see every element in G has a Smarandache inverse. This example leads us to the following engaging result about group of prime order; that is for all cyclic groups of prime order p.

**THEOREM 6.1.2:** *Let G be a cyclic group of prime order p; p an odd prime. Every element in $G \setminus \{1\}$ has a Smarandache inverse.*

*Proof:* Given $G = \langle g / g^p = 1 \rangle$ where p is a prime. Now $G \setminus \{1\}$ has exactly p-1 elements which we will pair in the form $(g, g^{p-1})$, $(g^2, g^{p-2})$, $(g^3, g^{p-3})$, ... $\left( g^{\frac{p-1}{2}}, g^{\frac{p+1}{2}} \right)$. The pairs are inverses of each other that is $g \cdot g^{p-1} = 1$, $g^2 \cdot g^{p-2} = 1$, $g^3 \cdot g^{p-3} = 1$, ... , $g^{\frac{p-1}{2}} \cdot g^{\frac{p+1}{2}} = 1$. Now for each pair $(g, g^{p-1})$; $g \cdot g^{p-1} = 1$ the pair $(g^2, g^{p-2})$ acts as a Smarandache inverse, for $(g, g^{p-1})$ as we have $g \cdot g^{p-2} = g^{p-1}$ and $g^{p-1} \cdot g^2 = g$ with $g^2 \cdot g^{p-2} = 1$.

Similarly, we can show for the element $g^2 \in G$, $g^2 \cdot g^{p-2} = 1$ for this pair $(g^3, g^{p-3})$ acts as the Smarandache inverse and so on. Thus finally for the pair $\left( g^{\frac{p-1}{2}}, g^{\frac{p+2}{2}} \right)$ we have $g^{\frac{p-1}{2}} \bullet g^{\frac{p+2}{2}} = 1$; we have the pair $(g, g^{p-1})$ acts as the Smarandache inverse for

$$g^{\frac{p-1}{2}} \bullet g = g^{\frac{p-1+2}{2}} = g^{\frac{p+1}{2}}$$

and

$$g^{\frac{p+1}{2}} \bullet g^{p-1} = g^{\frac{p+1+2p-2}{2}} = g^p \bullet g^{\frac{p-1}{2}} = g^{\frac{p-1}{2}}$$

Hence the claim. This make us to define the following.

**DEFINITION:** *Let G be a group. If every element in G has a Smarandache inverse then we say G is a Smarandache inverse group.*

**THEOREM 6.1.3:** *All symmetric groups $S_n$ of degree n, are not Smarandache inverse groups ($n \geq 4$).*

*Proof:* Given $S_n$ symmetric group of degree n; $n \geq 4$. Clearly $g \in S_n$ is such that $g = \begin{pmatrix} 1 & 2 & 3 & 4 & 5 & \ldots & n \\ 2 & 3 & 4 & 1 & 5 & \ldots & n \end{pmatrix}$. Now g generates a cyclic group of order 4 and $g^2 = \begin{pmatrix} 1 & 2 & 3 & 4 & 5 & \ldots & n \\ 3 & 4 & 1 & 2 & 5 & \ldots & n \end{pmatrix}$ and $g^2$ has no Smarandache inverse in $S_n$ as $g^2 \bullet g^2 = 1 = \begin{pmatrix} 1 & 2 & 3 & \ldots & n \\ 1 & 2 & 3 & \ldots & n \end{pmatrix}$ and we do not have $x \in S_n \setminus \{1, g^2\}$ such that $xg^2 = g^2$. Hence the theorem.

**COROLLARY 6.1.4:** *$S_n$ has elements which have Smarandache inverses.*

*Proof:* By Cauchy's theorem $S_n$ has elements $x (x \neq 1)$ such that $x^p = 1$ where p is an odd prime and $p < n$. By theorem 6.1.2 we have the cyclic group G of order p generated by x is such that $G \subset S_n$ and every element in G has a Smarandache inverse so; G is a Smarandache inverse group.

We can still generalize this to the following theorem.

**THEOREM 6.1.5:** *Let G be a group of finite order. If $p / o(G)$ where p is a prime ($p \geq 5$). Then G has Smarandache inverse elements.*

*Proof:* Given G is a finite group such that p is a prime which divides order of G. Now let $x \in G$ by Cauchy's theorem $x^p = 1$. Then by theorem 6.1.2 every element in the subgroup generated by x has Smarandache inverse baring the identity. Hence the claim.

In the theorem $p \geq 5$ is essential for if $p = 3$ or 2 we see when $p = 2$ no element in the group G has Smarandache inverse. When $p = 3$ we see no element of G has Smarandache inverse.

**THEOREM 6.1.6:** *Let G be a group. If $x \in G$ is such that $x^2 = 1$. Then x has no Smarandache inverse.*

*Proof:* Since $x \in G \setminus \{1\}$ is such that $x^2 = 1$ that is x is a sef inverse element. We see there is no y in $G \setminus \{1, x\}$, suppose we have $y \in G \setminus \{1, x\}$ we will arrive at a contradiction. Given $y \neq 1$ and $y \neq x$ but $xy = x$ multiply by x on the left and use the fact $x^2 = 1$; $x^2y = x^2 = 1$ so $y = 1$ a contradiction to our assumption $y \in G \setminus \{1, x\}$. So if $x \in G$ is such that $x^2 = 1$ then x has no Smarandache inverse.

From this, we arrive at the following conclusions.

**THEOREM 6.1.7:** *The dihedral group $D_{2n}$ has elements, which have no Smarandache inverses.*

*Proof:* We know the dihedral group $D_{2n} = \{a, b | a^2 = b^n = 1 \text{ with } bab = a\}$. We see $a \in D_{2n}$ has no Smarandache inverse as $a^2 = 1$ by theorem 6.1.6. Further every element of the form $ab^i$; $i \leq n-1$ have no Smarandache inverse as $ab^i ab^i = 1$ when $i \leq n-1$. Hence the claim.

**THEOREM 6.1.8:** *The symmetric group $S_3$ has no element, which has Smarandache inverse.*

*Proof:* $S_3 = \{1, p_1, p_2, p_3, p_4, p_5\}$ where

$$p_1 = \begin{pmatrix} 1 & 2 & 3 \\ 1 & 3 & 2 \end{pmatrix}, \quad p_2 = \begin{pmatrix} 1 & 2 & 3 \\ 3 & 2 & 3 \end{pmatrix}, \quad p_3 = \begin{pmatrix} 1 & 2 & 3 \\ 2 & 1 & 3 \end{pmatrix},$$

$$p_4 = \begin{pmatrix} 1 & 2 & 3 \\ 2 & 3 & 1 \end{pmatrix} \text{ and } p_5 = \begin{pmatrix} 1 & 2 & 3 \\ 3 & 1 & 2 \end{pmatrix} \text{ with } 1 = \begin{pmatrix} 1 & 2 & 3 \\ 1 & 2 & 3 \end{pmatrix}.$$

Clearly $p_1^2 = p_2^2 = p_3^2 = 1$, so by theorem 6.1.7 $p_1$, $p_2$ and $p_3$ have no Smarandache inverse. Finally $p_4 \cdot p_5 = 1$ so $p_4$ has no Smarandache inverse for we cannot find a $p_i$, $i \leq 3$ such that $p_u p_i = p_5$. Hence the claim.

**COROLLARY:** *The symmetric group of degree 4 that is $S_4$ has Smarandache inverses.*

*Proof:* Consider the element

$$g = \begin{pmatrix} 1 & 2 & 3 & 4 \\ 2 & 3 & 4 & 1 \end{pmatrix}$$

clearly $g^4 = 1$ and we have $g \cdot g^3 = 1$ and $g^2 \in S_4$ is such that $g^2 \cdot g^3 = g$ and $g^2 \cdot g = g^3$, with $g^2 \cdot g^2 = 1$. Hence the claim. Thus $g \in S_4$ has Smarandache inverse $g^3$ and $(g^2, g^2)$ is the related pair.

**DEFINITION:** *Let G be a group if no element in G has a Smarandache inverse call G a Smarandache inverse free group.*

Clearly, the group $S_3$ is a Smarandache inverse free group. It is clear from the earlier theorems $S_n$ ($n \geq 4$) are not Smarandache inverse free groups.

**DEFINITION:** *Let $x \in G$ have a Smarandache inverse $y$ for the pair $(x, y)$, $(a, b)$ is the related pair. We say $(x, y)$ is a Smarandache self inversed pair if $(a, b)$ has the Smarandache inverse and the related pair is $(x, y)$.*

***Example 6.1.4:*** Let $G = \langle g/g^{25} = 1 \rangle$. Now $g^{15} \in G$ has a Smarandache inverse $g^{10} \in G$ such that $g^{15} \bullet g^{10} = 1$ and $(g^{20}, g^5) \in G$ is such that $g^{15} \bullet g^{20} = g^{10}$, $g^{10} \bullet g^5 = g^{15}$ with $g^{20} \bullet g^5 = 1$.

Now for $g^{20} \in G$ and $g^{20} \bullet g^5 = 1$ we have the pair $(g^{10}, g^{15})$ in $G$ such that $g^{20} \bullet g^{10} = g^5$, $g^5 \bullet g^{15} = g^{20}$, with $g^{20} \bullet g^5 = 1$. Thus we see the pair $(g^{20}, g^5)$ is the Smarandache self inversed pair and its self inverse is $(g^{10}, g^{15})$ and vice versa.

**THEOREM 6.1.9:** *Every Smarandache inverse pair in a group in general need not be a Smarandache self inversed pair.*

*Proof:* We prove this by an example. Let $Z'_5 = Z_5 \setminus \{0\}$ be the group of integers under modulo multiplication 5.

Clearly for $2 \in Z'_5$ we have $3 \in Z'_5$ with $2.3 \equiv 1 \pmod{5}$ and $(4, 4)$ acts as the Smarandache co inverse. But 4 has no Smarandache inverse as $4^2 \equiv 1 \pmod 5$. Hence the claim.

We now define the concept of Smarandache conjugate elements in a group G.

## 6.2 Smarandache Conjugate in Groups

In this section, we introduce the concept of Smarandache conjugate in groups. Throughout this section by a group G we mean only a non commutative group as the concept of conjugates has no meaning in commutative groups. We define a new concept called Smarandache conjugates in a group as follows.

**DEFINITION:** *Let G be a group let $x \in G$; $x$ is said to have a Smarandache conjugate $y$ in G if*

*1. $x$ is conjugate to $y$ (that is there exist $a \in G$ such that $x = aya^{-1}$).*
*2. $a$ is conjugate with $x$ and $a$ is conjugate with $y$.*

Example 6.2.1: Let

$$S_3 = \left\{ 1 = \begin{pmatrix} 1 & 2 & 3 \\ 1 & 2 & 3 \end{pmatrix}, p_1 = \begin{pmatrix} 1 & 2 & 3 \\ 1 & 3 & 2 \end{pmatrix}, p_2 = \begin{pmatrix} 1 & 2 & 3 \\ 3 & 2 & 1 \end{pmatrix}, \right.$$

$$\left. p_3 = \begin{pmatrix} 1 & 2 & 3 \\ 2 & 1 & 3 \end{pmatrix}, p_4 = \begin{pmatrix} 1 & 2 & 3 \\ 2 & 3 & 1 \end{pmatrix} \text{ and } p_5 = \begin{pmatrix} 1 & 2 & 3 \\ 3 & 1 & 2 \end{pmatrix} \right\}$$

be the symmetric group of degree 3. Now $p_1 \in S_3$ has Smarandache conjugate. For we have $p_1$ is conjugate with $p_3$ as $p_1 = p_2 p_3 p_2^{-1}$ that is $p_1 \sim p_3$. Now we have $p_2$ is conjugate with $p_3$ as

$$\begin{pmatrix} 1 & 2 & 3 \\ 3 & 2 & 1 \end{pmatrix} = \begin{pmatrix} 1 & 2 & 3 \\ 1 & 3 & 2 \end{pmatrix} \begin{pmatrix} 1 & 2 & 3 \\ 2 & 1 & 3 \end{pmatrix} \begin{pmatrix} 1 & 2 & 3 \\ 1 & 3 & 2 \end{pmatrix} \text{ that is } p_2 \sim p_3.$$

Also it can be verified $p_1 \sim p_2$ as

$$p_1 = \begin{pmatrix} 1 & 2 & 3 \\ 1 & 3 & 2 \end{pmatrix} = \begin{pmatrix} 1 & 2 & 3 \\ 2 & 1 & 3 \end{pmatrix} \begin{pmatrix} 1 & 2 & 3 \\ 3 & 2 & 1 \end{pmatrix} \begin{pmatrix} 1 & 2 & 3 \\ 2 & 1 & 3 \end{pmatrix}.$$

Thus $p_1$ has Smarandache conjugate $p_3$. Further every element in $S_3$ need not have Smarandache conjugate; for take $p_4 \in S_3$; clearly $p_4$ has no Smarandache conjugate we know $p_4 \sim p_5$ as

$$p_4 = \begin{pmatrix} 1 & 2 & 3 \\ 2 & 1 & 3 \end{pmatrix} = \begin{pmatrix} 1 & 2 & 3 \\ 1 & 3 & 2 \end{pmatrix} \begin{pmatrix} 1 & 2 & 3 \\ 3 & 1 & 2 \end{pmatrix} \begin{pmatrix} 1 & 2 & 3 \\ 1 & 3 & 2 \end{pmatrix}.$$

It is left for the reader to verify that $p_4$ has no Smarandache conjugate but $p_5$ acts only as its conjugate.

This leads to the following theorem.

**THEOREM 6.2.1:** *Let G be a non abelian group. If $x \in G$ has a Smarandache conjugate then x has conjugate; conversely if $x \in G$ has conjugate then $x \in G$ in general have no Smarandache conjugate.*

*Proof:* Clearly by the very definition of Smarandache conjugate we see if x has Smarandache conjugate then it obviously has conjugate.

To prove the fact that if an element has conjugate then it in general need not have Smarandache conjugate. We prove this by the example 6.2.1.

Clearly $p_4 \in S_3$ has conjugate $p_5$ but $p_4$ has no Smarandache conjugate for we see $p_4 = p_1 p_5 p_1$ but $p_1$ can never be conjugate with $p_4$ or $p_5$. Hence the claim.

**THEOREM 6.2.2:** *Let $S_n$ be the symmetric group of degree n, $n \geq 3$. $S_n$ has Smarandache conjugates.*

Proof: Let $x = \begin{pmatrix} 1 & 2 & 3 & 4 & \ldots & n \\ 2 & 1 & 3 & 4 & \ldots & n \end{pmatrix} \in S_n$ we have $y = \begin{pmatrix} 1 & 2 & 3 & 4 & \ldots & n \\ 3 & 2 & 1 & 4 & \ldots & n \end{pmatrix} \in S_n$

is such that x is Smarandache conjugate with y for take $z = \begin{pmatrix} 1 & 2 & 3 & 4 & \ldots & n \\ 1 & 3 & 2 & 4 & \ldots & n \end{pmatrix}$ we have $zyz^{-1} = x$. Also we can easily verify $z \sim x$ and $y \sim z$.

This leads us to find conditions for elements in $S_n$ to be Smarandache conjugate with each other. This is given by the following theorem.

**THEOREM 6.2.3:** *Let $S_n$ be the symmetric group of degree n. Let $x, y \in S_n$ be the Smarandache conjugate by the conjugating element a that is $x = aya^{-1}$. Then all the three elements $a, x, y \in S_n$ have the same cycle decomposition.*

*Proof:* Now to show x is Smarandache conjugate with y, $x, y \in S_n$ we have to show $x = aya^{-1}$ and $a \in S_n$ with $x \sim a$ that is $x = bab^{-1}$ with $b \in S_n$. To do this we need to prove first.

Two permutations in $S_n$ are Smarandache conjugate if and only if x, y and a have the same cycle decomposition. that is if $x = aya^{-1}$ then x is Smarandache conjugate with y if and only if all the three permutations x, y and a in $S_n$ have the same cycle decomposition.

To prove this first we prove the permutations x and y in $S_n$ are conjugate if the permutations x and y have same cycle decomposition. Suppose $x \in S_n$ and that x sends $i \to j$. How do we find $\theta^{-1} x \theta$ where $\theta \in S_n$? Suppose that $\theta$ sends $i \to s$ and $j \to t$ then $\theta^{-1} x \theta$ sends $s \to t$.

In other words, to compute $\theta^{-1} x \theta$ replace every symbol in x by its image under $\theta$. For example to determine $\theta^{-1} x \theta$ where $\theta = (1, 2, 3)(4, 7)$ and $x = (5, 6, 7)(3, 4, 2)$, then since $\theta : 5 \to 5, 6 \to 6, 7 \to 4, 3 \to 1, 4 \to 7, 2 \to 3$, $\theta^{-1} x \theta$ is obtained from x by replacing in x, 5 by 5, 6 by 6, 7 by 4, 3 by 1, 4 by 7 and 2 by 3 so $\theta^{-1} x \theta$ is obtained from x by replacing in x 5 by 5, 6 by 6, 7 by 4, 3 by 1, 4 by 7 and 2 by 3 so that $\theta^{-1} x \theta = (5, 6, 4)(1, 7, 3)$.

With this algorithm for computing conjugates it becomes clear that two permutations having same cycle decomposition are conjugate. For if $x = (a_1, a_2, \ldots, a_{n_1})(b_1, b_2, \ldots, b_{n_2}) \ldots (x_1, x_2, \ldots, x_{n_r})$ and $y = (\alpha_1, \alpha_2, \ldots, \alpha_{n_1})(\beta_1, \beta_2, \ldots, \beta_{n_1}) \ldots (z_1, z_2, \ldots, z_{n_r})$ then $y = \theta^{-1} x \theta$, where one could get as $\theta$ the permutation

$$\begin{pmatrix} a_1 & a_2 & \ldots & a_{n_1} & b_1 & \ldots & b_{n_2} & \ldots & x_1 & \ldots & x_{n_r} \\ \alpha_1 & \alpha_2 & \ldots & \alpha_{n_1} & \beta_1 & \ldots & \beta_{n_2} & \ldots & z_1 & \ldots & z_{n_r} \end{pmatrix}.$$

Thus for instance $(1, 2)(3, 4, 5)(6, 7, 8)$ and $(7, 5)(1, 3, 6)(2, 4, 8)$ can be exhibited as conjugates by using the conjugating permutation

$$\begin{pmatrix} 1 & 2 & 3 & 4 & 5 & 6 & 7 & 8 \\ 7 & 5 & 1 & 3 & 6 & 2 & 4 & 8 \end{pmatrix}.$$

That two conjugates have the same cycle decomposition is now trivial for, by our rule, to compute a conjugate, replace every element in a given cycle by its unique image under the conjugating permutation. Now x ~ y but for x to be Smarandache conjugate we need x = aya$^{-1}$ and x and a are conjugate this once again forces x and a should have the same cycle decomposition. Thus x is Smarandache conjugate with y if and only if all the 3 elements x, y and a have same cycle decomposition. Hence the claim.

Thus we see a pair is Smarandache conjugate if and only if even the conjugating permutation have the same cycle decomposition.

Now we are trying to see whether the concept of Smarandache conjugating is first of all an equivalence relation on the group. Unfortunately even at this stage we see x is Smarandache conjugate with itself for if x ~ x if we have a ∈ G \ {e} with x = axa$^{-1}$ that is ax = xa and a ~ x. Such things happen in reality also for consider $S_4$ the symmetric group of degree 4.

$$x = \begin{pmatrix} 1 & 2 & 3 & 4 \\ 4 & 3 & 2 & 1 \end{pmatrix} \text{ and } a = \begin{pmatrix} 1 & 2 & 3 & 4 \\ 2 & 1 & 4 & 3 \end{pmatrix}$$

xa = ax also;

$$axa^{-1} = \begin{pmatrix} 1 & 2 & 3 & 4 \\ 2 & 1 & 4 & 3 \end{pmatrix}$$

$$\begin{pmatrix} 1 & 2 & 3 & 4 \\ 4 & 3 & 2 & 1 \end{pmatrix} \begin{pmatrix} 1 & 2 & 3 & 4 \\ 2 & 1 & 4 & 3 \end{pmatrix} = \begin{pmatrix} 1 & 2 & 3 & 4 \\ 4 & 3 & 2 & 1 \end{pmatrix}$$

so x ~ x now

$$xax^{-1} = \begin{pmatrix} 1 & 2 & 3 & 4 \\ 4 & 3 & 2 & 1 \end{pmatrix} \begin{pmatrix} 1 & 2 & 3 & 4 \\ 2 & 1 & 4 & 3 \end{pmatrix} \begin{pmatrix} 1 & 2 & 3 & 4 \\ 4 & 3 & 2 & 1 \end{pmatrix} = \begin{pmatrix} 1 & 2 & 3 & 4 \\ 2 & 1 & 4 & 3 \end{pmatrix} = a \neq e = \begin{pmatrix} 1 & 2 & 3 & 4 \\ 1 & 2 & 3 & 4 \end{pmatrix}.$$ Hence the claim.

So x ~ x (~ Smarandache conjugate if there exists a ≠ e with x = axa$^{-1}$ that is ax = xa and implies a = xax$^{-1}$). Unless x commutes with an element other than identity element of G, reflexive property cannot hold good. Clearly if x is Smarandache conjugate with y then obviously y is Smarandache conjugate with x. Thus Smarandache conjugacy is always symmetric. If x ~ y and y ~ z then we have x = aya$^{-1}$ and a = bxb$^{-1}$ and y = czc$^{-1}$, c = dyd$^{-1}$ thus x ~ y ~ a and y ~ z ~ c so to show x Smarandache conjugate with z we have x = tzt$^{-1}$ with x ~ t. Now there are examples in which such result is true so as in the case of reflexive property transitivity many or may not be true.

Thus on the whole we cannot call Smarandache conjugate relation an equivalence relation on a group.

**PROBLEM 1:** Let $S_8$ be the symmetric group of degree 8. Does $(1, 2, 3) (4, 8) (5, 6, 7) \in S_8$ have a Smarandache inverse?

**PROBLEM 2:** Find the Smarandache inverse if it exists for the element

$$\begin{pmatrix} 1 & 1 & 1 \\ 0 & 1 & 1 \\ 0 & 0 & 1 \end{pmatrix} \in S_{3 \times 3}$$

where $S_{3 \times 3}$ is a group of all invariable matrices with entries from $Z_2 = (0, 1)$, under matrix multiplication.

**PROBLEM 3:** Let $G = D_{2.9} = \{a, b / a^2 = b^9 = 1 \text{ bab} = a\}$. Find the Smarandache inverse of $ab^5$?

**PROBLEM 4:** Find all elements in $S_4$, which have Smarandache inverses. Do the collection from a subgroup?

**PROBLEM 5:** Find all the elements in G which have Smarandache inverses where $G = A_5$, the alternating group of degree 5.

**PROBLEM 6:** Find the Smarandache inverse for $g^{15}$ in $G = <g/g^{21} = 1>$.

**PROBLEM 7:** $S_{10}$ be the symmetric group of degree 10. Find the Smarandache inverse for $g = (1, 2, 3, 4) (9, 5) (10, 6, 7)$.

**PROBLEM 8:** Does the Smarandache inverse exist for $(1, 2, 3) (4, 5, 6) (7, 8, 9) (10, 11, 12)$ in $S_{12}$?

**PROBLEM 9:** Does the Smarandache inverse exist for $(1, 2) (3, 4) (5, 6) (7, 8) (9, 0) \in S_{10}$? Justify your answer!

**PROBLEM 10:** Find the Smarandache inverse of $(1, 4, 5, 7) (3, 6, 8, 11) \in S_{11}$.

**PROBLEM 11:** Find the Smarandache conjugate of $(1, 2, 3) (4, 5, 6) (7, 8) \in S_8$.

**PROBLEM 12:** Does there exist a Smarandache conjugate for $(1, 7) (3, 4, 6) (9, 12) \in S_{12}$?

**PROBLEM 13:** Find all elements in $S_4$, which have Smarandache conjugates.

**PROBLEM 14:** Find all elements in $A_6$ which has Smarandache inverses and Smarandache conjugates? Does there exist any relation between these two sets?

**PROBLEM 15:** Find a Smarandache conjugate for $x = (1, 2) (3, 4, 5, 6) \in S_6$. How many elements can be Smarandache conjugate with x?

**PROBLEM 16:** Can $(1, 5, 7) (3, 8, 9) (4, 2)$ and $(8, 2, 7) (9, 3, 4) (1, 5)$ be Smarandache conjugates in $S_9$? Justify your answer.

**PROBLEM 17:** Does $(1, 2, 3, 4) (5) (6, 7, 8) \in S_8$?

1. have Smarandache conjugate in $S_8$.
2. have Smarandache inverse in $S_8$.

**PROBLEM 18:** Can $(1, 2, 3) (4, 5, 6) (7, 8, 9) \in S_9$ have Smarandache inverse?

**PROBLEM 19:** Find the Smarandache inverse of $x = (1, 2, 3, 4) (5, 6) (7, 8) \in S_9$. Is $y = (4, 6, 7, 8) (3, 9) (5, 3)$ a Smarandache conjugate of x?

**PROBLEM 20:** Is $(2, 3) (4, 7) (6, 8)$ is Smarandache conjugate with $(8, 3) (1, 5) (2, 7)$ in $S_8$. Justify your answer.

### 6.3 Smarandache Double Cosets.

In the section we introduce the concept of Smarandache double cosets in a S-semigroup S and prove that Smarandache double coset relation in general is not an equivalence relation on S.

**DEFINITION:** *Let S be a S-semigroup. Let $A \subset S$ and $B \subset S$ be two proper subgroups in S under the operations of S. For some $x \in S$ we define the Smarandache double coset as AxB = {axb | $a \in A$, $b \in B$} AxB is called the Smarandache double coset of A and B for $x \in S$.*

*Remark:* If $x \in A$ or $x \in B$ we do not have any nice special results. Only when we take $x \in S$ and $x \notin A$ and $x \notin B$ we obtain many nice and fascinating properties which does not allow us to extend classically the result, the double coset relation is an equivalence relation on S.

*Example 6.3.1:* Let $Z_{10} = \{0, 1, 2, \ldots, 9\}$ be the S-semigroup under multiplication modulo 10. Take $A = \{1, 9\} \subset Z_{10}$ and $B = \{2, 4, 6, 8\} \subset Z_{10}$. Now take $x = 3$.

A3B = {2, 4, 6, 8} = B. Similarly for x = 7 we get A7B = {2, 4, 6, 8} = B. Consider A5B = {0}. AxB does not divide $Z_{10}$ into disjoint sets.

Now for x = 2, A2A = {2, 8}, A0A ={0}, A5A = {5}, A4A = {4, 6}, A3A = {3, 7}, A1A = {1, 9}. Thus, AxA unlike AxB partition $Z_{10}$ but the number of elements in each class is not the same.

Finally B0B = B5B ={0} BxB = B. Thus this double coset does not even partition $Z_{10}$.

All this study enables us to formulate the following theorem.

**THEOREM 6.3.1:** *Double coset relation on S-semigroup $Z_n$ in general does not partition $Z_n$ for all subgroups in $Z_n$.*

*Proof:* We prove this by an example. Consider the example 6.3.1, the S-semigroup $Z_{10}$. We see AxB gives only two sets viz. {0} and {2, 4, 6, 8} where as BxB = {2, 4, 6, 8} or {0}. Hence the claim.

***Example 6.3.2:*** Consider the semigroup of $2 \times 2$ matrices under matrix multiplication with entries from $Z_2 = \{0, 1\}$.

$$S_{2\times 2} = \left\{ \begin{pmatrix} 0 & 0 \\ 0 & 0 \end{pmatrix}, \begin{pmatrix} 1 & 0 \\ 0 & 0 \end{pmatrix}, \begin{pmatrix} 0 & 1 \\ 0 & 0 \end{pmatrix}, \begin{pmatrix} 0 & 0 \\ 1 & 0 \end{pmatrix}, \begin{pmatrix} 0 & 0 \\ 0 & 1 \end{pmatrix}, \begin{pmatrix} 1 & 1 \\ 0 & 0 \end{pmatrix}, \begin{pmatrix} 1 & 0 \\ 1 & 0 \end{pmatrix}, \begin{pmatrix} 0 & 1 \\ 0 & 1 \end{pmatrix}, \right.$$

$$\left. \begin{pmatrix} 0 & 0 \\ 1 & 1 \end{pmatrix}, \begin{pmatrix} 1 & 0 \\ 0 & 1 \end{pmatrix}, \begin{pmatrix} 0 & 1 \\ 1 & 0 \end{pmatrix}, \begin{pmatrix} 1 & 1 \\ 0 & 1 \end{pmatrix}, \begin{pmatrix} 1 & 0 \\ 1 & 1 \end{pmatrix}, \begin{pmatrix} 1 & 1 \\ 1 & 0 \end{pmatrix}, \begin{pmatrix} 0 & 1 \\ 1 & 1 \end{pmatrix} \text{ and } \begin{pmatrix} 1 & 1 \\ 1 & 1 \end{pmatrix} \right\}.$$

Now $S_{2\times 2}$ is a S-semigroup under multiplication. Consider the subgroups

$$A = \left\{ \begin{pmatrix} 1 & 0 \\ 0 & 1 \end{pmatrix}, \begin{pmatrix} 0 & 1 \\ 1 & 0 \end{pmatrix} \right\} \text{ and } B = \left\{ \begin{pmatrix} 1 & 0 \\ 0 & 1 \end{pmatrix}, \begin{pmatrix} 1 & 0 \\ 1 & 1 \end{pmatrix} \right\}.$$

Calculate AxB by varying $x \in S_{2\times 2}$.

$$AxB = \left\{ \begin{pmatrix} 0 & 0 \\ 0 & 0 \end{pmatrix} \right\} \text{ when } x = \begin{pmatrix} 0 & 0 \\ 0 & 0 \end{pmatrix}$$

$$AxB = \left\{ \begin{pmatrix} 1 & 0 \\ 0 & 1 \end{pmatrix}, \begin{pmatrix} 0 & 1 \\ 1 & 0 \end{pmatrix}, \begin{pmatrix} 1 & 1 \\ 1 & 0 \end{pmatrix}, \begin{pmatrix} 1 & 0 \\ 1 & 1 \end{pmatrix} \right\} \text{ when } x = \begin{pmatrix} 1 & 0 \\ 0 & 1 \end{pmatrix}.$$

$$AxB = \left\{ \begin{pmatrix} 1 & 0 \\ 0 & 0 \end{pmatrix} \begin{pmatrix} 0 & 0 \\ 1 & 0 \end{pmatrix} \right\} \text{ when } x = \begin{pmatrix} 1 & 0 \\ 0 & 0 \end{pmatrix}.$$

$$AxB = \left\{ \begin{pmatrix} 1 & 0 \\ 0 & 0 \end{pmatrix} \begin{pmatrix} 0 & 0 \\ 1 & 0 \end{pmatrix} \right\} \text{ when } x = \begin{pmatrix} 0 & 0 \\ 1 & 0 \end{pmatrix}.$$

$$A \begin{pmatrix} 1 & 1 \\ 0 & 0 \end{pmatrix} B = \left\{ \begin{pmatrix} 1 & 1 \\ 0 & 0 \end{pmatrix}, \begin{pmatrix} 0 & 0 \\ 1 & 1 \end{pmatrix}, \begin{pmatrix} 0 & 1 \\ 0 & 0 \end{pmatrix}, \begin{pmatrix} 0 & 0 \\ 0 & 1 \end{pmatrix} \right\}.$$

$$A \begin{pmatrix} 1 & 0 \\ 1 & 0 \end{pmatrix} B = \left\{ \begin{pmatrix} 1 & 0 \\ 1 & 0 \end{pmatrix} \right\}.$$

$$A \begin{pmatrix} 0 & 1 \\ 0 & 1 \end{pmatrix} B = \left\{ \begin{pmatrix} 0 & 1 \\ 0 & 1 \end{pmatrix}, \begin{pmatrix} 1 & 1 \\ 1 & 1 \end{pmatrix} \right\}.$$

Thus AxB divides the S-semigroup $S_{2\times 2}$ into disjoint set of order 1, 2 and 4. The natural question would be "Is Smarandache double coset relation on a S-semigroup S an equivalence relation on S". The major observation made by us is that Smarandache double cosets in general does not divide S into disjoint sets such that their union is S. Clearly from examples.

To test whether Smarandache double coset relation is an equivalence relation on S, we define a relation called Smarandache double coset relation on S.

**DEFINITION:** *Let S be a S-semigroup. Let $A \subset S$ and $B \subset S$ be subgroups of S under the same operation as in the semigroup S. If for $x, y \in S$ we define $x \underset{s}{\sim} y$ if $y = axb$ for some $a \in A$ and $b \in B$. This relation $\underset{s}{\sim}$ we call as Smarandache double coset relation of S relative to the subgroup A and B.*

Now $\underset{s}{\sim}$ relation may be an equivalence relation on S. Still it may not be a "Smarandache equivalence relation" on S to this end we are going to define now the following.

**DEFINITION:** *Let S be a S-semigroup; the relation $\underset{s}{\sim}$ on S is said to be a Smarandache equivalence relation on S if*

1. *$\underset{s}{\sim}$ is reflexive that is $x \underset{s}{\sim} x$,*
2. *If $x \underset{s}{\sim} y$ then $y \underset{s}{\sim} x$ (symmetric)*
3. *If $x \underset{s}{\sim} y$ and $y \underset{s}{\sim} z$ then $x \underset{s}{\sim} z$ transitivity.*
4. *If $A_1, \ldots, A_n$ are the disjoint sets then we need $S = \bigcup_{i=1}^{n} A_I$.*

*Then only $\underset{s}{\sim}$ is said to be a Smarandache relation on the S-semigroup S.*

The condition 4 is important as we see from example 6.3.1 $\cup A_i \neq S$.

**THEOREM 6.3.2:** *If $\underset{s}{\sim}$ is a Smarandache equivalence relation on S then $\underset{s}{\sim}$ is an equivalence relation on S, but if $\underset{s}{\sim}$ is not an equivalence relation on S then it need not be a Smarandache equivalence relation on S.*

*Proof:* Clearly by the very definition of $\underset{s}{\sim}$ we see if $\underset{s}{\sim}$ is a Smarandache equivalence relation on S then it is an equivalence relation on S.

Now to prove every equivalence relation on the S-semigroup need not in general be a Smarandache equivalence relation on S.

The proof of the second part follows from the example 6.3.1 constructed using the S-semigroup $Z_{10}$. We see if $A = \{2, 4, 6, 8\}$ then $x \underset{s}{\sim} x$ for element in a that is $x = axa$ can never occur for when x is 3 and 5. Hence $\underset{s}{\sim}$ is not a Smarandache relation means it cannot be even an equivalence relation.

Regarding these notions, we suggest the reader to refer the chapter 7 on research problems.

***Example 6.3.3:*** Let S(3) be the set of all mappings of the set (1, 2, 3) to itself. We have o(S(3)) = $3^3$. Now let

$$A = \left\{ \begin{pmatrix} 1 & 2 & 3 \\ 1 & 2 & 3 \end{pmatrix}, \begin{pmatrix} 1 & 2 & 3 \\ 1 & 3 & 2 \end{pmatrix} \right\} \text{ and } B = \left\{ \begin{pmatrix} 1 & 2 & 3 \\ 1 & 2 & 3 \end{pmatrix}, \begin{pmatrix} 1 & 2 & 3 \\ 3 & 2 & 1 \end{pmatrix} \right\}.$$

For $x = \begin{pmatrix} 1 & 2 & 3 \\ 1 & 1 & 1 \end{pmatrix}$ we have $AxB = \left\{ \begin{pmatrix} 1 & 2 & 3 \\ 1 & 1 & 1 \end{pmatrix}, \begin{pmatrix} 1 & 2 & 3 \\ 3 & 3 & 3 \end{pmatrix} \right\}.$

For $x = \begin{pmatrix} 1 & 2 & 3 \\ 2 & 2 & 2 \end{pmatrix}$ we have $AxB = \left\{ \begin{pmatrix} 1 & 2 & 3 \\ 2 & 2 & 2 \end{pmatrix} \right\};$

For $x = \begin{pmatrix} 1 & 2 & 3 \\ 2 & 2 & 1 \end{pmatrix}$ we have

$$AxB = \left\{ \begin{pmatrix} 1 & 2 & 3 \\ 2 & 2 & 1 \end{pmatrix}, \begin{pmatrix} 1 & 2 & 3 \\ 2 & 1 & 2 \end{pmatrix}, \begin{pmatrix} 1 & 2 & 3 \\ 2 & 2 & 3 \end{pmatrix}, \begin{pmatrix} 1 & 2 & 3 \\ 2 & 3 & 2 \end{pmatrix} \right\}$$

For $x = \begin{pmatrix} 1 & 2 & 3 \\ 1 & 2 & 1 \end{pmatrix},$

$$A\begin{pmatrix} 1 & 2 & 3 \\ 1 & 2 & 1 \end{pmatrix}B = \left\{ \begin{pmatrix} 1 & 2 & 3 \\ 1 & 2 & 1 \end{pmatrix}, \begin{pmatrix} 1 & 2 & 3 \\ 1 & 1 & 2 \end{pmatrix}, \begin{pmatrix} 1 & 2 & 3 \\ 3 & 3 & 2 \end{pmatrix}, \begin{pmatrix} 1 & 2 & 3 \\ 3 & 2 & 3 \end{pmatrix} \right\}.$$

For $x = \begin{pmatrix} 1 & 2 & 3 \\ 3 & 2 & 1 \end{pmatrix},$

$$AxB = \left\{ \begin{pmatrix} 1 & 2 & 3 \\ 3 & 2 & 1 \end{pmatrix}, \begin{pmatrix} 1 & 2 & 3 \\ 3 & 1 & 2 \end{pmatrix}, \begin{pmatrix} 1 & 2 & 3 \\ 1 & 2 & 3 \end{pmatrix}, \begin{pmatrix} 1 & 2 & 3 \\ 1 & 3 & 2 \end{pmatrix} \right\}.$$

$$A\begin{pmatrix} 1 & 2 & 3 \\ 2 & 3 & 1 \end{pmatrix}B = \left\{ \begin{pmatrix} 1 & 2 & 3 \\ 2 & 3 & 1 \end{pmatrix}, \begin{pmatrix} 1 & 2 & 3 \\ 2 & 1 & 3 \end{pmatrix} \right\}$$

For $x = \begin{pmatrix} 1 & 2 & 3 \\ 1 & 3 & 3 \end{pmatrix},$ $A\begin{pmatrix} 1 & 2 & 3 \\ 1 & 3 & 3 \end{pmatrix}B = \left\{ \begin{pmatrix} 1 & 2 & 3 \\ 1 & 3 & 3 \end{pmatrix}, \begin{pmatrix} 1 & 2 & 3 \\ 3 & 1 & 1 \end{pmatrix} \right\}.$

$$A\begin{pmatrix} 1 & 2 & 3 \\ 3 & 3 & 1 \end{pmatrix}B = \left\{ \begin{pmatrix} 1 & 2 & 3 \\ 3 & 3 & 1 \end{pmatrix}, \begin{pmatrix} 1 & 2 & 3 \\ 3 & 1 & 3 \end{pmatrix}, \begin{pmatrix} 1 & 2 & 3 \\ 1 & 3 & 1 \end{pmatrix}, \begin{pmatrix} 1 & 2 & 3 \\ 1 & 1 & 3 \end{pmatrix} \right\}.$$

For $x = \begin{pmatrix} 1 & 2 & 3 \\ 2 & 3 & 3 \end{pmatrix},$ $AxB = \left\{ \begin{pmatrix} 1 & 2 & 3 \\ 2 & 3 & 3 \end{pmatrix}, \begin{pmatrix} 1 & 2 & 3 \\ 2 & 1 & 1 \end{pmatrix} \right\}.$

Clearly the Smarandache equivalence relation on S(3) by A and B. Suppose

$$D = \left\{ \begin{pmatrix} 1 & 2 & 3 \\ 1 & 2 & 3 \end{pmatrix}, \begin{pmatrix} 1 & 2 & 3 \\ 2 & 3 & 1 \end{pmatrix}, \begin{pmatrix} 1 & 2 & 3 \\ 3 & 1 & 2 \end{pmatrix} \right\}$$

and

$$C = \left\{ \begin{pmatrix} 1 & 2 & 3 \\ 2 & 1 & 3 \end{pmatrix}, \begin{pmatrix} 1 & 2 & 3 \\ 1 & 2 & 3 \end{pmatrix} \right\}.$$

For $x = \begin{pmatrix} 1 & 2 & 3 \\ 1 & 1 & 1 \end{pmatrix}$, $DxC = \left\{ \begin{pmatrix} 1 & 2 & 3 \\ 1 & 1 & 1 \end{pmatrix}, \begin{pmatrix} 1 & 2 & 3 \\ 2 & 2 & 2 \end{pmatrix} \right\}.$

For $x = \begin{pmatrix} 1 & 2 & 3 \\ 1 & 2 & 3 \end{pmatrix}$ we have

$$DxC = \left\{ \begin{pmatrix} 1 & 2 & 3 \\ 1 & 2 & 3 \end{pmatrix}, \begin{pmatrix} 1 & 2 & 3 \\ 2 & 3 & 1 \end{pmatrix}, \begin{pmatrix} 1 & 2 & 3 \\ 3 & 1 & 2 \end{pmatrix}, \begin{pmatrix} 1 & 2 & 3 \\ 2 & 1 & 3 \end{pmatrix}, \begin{pmatrix} 1 & 2 & 3 \\ 1 & 3 & 2 \end{pmatrix}, \begin{pmatrix} 1 & 2 & 3 \\ 3 & 2 & 1 \end{pmatrix} \right\}.$$

For $x = \begin{pmatrix} 1 & 2 & 3 \\ 1 & 1 & 2 \end{pmatrix}$ we have $D\begin{pmatrix} 1 & 2 & 3 \\ 1 & 1 & 2 \end{pmatrix}C = \left\{ \begin{pmatrix} 1 & 2 & 3 \\ 1 & 1 & 2 \end{pmatrix}, \begin{pmatrix} 1 & 2 & 3 \\ 1 & 2 & 1 \end{pmatrix}, \right.$

$$\left. \begin{pmatrix} 1 & 2 & 3 \\ 2 & 1 & 1 \end{pmatrix}, \begin{pmatrix} 1 & 2 & 3 \\ 2 & 1 & 2 \end{pmatrix}, \begin{pmatrix} 1 & 2 & 3 \\ 1 & 2 & 2 \end{pmatrix}, \begin{pmatrix} 1 & 2 & 3 \\ 2 & 2 & 1 \end{pmatrix} \right\}.$$

$$D\begin{pmatrix} 1 & 2 & 3 \\ 3 & 3 & 1 \end{pmatrix}C = \left\{ \begin{pmatrix} 1 & 2 & 3 \\ 3 & 3 & 1 \end{pmatrix}, \begin{pmatrix} 1 & 2 & 3 \\ 3 & 1 & 3 \end{pmatrix}, \begin{pmatrix} 1 & 2 & 3 \\ 1 & 3 & 3 \end{pmatrix}, \begin{pmatrix} 1 & 2 & 3 \\ 3 & 2 & 3 \end{pmatrix}, \begin{pmatrix} 1 & 2 & 3 \\ 1 & 2 & 2 \end{pmatrix}, \right.$$

$$\left. \begin{pmatrix} 1 & 2 & 3 \\ 3 & 3 & 2 \end{pmatrix} \right\}.$$

$$D\begin{pmatrix} 1 & 2 & 3 \\ 1 & 1 & 3 \end{pmatrix}C = \left\{ \begin{pmatrix} 1 & 2 & 3 \\ 1 & 1 & 3 \end{pmatrix}, \begin{pmatrix} 1 & 2 & 3 \\ 1 & 3 & 1 \end{pmatrix}, \begin{pmatrix} 1 & 2 & 3 \\ 2 & 3 & 2 \end{pmatrix}, \begin{pmatrix} 1 & 2 & 3 \\ 3 & 1 & 1 \end{pmatrix}, \begin{pmatrix} 1 & 2 & 3 \\ 3 & 2 & 2 \end{pmatrix}, \right.$$

$$\left. \begin{pmatrix} 1 & 2 & 3 \\ 2 & 2 & 3 \end{pmatrix} \right\}.$$

$$D\begin{pmatrix} 1 & 2 & 3 \\ 3 & 3 & 3 \end{pmatrix}C = \left\{ \begin{pmatrix} 1 & 2 & 3 \\ 3 & 3 & 3 \end{pmatrix} \right\}.$$

Here also the subgroups of S(3) partition S(3) into Smarandache equivalence classes and it is a Smarandache equivalence relation.

**PROBLEM 1:** Let S(7) be the Smarandache symmetric group. For $A_7$ and $B = \{<(1, 2, 3)(4, 5, 6, 7)>\}$ the group generated by the permutation $(1, 2, 3)(4, 5, 6, 7)$. Is $A_7 \times B$ a Smarandache equivalence relation on S(7)?

**PROBLEM 2:** Let S(4) be the Smarandache symmetric semigroup. Take

$$A = \left\{ \begin{pmatrix} 1 & 2 & 3 & 4 \\ 2 & 3 & 1 & 4 \end{pmatrix}, \begin{pmatrix} 1 & 2 & 3 & 4 \\ 3 & 1 & 2 & 4 \end{pmatrix}, \begin{pmatrix} 1 & 2 & 3 & 4 \\ 1 & 2 & 3 & 4 \end{pmatrix} \right\}$$

and

$$B = \left\{ \begin{pmatrix} 1 & 2 & 3 & 4 \\ 2 & 1 & 4 & 3 \end{pmatrix}, \begin{pmatrix} 1 & 2 & 3 & 4 \\ 1 & 2 & 3 & 4 \end{pmatrix} \right\}$$

as subgroups of S(4). Is $A \times B$ a Smarandache equivalence relation on S(4)?

**PROBLEM 3:** In S(4) if $A = A_4$ and

$$B = \left\{ \left\langle \begin{pmatrix} 1 & 2 & 3 & 4 \\ 2 & 3 & 4 & 1 \end{pmatrix} \right\rangle \right\}$$

that is a group generated by

$$\begin{pmatrix} 1 & 2 & 3 & 4 \\ 2 & 3 & 4 & 1 \end{pmatrix}.$$

Is $A_4 \times B$ a Smarandache equivalence relation on S(4)?

**PROBLEM 4:** Let $S_{3\times 3} = \{a_{ij} \mid a_{ij} \in Z_2 = \{0, 1\}\}$ set of all $3 \times 3$ matrixes, be the semigroup under matrix multiplication. For

$$A = \left\{ \begin{pmatrix} 1 & 0 & 0 \\ 0 & 1 & 0 \\ 0 & 0 & 1 \end{pmatrix}, \begin{pmatrix} 0 & 0 & 1 \\ 0 & 1 & 0 \\ 1 & 0 & 0 \end{pmatrix} \right\}$$

and

$$B = \left\{ \begin{pmatrix} 1 & 0 & 0 \\ 1 & 1 & 0 \\ 1 & 1 & 1 \end{pmatrix}, \begin{pmatrix} 1 & 0 & 0 \\ 0 & 1 & 0 \\ 1 & 0 & 1 \end{pmatrix}, \begin{pmatrix} 1 & 0 & 0 \\ 1 & 1 & 0 \\ 0 & 1 & 1 \end{pmatrix}, \begin{pmatrix} 1 & 0 & 0 \\ 0 & 1 & 0 \\ 0 & 0 & 1 \end{pmatrix} \right\},$$

2 subgroups of $S_{3\times 3}$. Does $A \times B$ divide $S_{3\times 3}$ into Smarandache equivalence classes?

**PROBLEM 5:** Let $Z_{18} = \{0, 1, 2, 3, \ldots, 17\}$ be the S-semigroup under multiplication mod 18. $A = \{1, 17\}$ and $B = \{10, 2, 4, 8, 14, 16\}$ are subgroups of $Z_{18}$. Find the Smarandache double cosets with elements 7, 6 and 5. Is $A \times B$ a Smarandache equivalence relation on $Z_{18}$?

**PROBLEM 6:** Let $Z_{20} = \{0, 1, 2, 3, \ldots, 19\}$ be the S-semigroup under multiplication modulo 20. Can $Z_{20}$ be divided by Smarandache double cosets by any suitable subgroup of $Z_{20}$?

**PROBLEM 7:** Let $Z_{125} = \{0, 1, 2, 3, \ldots, 124\}$ be the S-semigroup under multiplication modulo 125. Let $A = \{1, 124\}$ be a subgroup of $Z_{125}$. Does AxA divide $Z_{125}$ into Smarandache double coset equivalence relations?

**PROBLEM 8:** Let $Z_9 = \{0, 1, 2, \ldots, 8\}$ be the semigroup under multiplication mod 9. $A = \{1, 8\}$ and $B = \{1, 2, 4, 8, 5, 7\}$.

a.  Does the double coset AxB divide $Z_9$ into equivalence classes?

b.  If $A = \{1, 8\}$ and $B = \{1, 8\}$ be subgroups of $Z_9$. Does AxB that is AxA divide $Z_9$ into Smarandache equivalence classes?

c.  Let $A = \{1, 8\}$ and $B = \{1, 7, 4\}$. Does the Smarandache double coset AxB divide $Z_9$ into Smarandache equivalence classes?

**PROBLEM 9:** Let $S_{2\times 2} = \{(a_{ij}) \mid a_{ij} \in Z_5\}$ be the set of all $2 \times 2$ matrixes with entries in $Z_5$, $Z_5$ the prime field of characteristic 5. $S_{2\times 2}$ is a S-semigroup under matrix multiplication. Does there exist subgroups in $S_{2\times 2}$ such that there Smarandache double coset divides $S_{2\times 2}$ into Smarandache equivalence classes?

**PROBLEM 10:** Let $S(4)$ be the Smarandache symmetric group. $A_4$ be the subgroup of $S(4)$. Find $A_4xA_4$ the Smarandache double coset representation of $S(4)$? Does $A_4xA_4$ divide $S(4)$ into Smarandache equivalence classes?

### 6.4 Smarandache Normal subgroups

In this section we introduce the concept of Smarandache normal subgroups to a S-semigroup S and obtain some interesting results about them. This concept leads us to the definition of Smarandache quotient groups.

**DEFINITION:** *Let S be a S-semigroup. Let A be a proper subset of S which is a group under the operation of S. We say A is a Smarandache normal subgroup of the S-semigroup S if $xA \subseteq A$ and $Ax \subseteq A$ or $xA = \{0\}$ and $Ax = \{0\}$ for all $x \in S$ and if 0 is an element in S then we have $xA = \{0\}$ and $Ax = \{0\}$.*

*Remark:* Now we have to define $Ax \subseteq A$, $xA \subset A$ as we have for $x \in S$ we may or may not have $x^{-1}$ to exist in S. That is why we cannot define $xAx^{-1} = A$. Secondly if we restrict ourselves to the subgroup of S then it has nothing to do with S-semigroup.

*Example 6.4.1:* Let $Z_{10} = \{0, 1, 2, \ldots, 9\}$ be a S-semigroup under multiplication modulo 10. Let $A = \{2, 4, 6, 8\} \subset Z_{10}$ be the subgroup of $Z_{10}$. Now 6 is the identity element of A under multiplication modulo 10. Clearly $Ax = xA = A$ for all $x \in Z_{10} \setminus \{0, 5\}$. We have $5A = \{0\}$ and $0A = \{0\}$. Thus A is a Smarandache normal subgroup of the S-semigroup $Z_{10}$.

It is important to note that if the S-semigroup contains 1 as its identity element under multiplication and if the proper subgroup A has the same 1 as the identity element in general A may not be Smarandache normal subgroup of S. We consider the following example.

*Example 6.4.2:* Let $Z_7 = \{0, 1, 2, \ldots, 6\}$ be the S-semigroup $A = \{1, 6\}$. Now A is a subgroup in $Z_7$ but A is not a Smarandache normal subgroup of $Z_7$ as $xA \not\subset A$ for $x \in Z_7 \setminus \{A\} \cup \{0\}$.

So even if the S-semigroup is a commutative group we many not have the subgroups of S to be Smarandache normal subgroups.

This following result allows to state which groups in a certain S-semigroup S are Smarandache normal subgroups of S.

**THEOREM 6.4.1:** *Let $Z_p = \{0, 1, 2, \ldots, p-1\}$ be the S-semigroup of order p under multiplication where p is an odd prime. Then $Z_p$ has only two subgroups $A = \{1, p-1\}$ and $B = \{1, 2, 3, \ldots, p-1\}$ of which A is not a Smarandache normal subgroup of $Z_p$, and B is a Smarandache normal subgroup of $Z_p$.*

*Proof:* Now given $Z_p$ is a S-semigroup of order p, p an odd prime. We have only two subgroups in $Z_p$ viz $A = \{1, p-1\}$ and $B = \{1, 2, \ldots, p-1\}$. Clearly, A is not a Smarandache normal subgroup of $Z_p$ as if we take $0 \neq x \in Z_p \setminus \{1, p-1\}$ we $xA \not\subset A$. Hence the claim.

Now $B = \{1, 2, \ldots, p-1\} \subset Z_p$ and $B = Z_p \setminus \{0\}$. Clearly B is a subgroup and $xB = B$ $x \in Z_p \setminus \{0\}$ and $xB = \{0\}$ if $x = 0$. Hence, B is a Smarandache normal subgroup of $Z_p$. It is still interesting to note that in general when S is a finite S-semigroup then the Smarandache normal subgroup A of S need not divide the order of S.

This is a very different and distinct from the behaviour of groups.

*Example 6.4.3:* Let $S(3)$ be the Smarandache symmetric semigroup.

Let

$$A = \left\{ \begin{pmatrix} 1 & 2 & 3 \\ 1 & 2 & 3 \end{pmatrix}, \begin{pmatrix} 1 & 2 & 3 \\ 2 & 1 & 3 \end{pmatrix} \right\}$$

be a subgroup in $S(3)$. $xA \not\subseteq A$ for $x \in S(3) \setminus S_3$. A is not normal in $S_3$. Hence the claim. So we get a nice theorem about $S(n)$.

**THEOREM 6.4.2:** *Let $S(m)$ be the Smarandache symmetric semigroup. Then $S(m)$ has no subgroup which is Smarandache normal in $S(m)$.*

*Proof:* $S(m)$ is the S-semigroup got by mapping elements of the set $S = \{1, 2, 3, \ldots, m\}$ to itself. Now any proper subset of $S(m)$, which is a subgroup, has only

$$\begin{pmatrix} 1 & 2 & 3 & \ldots & m \\ 1 & 2 & 3 & \ldots & m \end{pmatrix}$$

as the identity element. Let A be proper subset which is a subgroup of S(m). Now for any x ∈ S(m) \ A, xA ⊄ A. Hence S(m) has no subgroups which are Smarandache normal in S(m). Hence the claim.

This leads us to the following definition.

**DEFINITION:** *Let S be a S-semigroup if A has no Smarandache normal subgroup then S is called Smarandache pseudo simple.*

Using the above theorem and definition we can have the following theorem.

**THEOREM 6.4.3:** *S(n) the Smarandache symmetric semigroup is a Smarandache pseudo simple semigroup.*

***Example 6.4.4:*** Let $S_{2\times 2} = \{a_{ij} \mid a_{ij} \in Z_2 = \{0, 1\}\}$ be the set of all $2 \times 2$ matrixes with entries from $Z_2 = \{0, 1\}$. $S_{2\times 2}$ is a S-semigroup under matrix multiplication. Now if

$$A = \left\{ \begin{pmatrix} 0 & 1 \\ 1 & 0 \end{pmatrix}, \begin{pmatrix} 1 & 0 \\ 0 & 1 \end{pmatrix} \right\}$$

is a subgroup of $S_{2\times 2}$. But A is a not Smarandache normal subgroup of $S_{2\times 2}$.

For

$$\begin{pmatrix} 1 & 0 \\ 0 & 0 \end{pmatrix} A \neq A .$$

Hence the claim.

**THEOREM 6.4.4:** $S_{n\times n} = \{a_{ij} \mid a_{ij} \in Z_2 = \{0, 1\}\}$ *set of all $n \times n$ matrices with entries from $Z_2 = \{0, 1\}$ under matrix multiplication is a S-semigroup which is a Smarandache pseudo simple semigroup.*

*Proof:* Now let

$$A = \left\{ \begin{pmatrix} 1 & 0 & \ldots & 0 & 0 \\ 0 & 1 & \ldots & 0 & 0 \\ \vdots & \vdots & \ldots & \vdots & \vdots \\ 0 & 0 & \ldots & 0 & 0 \\ 0 & 0 & \ldots & 0 & 1 \end{pmatrix} = I_{n\times n}, \begin{pmatrix} 0 & 0 & \ldots & 0 & 1 \\ 0 & 0 & \ldots & 1 & 0 \\ \vdots & \vdots & \ldots & \vdots & \vdots \\ 0 & 1 & \ldots & 0 & 0 \\ 1 & 0 & \ldots & 0 & 0 \end{pmatrix} \right\}$$

is a subgroup of $S_{n\times n}$. Clearly A is not a Smarandache normal subgroup of $S_{n\times n}$.

In general $I_{n\times n}$ is the identity element for every subgroup A in $S_{n\times n}$. Now if we take a matrix $B \in S_{n\times n} \setminus \{$set of all invertible n×n matrices with entries from $Z_2 = \{0, 1\}\}$. We get BA ⊄ A. Hence the claim.

Thus we are justified in using the terminology pseudo simple instead of simple as S(n) is pseudo simple for S(n) in the classical sense has the alternating group $A_n \subset S_n \subset S(n)$ to be the normal subgroup. So if we want to define Smarandache quotient group of a S-semigroup S we cannot take S(n) or $S_{n \times n}$; we have to go for only other classes of S-semigroups.

**DEFINITION:** *Let S be a S-semigroup. A be a Smarandache normal subgroup of S. We define the Smarandache quotient group of the S-semigroup S by S/A = {Ax | x ∈ S}.*

**THEOREM 6.4.5:** *Let S be a S-semigroup $A \subset S$ be the Smarandache normal subgroup. The Smarandache quotient group S/A is a semigroup.*

*Proof:* Let S be a S-semigroup. $A \subset S$ be a Smarandache normal subgroup of S. S/A = {Ax | x ∈ S} if X = Ax and Y = Ay we have XY = AxAy we know Ax ⊂ A and Ax ⊂ A. So XY = AxAy ⊂ A. Hence the claim. S/A gives the number of distinct elements of the form Ax. From example 6.4.1 in $Z_{10}$, $\left|\dfrac{Z_{10}}{A}\right| = 2$ viz A and {0}.

**PROBLEM 1:** Let $Z_{23} = \{0, 1, 2, \ldots, 22\}$ be a S-semigroup under multiplication mod 23. Find a Smarandache normal subgroup of $Z_{23}$.

**PROBLEM 2:** Does $Z_{14}$ have Smarandache normal subgroup?

**PROBLEM 3:** Can $Z_8$ have Smarandache normal subgroups?

**PROBLEM 4:** Prove in every $Z_n = \{0, 1, 2, \ldots, n-1\}$ S-semigroup under multiplication mod n, $Z_n$ has a subgroup of order 2 which is never a Smarandache normal subgroup.

**PROBLEM 5:** Let $Z_{12} = \{0, 1, 2, \ldots, 12\}$. Is A = {3, 9} a Smarandache normal subgroup of $Z_{12}$? Can B = {4, 8} be a Smarandache normal subgroup of $Z_{12}$. Is $Z_{12}$ a Smarandache pseudo simple semigroup?

**PROBLEM 6:** Is $Z_8 \times Z_9$ Smarandache pseudo simple semigroup?

**PROBLEM 7:** Find Smarandache normal subgroup of $Z_7 \times Z_3$ where $Z_7 = \{0, 1, 2, \ldots, 6\}$ and $Z_3 = \{0, 1, 2\}$.

**PROBLEM 8:** Find a Smarandache normal subgroup of $Z_8 \times Z_{16} \times Z_6$. Find Smarandache quotient group for any Smarandache normal subgroup.

**PROBLEM 9:** Find Smarandache normal subgroup of $Z_8 \times Z_7 \times Z_{10}$. Find the Smarandache quotient group for

$$\dfrac{Z_8 \times Z_7 \times Z_{10}}{\{\{1\} \times Z_7 \setminus \{0\} \times \{1,9\}\}}.$$

**PROBLEM 10:** Find Smarandache quotient group of the S-semigroup

$$\frac{Z_8 \times Z_9 \times Z_{17}}{\{1 \times 1 \times Z_{17} \setminus \{0\}\}}.$$

## 6.5 Smarandache Direct Product in S-semigroups

In this section we introduce the concept of Smarandache direct product of S-semigroup and obtain some stunning results about them. Now to define the concept of Smarandache direct product in S-semigroup we need the notion of maximal subgroup of a S-semigroup.

**DEFINITION:** *Let S be a S-semigroup we say the proper subset $M \subset S$ is the maximal subgroup of S that is N if a subgroup such that $M \subset N$ then $N = M$ is the only possibility.*

The concept of maximal subgroup in a S-semigroup is such that a S-semigroup can have more than one maximal subgroup.

*Example 6.5.1:* Let $Z_7 = \{0, 1, ... , 6\}$ the S-semigroup under multiplication modulo 7. The only maximal subgroup of $Z_7$ is $G = \{1, 2, 3, ... , 6\} \subset Z_7$.

*Example 6.5.2:* Let $Z_{12} = \{0, 1, 2, 3, ... , 11\}$ be the S-semigroup under multiplication modulo 12. The maximal subgroups of $Z_{12}$ are $A_1 = \{4, 8\}$, $A_2 = \{9, 3\}$ and $A_3 = \{1, 5, 7, 11\}$.

Thus the number of maximal subgroups need not be one that is why we use maximal and not the term "largest". These examples leads to interesting definition.

**DEFINITION:** *Let S be a S-semigroup. If S has only one maximal subgroup we call S a Smarandache maximal semigroup.*

This definition of maximal subgroup of S-semigroup paves way for the following theorems.

**THEOREM 6.5.1:** *The S-semigroup S(n) is a Smarandache maximal semigroup with the maximal group $S_n$.*

*Proof:* The claim is true from the basic fact that $S_n$ contains the set of all 1-1 mapping of the set $S = (1, 2, ... , n)$ onto itself and it is the only maximal subgroup in $S(n)$, as $S(n) \setminus S_n$ has no elements which has inverse. So $S(n)$ is a Smarandache maximal semigroup.

**THEOREM 6.5.2:** *Let $Z_p = \{0, 1, 2, ... , p-1\}$ be the S-semigroup under multiplication mod p where p is a prime. $Z_p$ is a Smarandache maximal semigroup.*

*Proof:* Since p is a prime we know $Z_p \setminus \{0\}$ is a group under multiplication so $G = Z_p \setminus \{0\}$ is the largest subgroup in $Z_p$ under multiplication. Hence the claim.

Making use of the maximal subgroups of the S-semigroup we define the direct product of several S-semigroups as follows:

**DEFINITION:** *Let $S_1, S_2, \ldots, S_n$ be n S-semigroups $S = S_1 \times S_2 \times \ldots \times S_n = \{(s_1, s_2, \ldots, s_n) \mid s_i \in S_n$ for $i = 1, 2, \ldots, n\}$ is called the Smarandache direct product of the S-semigroups $S_1, S_2, \ldots, S_n$ if S is a Smarandache maximal semigroup, and G is got from the $S_1, S_2, \ldots, S_n$ as $G = G_1 \times G_2 \times \ldots \times G_n$ where each $G_i$ is the maximal subgroup of the S-semigroup $S_i$ for $i = 1, 2, \ldots, n$.*

*Example 6.5.3:* Let $S(5)$ be the symmetric S-semigroup. $Z_5 = \{0, 1, 2, 3, 4\}$ be the S-semigroup under multiplication modulo $Z_5$. $S = S(5) \times Z_5$ is the Smarandache direct product for S has a Smarandache maximal semigroup with the largest subgroup G, $G = S_5 \times \{1, 2, 3, 4\}$.

*Example 6.5.4:* $S_1 = S(3)$ and $S_2 = Z_6$ the Smarandache direct product of $S_1 \times S_2 = \{(\sigma, n) / \sigma \in S(3)$ and $n \in Z_6\}$.

Clearly $S_1 \times S_2$ is a S-semigroup for $\{S_3\} \times \{1\}$ is a subgroup of $S_1 \times S_2$ which is group. Hence the claim.

Is $S_1 \times S_2$ the Smarandache direct product of $S_1 \times S_2$? What are the maximal subgroups of $S_1 \times S_2$ are $\{S_3\} \times \{2, 4\}$ and $S_3 \times \{1, 5\}$. Both of them are maximal subgroups of $S(3) \times Z_6 = S_1 \times S_2$. Hence, $S_1 \times S_2$ is not a Smarandache direct product.

**DEFINITION:** *Let S be a S-semigroup. $A_1, \ldots, A_n$ be nonempty subsets of S. We say $A_1, A_2, \ldots, A_n$ is the Smarandache internal direct product of S if $S = A_1, \ldots, A_n = \{a_1 \ldots a_n \mid a_i \in A_i$ $i = 1, 2, \ldots, n\}$ and accounts for all elements in S.*

<u>Remark:</u> We do not demand $A_i$'s to S-semigroups or even semigroup for all situation. So it is sufficient $A_1, A_2, \ldots, A_n$ are just non empty subsets and $S = \{a_1 \ldots a_n \mid a_i \in A_i$ $i = 1, 2, \ldots, n\}$ accounts for all elements of S.

*Example 6.5.5:* Consider $Z_7 = \{0, 1, 2, \ldots, 6\}$. Take $A_1 = \{0, 1\}$ and $A_2 = \{1, 2, \ldots, 6\}$. $A_1 A_2 = \{a_1 a_2 \mid a_1 \in A_1, a_2 \in A_2\} = \{0, 1, 2, 3, \ldots, 6\}$ Here $A_2$ happens to be a group and $A_1$ to be S-semigroup. $Z_7$ is the internal direct product of $A_1 A_2$.

*Example 6.5.6:* Let $Z_6$ be the S-semigroup having elements $\{0, 1, 2, \ldots, 5\}$. $Z_6 = A_1 \bullet A_2 \bullet A_3$ where $A_1 = \{1, 3, 0\}$, $A_2 = \{1, 5\}$ and $A_3 = \{1, 2, 4\}$ be the subsets of $Z_6$. Clearly $Z_6$ is the Smarandache internal direct product of the sets $A_1, A_2$ and $A_3$.

Now we have the following definition.

**DEFINITION:** *Let S be a S-semigroup. If $S = B \bullet A_1 \bullet A_2 \bullet A_3 \bullet A_n$ where B is a S-semigroup and $A_1 \ldots A_n$ are maximal subgroup of S. Then we say S is a Smarandache strong internal direct product.*

*Example 6.5.7:* Let $Z_9 = \{0, 1, 2, \ldots, 8\}$ be the S-semigroup under multiplication mod 9. $Z_9 = A_1 A_2$, where $A_1 = \{0, 1, 3, 8, 6\}$ and $A_2 = \{1, 2, 4, 8, 5, 7\}$ here $A_2$ is the maximal subgroup of $Z_9$ and $A_1$ is a S-semigroup as it contains $\{1, 8\}$ as a subgroup of $A_1$.. Thus, $Z_9$ is the Smarandache strong internal direct product of $A_1$ and $A_2$. Hence the claim.

**THEOREM 6.5.3:** *Let S be a S-semigroup. If S is the strongly Smarandache internal direct product viz. $S = A_1 \times ... \times A_n$ then S is trivially a Smarandache internal direct product and not conversely.*

*Proof:* Follows from the very definition of strongly Smarandache internal direct product and Smarandache internal direct product but converse is not true for S can made as the internal product of subsets $A_i$ for i = 1, 2, ... , n, none of them being subgroups. Hence the claim.

*Example 6.5.8:* Let $Z_{12}$ = {0, 1, 2, 3, ... , 11} be the S-semigroup under multiplication modulo 12. $Z_{12}$ = $A_1 A_2$ where $A_1$ = {0, 1, 2, 3, 4, 6, 8, 9, 10} and $A_2$ = {1, 5, 7, 11}. Here $A_1$ is a S-semigroup and $A_2$ is the maximal subgroup of $Z_{12}$. Thus $Z_{12}$ is Smarandache strong internal direct product.

This S-semigroup $Z_{12}$ has {4, 8} and {3, 9} as its maximal subgroups but in our Smarandache strong internal direct product of $Z_{12}$ we do not take all the maximal subgroup of $Z_{12}$.

Thus from this example one of the important question is whether in the definition of Smarandache strong internal product all maximal subgroups of S will have to be considered in the product. The answer is it cannot be, which is evident from example 6.5.8.

But an observation is important for if the maximal subgroups of a S-semigroup S have different identities then we take only those maximal subgroups which have the identity of S as its identity. This is illustrated by the following example.

*Example 6.5.9:* Let $Z_{20}$ = {0, 1, 2, ... , 19} be the S-semigroup under multiplication modulo 20.

The maximal subgroups of $Z_{20}$ are given by the following tables

| × | 5 | 15 |
|---|---|----|
| 5 | 5 | 15 |
| 15 | 15 | 5 |

*5 is the identity element*

| × | 16 | 4 | 8 | 12 |
|---|----|----|----|----|
| 16 | 16 | 4 | 8 | 12 |
| 4 | 4 | 16 | 12 | 8 |
| 8 | 8 | 12 | 4 | 16 |
| 12 | 12 | 8 | 16 | 4 |

*16 is the identity element*

| × | 1 | 3 | 7 | 9 | 11 | 13 | 17 | 19 |
|---|---|---|---|---|----|----|----|----|
| 1 | 1 | 3 | 7 | 9 | 11 | 13 | 17 | 19 |
| 3 | 3 | 9 | 1 | 7 | 13 | 19 | 11 | 17 |
| 7 | 7 | 1 | 9 | 3 | 17 | 11 | 19 | 13 |
| 9 | 9 | 7 | 3 | 1 | 19 | 17 | 13 | 11 |
| 11 | 11 | 13 | 17 | 19 | 1 | 3 | 7 | 9 |
| 13 | 13 | 19 | 11 | 17 | 3 | 9 | 1 | 7 |
| 17 | 17 | 11 | 19 | 13 | 7 | 1 | 9 | 3 |
| 19 | 19 | 17 | 13 | 11 | 9 | 7 | 3 | 1 |

*1 is the identity element*

Clearly $Z_{20}$ cannot be written as a product of these maximal subgroups. Thus $Z_{20} = A_1A_2$ where $A_1 = (Z_{20} \setminus A_2) \cup \{0\}$ where $A_2 = \{1, 3, 7, 9, 11, 13, 17, 19\}$ the maximal subgroup of $Z_{20}$ and $A_1$ is a S-semigroup as it contains the subgroups $\{5, 15\}$.

Hence $Z_{20}$ is the Smarandache strong internal direct product of $A_1$ and $A_2$.

Thus we see only maximal subgroups which have the same identity as the identity of the S-semigroup S will find its place in the Smarandache strong internal direct product. The following two interesting theorems give an insight into the Smarandache strong internal direct product.

As in the case of groups we may or may not get any proper relation between the Smarandache internal direct product and Smarandache external direct product but one nice relation.

**THEOREM 6.5.4:** *Let S(n) be the S-semigroup S(n) can be represented as the Smarandache strongly internal direct product of S(n).*

*Proof:* S(n) is a S-semigroup. By theorem 6.5.1 $S_n$ is the only maximal Smarandache subgroup of S(n). Hence $S(n) = G \bullet S_n$ where $G = (S(n) \setminus S_n) \cup \{1\}$. Thus $S(n) = G \bullet S_n$ is the Smarandache strong internal direct product of G and $S_n$.

**THEOREM 6.5.5:** *Let S be a S-semigroup. If every maximal subgroup of S contains the same unit as in S, as its identity then we can have all the maximal subgroups in the Smarandache strong internal direct product.*

*Proof:* Let $A_1, \ldots, A_n$ be the collection of all maximal subgroups of S with 1 as their identity for each $A_i$, $i = 1, 2, \ldots, n$. S is a S-semigroup with 1 as its identity. The Set B $= (S \setminus \{A_1 \cup A_2 \ldots \cup A_n\}) \cup \{1\}$.

Clearly $S = BA_1A_2\ldots A_n$ is a strongly Smarandache internal direct product of S.

**THEOREM 6.5.6:** $Z_p = \{0, 1, 2, \ldots, p-1\}$, *p an odd prime be the S-semigroup under multiplication modulo p.* $Z_p$ *has* $A_1 = \{1, p-1\}$ *and* $A_2 = \{1, 2, 3, \ldots, p-1\}$ *as* *subgroups, $Z_p$ is the Smarandache strongly internal direct product of $B \cdot A_2$ where $B$ is any Smarandache subsemigroup containing 0 and 1.*

*Proof:* It can be easily verified using the fact $A_2$ is the only maximal subgroup of $Z_p$ and B any Smarandache subsemigroup of $Z_p$ containing 0 and 1, we have $Z_p = BA_2$ to be the Smarandache strong internal direct product.

**PROBLEM 1:** Find the maximal subgroup of $S = S(9) \times Z_{20} \times Z_8$.

**PROBLEM 2:** Find the Smarandache internal direct product of $Z_{30}$.

**PROBLEM 3:** Find the Smarandache strong internal direct product of $Z_{75}$.

**PROBLEM 4:** Find the Smarandache strong internal direct product of $S(25)$.

**PROBLEM 5:** Find all the maximal subgroup of $S = Z_{12} \times A_5 \times Z_8$.

**PROBLEM 6:** Can $Z_{80}$ be represented as the Smarandache strong internal direct product?

**PROBLEM 7:** Represent $Z_{54}$ as the Smarandache internal direct product.

## Supplementary Reading

1. Herstein, I.N., Topic in Algebra, New York, Blaisdell (1964).

2. Padilla Raul, Smarandache algebraic structures, Bull of Pure and applied Sciences, Delhi, Vol. 17E, No. 1, 119-121, 1998.

3. W.B.Vasantha Kandasamy, Smarandache cosets, Smarandache Notions Journal, American Research Press, 2001. Internet address: http://www.gallup.unm.edu/~smarandache/Cosets.pdf

4. W.B.Vasantha Kandasamy, Smarandache loops, Smarandache Notions Journal, American Research Press, 2001. Internet address: http://www.gallup.unm.edu/~smarandache/Loops.pdf

# RESEARCH PROBLEMS

The study of S-semigroup and the Smarandache notions in groups is a fairly new subject and there are numerous unsolved problems as the very concept of Smarandache algebraic structure is very recent (1998). Some of the problems listed below may be simple but the main motivation for giving these set of research problems is mainly to attract researches and students and make them to contribute to Smarandache algebraic notions.

Any research book if it has research problems it has always a special place among students and researchers. Finally some of the problems are explained with examples.

**PROBLEM 1:** Let $Z_{p^n}$ be the S-semigroup (p a prime, n > 1) under multiplication mod $p^n$. Prove $Z_{p^n}$ has subset of order $p^n - p^r$, $1 < r < n$, which is a subgroup under multiplication mod $p^n$.

The reason for proposing problem 1 is $Z_8$, $Z_9$, and $Z_{25}$ have subsets. $\{1, 3, 5, 7\} \subseteq Z_8$, $\{1, 2, 4, 5, 7, 8\} \subseteq Z_9$ and $\{0, 5, 10, 15, 20\} \subseteq Z_{25}$ which are subgroups of order $p^n - p^r$ for the values of p = 2, 3 and 5 respectively. Problem 1 is a generalization of these examples.

**PROBLEM 2:** Find conditions on n, n a positive non-prime so that $Z_n$ the semigroup under multiplication modulo n is a Smarandache cyclic semigroup i.e. every subset of $Z_n$ which are subgroups of $Z_n$ are cyclic.

**PROBLEM 3:** Let $Z_{p^n}$ be the S-semigroup under multiplication modulo $p^n$, p an odd prime, n an integer greater than 1. Does $Z_{p^n}$ have only 2 proper subsets of order 2 and of order $p^n - p^r$, $(1 \le r \le n)$ which are subgroups of $Z_{p^n}$ under multiplication modulo $p^n$?

**PROBLEM 4:** Give an example of a Smarandache Lagrange semigroup S where S is a non commutative semigroup.

**PROBLEM 5:** Is $Z_{p^n}$ (where p is an odd prime) a S-semigroup which is not even Smarandache weakly Lagrange semigroup?

If problem 3 is true the answer for problem 5 is that $Z_{p^n}$ is not even a Smarandache weakly Lagrange semigroup.

**PROBLEM 6:** Let $S_{n \times n} = \{(a_{ij})/ a_{ij} \in Z_p\}$, the collection of all $n \times n$ matrixes with entries from $Z_p$, p a prime, is a S-semigroup under matrix multiplication. Prove or disprove $S_{n \times n}$ is a Smarandache Lagrange semigroup, when

1) $(p, n) = 1$, n not a prime

2) $(p, n) = p$,

3) $(p, n) = 1$, n is prime.

**PROBLEM 7:** Is $S_{n \times n} = \{(a_{ij})/a_{ij} / Z_m\}$, m a non-prime; a Smarandache Lagrange semigroup? (n – any number, no condition is imposed on it) ($S_{n \times n}$ is as in problem 6).

**PROBLEM 8:** Let $Z_{2^n} = \{0, 1, 2, 3, \ldots, 2^n – 1\}$ be the S-semigroup of order $2^n$ (n > 3) for n arbitrarily, large find the number of Smarandache 2-Sylow subgroups of $Z_{2^n}$. Does $Z_{2^n}$ have subsets of odd order which forms a subgroup under multiplication mod $2^n$? Justify your answer.

**PROBLEM 9:** Does there exists a S-semigroup of order n, $n \geq 24$, in which every subgroup is a Smarandache p-Sylow subgroup?

**PROBLEM 10:** Let S(n) be the symmetric S-semigroup; n an arbitrary integer. Find all subsets in S(n) which form a subgroup or equivalently how many subgroups does S(n) have?

**PROBLEM 11:** Does there exist a non-commutative S-semigroup of order p, p a prime (p > 3) (other than the S-semigroup $Z_p$) such that it has only 2 subsets which are subgroups of which one is of order 2 and the other is of order p – 1?

**PROBLEM 12:** Does there exists a S-semigroup S in which every element is a Smarandache Cauchy element in S? (S should not be taken as $Z_n$, n a composite number for which such result is true).

**PROBLEM 13:** Give an example of a finite S-semigroup S for which every subgroup H of S is such that H partitions S\\{0} into equivalence class of equal cardinality; (S\\{0} only if S contains {0}).

**PROBLEM 14:** Find / characterize all Smarandache inverse free groups.

**PROBLEM 15:** Does there exist Smarandache inverse groups other than cyclic groups of prime order p(p, a prime greater that or equal to 5)?

**PROBLEM 16:** Characterize those groups in which every Smarandache inverse pair is a Smarandache self-inversed pair.

**PROBLEM 17:** Does there exist a non-abelian group of finite order in which every element has a Smarandache conjugate?

**PROBLEM 18:** Characterize groups G in which Smarandache conjugate relation is an equivalence relation on G.

**PROBLEM 19:** Characterize non abelian groups G in which for every x in G the Smarandache conjugate relation is reflexive.

**PROBLEM 20:** Does there exist a group G in which no element in G has a Smarandache conjugate? ($o(G) > 20$).

**PROBLEM 21:** Obtain any interesting relation between Smarandache conjugate elements and Smarandache inverse elements in any group G.

**PROBLEM 22:** Suppose $x \in G$ has no Smarandache inverse does it imply x can be Smarandache conjugate with some element in G? Justify your answer with examples.

**PROBLEM 23:** Let S be a S-semigroup with identity. If every proper subset contained in S, which is a subgroup under the operations of S, contains the same multiplicative identity, which is the identity in the semigroups.

Then for any two subgroup A, B in S, $x \in S$, $AxB = \{ axb \,/\, a \in A \text{ and } b \in B\}$ can we say $x \underset{S}{\sim} y$ i.e. $x = ayb$ implies $\underset{S}{\sim}$ is Smarandache equivalence relation on S (or equivalently) if A and B have different multiplicative identity does it imply $x \underset{S}{\sim} y$ i.e. $x = ayb$ ($a \in A$ and $b \in B$) cannot be a Smarandache equivalence relation on S.

**PROBLEM 24:** Prove or disprove for two distinct groups A and B in S(n), the double coset AxB is not a Smarandache equivalence relation on S(n).

**PROBLEM 25:** Find for what values of n, $Z_n$, n not a prime have Smarandache normal subgroups.

**PROBLEM 26:** Let S be a S-semigroup. Suppose S contains $A_1, A_2, \ldots, A_n$ to be n maximal subgroups of S. Let B be a suitable subset of S which is a S-semigroup. Can we prove $S = BA_1A_2\ldots A_n$ is the Smarandache strong internal direct product in general for any S-semigroup.

**PROBLEM 27:** Give an example of S-semigroup which is simple (other than the class of semigroups given in this book).

**PROBLEM 28:** Give an example of a S-semigroup which has a Smarandache normal subgroup A and S/A is also a S-semigroup. ($S \neq S(n), S \neq Z_p$).

**PROBLEM 29:** Characterize S-semigroup S such that S has one and only one largest subgroup.

For example in case of the Smarandache symmetric semigroup S(n); the largest subgroup of S(n) in $S_n$. When we consider $Z_p$, p a prime the set $A = \{1, 2, \ldots, p-1\}$ is the largest subgroup of $Z_p$.

**PROBLEM 30:** Find the order of the largest subgroup in the S-semigroup $S_{n \times n} = \{(a_{ij}) / a_{ij} \in Z_p\} = \{0, 1, 2, \ldots, p-1\}\}$; p a prime and $(n, p) = 1$; when

  1) $n < p$.
  2) $n > p$.

**PROBLEM 31:** Find the order of the largest subgroup in the S-semigroup $S_{n \times n} = \{(a_{ij}) / a_{ij} \in Z_n = \{0, 1, 2, \ldots, n\}\}$;

  1) when n is a prime
  2) when n is not a prime

**PROBLEM 32:** Find interesting/ innovative results on S-semigroups.

# INDEX

# About the Author

Dr. W. B. Vasantha is an Associate Professor in the Department of Mathematics, Indian Institute of Technology Madras, Chennai, where she lives with her husband Dr. K. Kandasamy and daughters Meena and Kama. Her current interests include Smarandache algebraic structures, fuzzy theory, coding/ communication theory. In the past decade she has guided seven Ph.D. scholars in the different fields of non-associative algebras, algebraic coding theory, transportation theory, fuzzy groups, and applications of fuzzy theory of the problems faced in chemical industries and cement industries. Currently, five Ph.D. scholars are working under her guidance. She has to her credit 241 research papers of which 200 are individually authored. Apart from this she and her students have presented around 262 papers in national and international conferences. She teaches both undergraduate and post-graduate students and has guided over 41 M.Sc. and M.Tech projects. She has worked in collaboration projects with the Indian Space Research Organization and with the Tamil Nadu State AIDS Control Society.

She can be contacted at vasantha@iitm.ac.in
You can visit her on the web at: http://mat.iitm.ac.in/~wbv